\documentclass[10pt]{article}

\usepackage{fancyhdr}
\usepackage{extramarks}
\usepackage{amsmath}
\usepackage{amsthm}
\usepackage{amsfonts}
\usepackage{siunitx}
\usepackage{tikz}
\usepackage[plain]{algorithm}
\usepackage{algpseudocode}
\usepackage{multirow}
\usepackage{booktabs}
\usepackage{palatino}
\usepackage{graphicx}
\usepackage{subfigure}
\usepackage[colorlinks,linkcolor=black,anchorcolor=black,citecolor=black]{hyperref}
\usepackage{amsmath,bm}
\usepackage{booktabs}
\usepackage{mathtools}
\usepackage{amssymb}
\usepackage{caption}
\usepackage{hyperref}
\usepackage{capt-of}
\usepackage{mciteplus}
\usepackage{cite}
\usepackage{mathrsfs}
\usepackage[title,titletoc,toc]{appendix}
\usetikzlibrary{automata,positioning}

\renewcommand{\part}[1]{\textbf{\large Part \Alph{partCounter}}\stepcounter{partCounter}\\}

\DeclareMathOperator{\im}{im}

\newtheorem{thm}{Theorem}[section]
\newtheorem{remark}{Remark}[section]
\newtheorem{Def}{Defination}[section]

\newtheorem{pro}{Property}[section]

\begin{document}

\title{Persistent spectral graph}
\author{  Rui Wang$^1$, Duc Duy Nguyen$^1$ and Guo-Wei Wei$^{1,2,3}$ \footnote{Address correspondences to Guo-Wei Wei. E-mail:wei@math.msu.edu} \\ 
$^1$ Department of Mathematics, 
Michigan State University, MI 48824, USA\\
$^2$  Department of Biochemistry and Molecular Biology\\
Michigan State University, MI 48824, USA \\
$^3$ Department of Electrical and Computer Engineering \\
Michigan State University, MI 48824, USA \\
}
\date{\today}
\maketitle

\begin{abstract}
 Persistent homology is constrained to purely topological persistence  while multiscale graphs account only for geometric information. This work introduces persistent spectral theory to create a unified low-dimensional  multiscale paradigm for revealing topological persistence and extracting geometric shape from high-dimensional datasets.  For a point-cloud dataset, a filtration procedure is used to generate a sequence of chain complexes and associated families of simplicial complexes and chains, from which we construct persistent combinatorial Laplacian matrices.  We show that a full set of topological  persistence can be completely recovered from the harmonic persistent spectra, i.e., the spectra that have zero eigenvalues, of the persistent combinatorial Laplacian matrices. However, non-harmonic spectra of the Laplacian matrices induced by the filtration offer another power tool for data analysis, modeling, and prediction. In this work, non-harmonic persistent spectra are successfully devised to analyze the structure and stability of fullerenes and predict the B-factors of a protein, which cannot be  straightforwardly extracted from the current persistent homology. Extensive numerical experiments indicate the tremendous potential of the proposed  persistent spectral analysis in data science. 
\end{abstract}


\newpage

{\setcounter{tocdepth}{4} \tableofcontents}

\clearpage \pagebreak \setcounter{page}{1}
\renewcommand{\thepage}{{\arabic{page}}}

\section{Introduction}

Graph theory, a branch of discrete mathematics, concerns the relationship between objects. These objects can be either simple vertices, i.e., nodes and/or points (zero simplexes), or high-dimensional simplexes. Here, the relationship refers to connectivity with possible orientations. Graph theory has many branches, such as geometric graph theory, algebraic graph theory, and topological graph theory. The study of graph theory draws on many other areas of mathematics, including algebraic topology, knot theory, algebra, geometry, group theory, combinatorics, etc.  For example, algebraic graph theory can be investigated by using either linear algebra, group theory, or graph invariants. Among them, the use of learning algebra in graph study leads to spectral graph theory.

Precursors of the spectral theory have often had a geometric flavor. An interesting spectral geometry question asked by Mark Kac was ``Can one hear the shape of a drum?'' \cite{kac1966can}. The Laplace-Beltrami operator on a closed Riemannian manifold has been intensively studied \cite{kamber1987rham}. Additionally,  eigenvalues and isoperimetric properties of graphs are the foundation of the explicit constructions of expander graphs \cite{hoory2006expander}. Moreover, the study of random walks and rapidly mixing Markov chains utilized the discrete analog of the Cheeger inequality \cite{chung2005laplacians}. The interaction between spectral theory and differential geometry became one of the critical developments \cite{FanR.K.Chung1997b}. For example, the spectral theory of the Laplacian on a compact Riemannian manifold is a central object of de Rham-Hodge theory \cite{kamber1987rham}. Note that the Hodge Laplacian spectrum contains the topological information of the underlying manifold. Specifically, the harmonic part of the Hodge Laplacian spectrum corresponds to topological cycles.   Connections between topology and spectral graph theory also play a central role in understanding the connectivity properties of graphs \cite{grone1990laplacian,Kirkland2002a, zhang2011laplacian,wu2018weighted}. Similarly, as the topological invariants revealing the connectivity of a topological space, the multiplicity of $0$ eigenvalues of a  0-combinatorial Laplacian matrix is the number of connected components of a graph. Indeed,  the number of $q$-dimensional holes can also be unveiled from the number of $0$ eigenvalues of the $q$-combinatorial Laplacian \cite{serrano2019centrality, hernandez2019higher,Maletic2014, goldberg2002combinatorial}.  Nonetheless,  spectral graph theory offers additional non-harmonic spectral information beyond topological invariants.

The traditional topology and homology are independent of metrics and coordinates and thus, retain little geometric information. This obstacle hinders their practical applicability in data analysis. Recently, persistent homology has been introduced to overcome this difficulty by creating low-dimensional multiscale representations of a given object of interest \cite{frosini1992measuring,edelsbrunner2000topological,zomorodian2005computing,edelsbrunner2008persistent,mischaikow2013morse,carlsson2009zigzag}. Specifically, a filtration parameter is devised to induce a family of geometric shapes for a given initial data. Consequently, the study of the underlying topologies or homology groups of these geometric shapes leads to the so-called topological  persistence.  Like the de Rham-Hodge theory which bridges differential geometry and algebraic topology, persistent homology bridges multiscale analysis and algebraic topology. Topological persistence is the most important aspect of the popular topological data analysis (TDA) \cite{de2007coverage, YaoY:2009,bubenik2014categorification,dey2014computing} and has had tremendous success in computational biology \cite{KLXia:2014c,cang2017topologynet} and worldwide competitions in computer-aided drug design \cite{nguyen2019mathematical}.

Graph theory has been applied in various fields \cite{Garcia-Domenech2008}. For example, spectral graph theory is applied to the quantum calculation of $\pi$-delocalized systems.   The H\"uckel method, or H\"uckel molecular orbital theory, describes the quantum molecular orbitals of $\pi$-electrons in $\pi$-delocalized systems in terms of a kind of adjacency matrix that contains atomic connectivity information \cite{Balasubramanian1985a,gutman1972graph}. Additionally, the Gaussian network model (GNM) \cite{Bahar1997} and anisotropic network model (ANM) \cite{Atilgan2001} represent protein C$_\alpha$ atoms as an elastic mass-and-spring network by graph Laplacians. These approaches were influenced by the Flory theory of elasticity and the Rouse model \cite{Bahar1998}.
Like traditional topology, traditional graph theory extracts very limited information from data. In our earlier work, we have proposed multiscale graphs, called multiscale flexibility rigidity index (mFRI),  to describe the multiscale nature of biomolecular interactions \cite{opron2015communication}, such as hydrogen bonds, electrostatic effects, van der Waals interactions, hydrophilicity, and hydrophobicity. A multiscale spectral graph method has also been proposed as generalized GNM and generalized ANM \cite{xia2015multiscale}. Our essential idea is to create a family of graphs with different characteristic length scales for a given dataset.  We have demonstrated that our multiscale weighted colored   graph (MWCG) significantly outperforms traditional spectral graph methods in protein flexibility analysis \cite{bramer2018multiscale}. More recently, we demonstrate that our MWCG outperforms other existing approaches in protein-ligand binding scoring, ranking, docking, and screening \cite{nguyen2019agl}.

 The objective of the present work is to introduce persistent spectral graph as a new paradigm for the multiscale analysis of the topological invariants and geometric shapes of high-dimensional datasets. Motivated by the success of persistent homology \cite{cang2017topologynet} and   multiscale graphs \cite{nguyen2019agl} in dealing with complex biomolecular data, we construct a family of spectral graphs induced by a filtration parameter. In the present work, we consider the radius filtration via the Vietoris-Rips complex while other filtration methods can be implemented as well. As the filtration radius is increased, a family of persistent $q$-combinatorial Laplacians are constructed for a given point-cloud dataset. The diagonalization of these persistent  $q$-combinatorial Laplacian matrices gives rise to persistent spectra. It is noted that our harmonic persistent spectra of 0-eigenvalues fully recover the persistent barcode or persistent diagram of persistent homology. Additional information is generated from non-harmonic persistent spectra, namely, the non-zero eigenvalues  and associated  eigenvectors. In a combination with a simple machine learning algorithm,  this additional spectral information is found to provide a powerful new tool for the quantitative analysis of molecular data.

\section{Theories and methods}
In this section, we give a brief review of spectral graph theory and simplicial complex to establish
notations and provide essential background. Subsequently, we introduce persistent spectral analysis.

\subsection{Spectral graph theory}

Graph structure encodes inter-dependencies among constituents and provides low-dimensional representations of high-dimensional datasets. One of the representations frequently used in spectral graph theory (SGT) is to associate graphs with matrices, such as the Laplacian matrix and adjacency matrix. Analyzing the spectra from such matrices leads to the understanding of the topological and spectral properties of the graph.

Let $V$ be the vertex set, and $E$ be the edge set. For a given simple graph $G(V, E)$ (A simple graph can be either connected or disconnected), the degree of the vertex $v\in V$ is the number of edges that are adjacent to $v$, denoted $\text{deg}(v)$. The adjacency matrix $\mathcal{A}$ is defined by 
    \begin{equation}
        \mathcal{A}(G)=
        \begin{cases}\label{equ:adj}
            1 & \mbox{if $v_i$ and $v_j$ are adjacent,}   \\
            0 & \mbox{otherwise}.
        \end{cases}
    \end{equation}
    and the Laplacian matrix $   \mathcal{L}$ is given by
    \begin{equation}
        \mathcal{L}(G)=
        \begin{cases}\label{equ:lap}
            \text{deg}(v_i) & \mbox{if $v_i=v_j$,} \\
            -1 & \mbox{if $v_j$ and $v_j$ are adjacent,}   \\
            0 & \mbox{otherwise}.
        \end{cases}
    \end{equation}
Obviously, the adjacency matrix characterizes the graph connectivity.
    The above two matrices are related through diagonal matrix $\mathcal{D}$
    \[
        \mathcal{L} = \mathcal{D} - \mathcal{A}
    \]
   Assuming $G(V,E)$ has $N$ nodes, then adjacency matrix $\mathcal{A}$ and Laplacian matrix $\mathcal{L}$ are both real symmetric $N\times N$ matrices. The eigenvalues of adjacency and Laplacian matrices are denoted and ordered as
    \begin{equation}
        \begin{split}
            \alpha_{\min} = \alpha_{N} \leq \cdots \leq \alpha_{2} \leq \alpha_{1} = \alpha_{\max} \\
            \lambda_{\min} = \lambda_{1} \leq \lambda_{2} \leq \cdots \leq \lambda_{N} = \lambda_{\max}.
        \end{split}
    \end{equation}
    The spectra of $\mathcal{A}$ and $\mathcal{L}$ have several interesting proprieties as following  \cite{zhang2011laplacian}.
    \begin{pro}\label{prop1}
        The eigenvalues of adjacency matrix $\mathcal{A}$ lie in the interval $[-d,d ]$ with $d$ being the largest vertex degree of graph $G(V,E)$.
    \end{pro}
    \begin{pro}\label{prop2}
        A graph is bipartite if and only if its adjacency spectrum is symmetric about 0.
    \end{pro}
    \begin{pro}\label{prop3}
        The Laplacian matrix $\mathcal{L}$ is positive semi-defined, all eigenvalues of $\mathcal{L}$ are non-negative and lie in the interval $[0,2d]$ with $d$  being the largest vertex degree of graph $G(V,E)$.
    \end{pro}
     Except for the aforementioned proprieties for adjacency and Laplacian matrices, one also analyzes the upper and lower bounds for algebraic connectivity $\lambda_2$ and the largest Laplacian eigenvalue $\lambda_{N}$ of graphs. This analysis helps to understand the robustness and connectivity of a graph. For more detailed theorems and proofs, we refer the interested reader to a survey on Laplacian eigenvalues of graphs  \cite{zhang2011laplacian}.
    \begin{thm}\label{thm1}
        Let $G(V,E)$ be a simple graph of order $N$, then the multiplicity of $0$ eigenvalue for Laplacian matrix is the number of connected components of $G(V,E)$. The vertex degree is the value of the diagonal entry.
    \end{thm}
    \begin{thm}\label{thm2}
    Let $G(V, E)$ be a simple graph of order $N$, then the largest eigenvalue
    \begin{equation}
        \lambda_{N} \le \max \{\text{deg}(u)+\text{deg}(v) | (u,v)\in E\},
    \end{equation}
    where $\text{deg}(u)$ is the degree of vertex $u\in V$.
    \end{thm}
    \begin{thm}\label{thm3}
        Let $G(V,E)$ be a simple graph of order $N$ rather than a complete graph with vertex connectivity $\kappa(G)$ and edge connectivity $\kappa^{\prime}(G)$. Then
    \begin{equation}
        2\kappa^{\prime}(G)(1-\cos(\pi/N)) \leq \lambda_{2}(G) \leq \kappa(G) \leq \kappa^{\prime}(G).
    \end{equation}
    The vertex connectivity $\kappa(G)$ is the minimum number of nodes whose deletion disconnects $G$ and edge connectivity $\kappa^{\prime}(G)$ to be the minimum number of edges whose deletion from a graph $G$ disconnects $G$.
    \end{thm}
    \begin{remark}\label{rem1}
        In this section, the Laplacian spectrum of simple graph $G(V,E)$ is concerned. For
        \[
        G(V,E) = G_1(V_1, E_1) \bigcup \cdots \bigcup G_m(V_m,E_m), m\ge 1, m\in \mathbb{Z},
        \]
        where $G_i(V_i,E_i) \subset G(V, E),i=1,\cdots,m$ is a connected simple graph. If $m \ge 2$,  the zero eigenvalue of $L(G)$ has multiplicity $m$, which results in algebraic connectivity $\lambda_2 = 0$. However, $\lambda_2 = 0$ cannot give  any information about the $G_i(V_i, E_i)$. Therefore, we    study the smallest non-zero eigenvalue of $L(G)$, which is actually the smallest algebraic connectivity of $\mathcal{L}(G_i), i=1, \cdots, m$.  In \autoref{sec: PSA} and \autoref{sec: Application}, we  analyze the smallest non-zero eigenvalue of the Laplacian matrix. To make the expression more concise, we still use $\mathbf{\tilde{\lambda}_2}$ as the smallest non-zero eigenvalue. If $G(V,E)$ is a connected simple graph, i.e., $m=1$, one has  $\lambda_2 = \tilde{\lambda}_2$.
    \end{remark}

\begin{figure}[H]
        \centering
        \includegraphics[width=14cm]{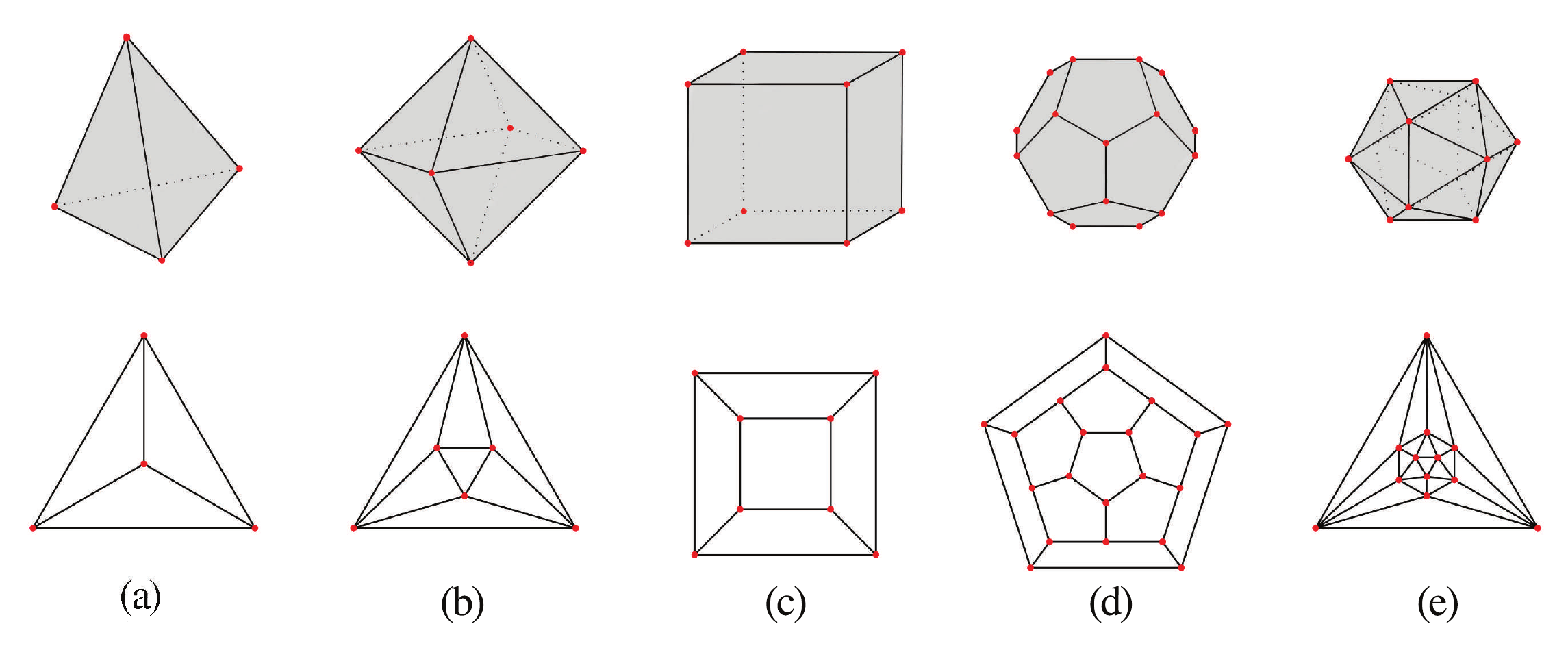}
        \captionsetup{margin=0.8cm}  
        \caption{{Platonic solids (top row) and its platonic graphs (bottom row). (a) Tetrahedron and tetrahedral graph. (b) Octahedron and octahedral graph. (c) Cube and cubical graph. (d) Dodecahedron and dodecahedral graph. (e) Icosahedron and icosahedral graph.}}
        \label{fig:platonic solid}
    \end{figure}

    The charts in the top row of  \autoref{fig:platonic solid} show  $5$ different types of regular convex polyhedrons, which are called platonic solids. The charts in the bottom row are platonic graphs that intuitively describe the vertices and edges as points and line segments in the Euclidean plane. In the three-dimensional (3D) space, objects (vertices) and the relationship (edges) between objects can be expressed by the Laplacian matrix and vice versa. Taking tetrahedron as an example, we denote the top vertex as $v_1$, while the other $3$ vertices on the plane are denoted as $v_2, v_3$, and $v_4$.  The Laplacian matrix of the tetrahedron can be expressed as:
    \[
    \mathcal{L}_{\text{Tetra}} =
    \left[\begin{array}{cccc}
         3 & -1 & -1 & -1 \\
        -1 &  3 & -1 & -1 \\
        -1 & -1 &  3 & -1 \\
        -1 & -1 & -1 &  3
    \end{array}\right]
    \]
    with eigenvalues being $\lambda_1 = 0, \lambda_2 = 4, \lambda_3 = 4,$ and $\lambda_4 = 4$. Topological information can be extracted from this Laplacian matrix. First, the multiplicity of $0$ eigenvalue is $1$, which means there is only  one connected component. Secondly, all of the entries in $L_{\text{Tetra}}$ are non-zero and all the entries in the diagonal equal to $3$, which means the corresponding graph is a complete $3$-regular graph. Moreover, as stated by Theorems \ref{thm2} and \ref{thm3}, the largest eigenvalue $\lambda_4 \le \max \{\text{deg}(u)+\text{deg}(v) | (u,v)\in E\}=6$. The second smallest eigenvalue $\lambda_2$ is controlled by $\kappa(G)=4$ and greater than $2\kappa^{\prime}(G)(1-\cos(\pi/N)) = 2\times 3 \times (1-\cos(\pi/                                                         4)) \approx 1.75735931288$. A similar analysis can be applied to other $4$ platonic solids. \autoref{table:platonic solid} shows some characteristics of platonic solids.

    \begin{table}[H]
        \centering
        \setlength\tabcolsep{20pt}
        \captionsetup{margin=0.8cm}
        \caption{Characteristics of platonic solids in \autoref{fig:platonic solid}. $V, E, \beta_0, \kappa, \kappa^{\prime}, d,$ and $\tilde{\lambda}_2$ stand for the number of vertices, the number of edges, the number of zero eigenvalues, vertex connectivity, edge connectivity, and the smallest non-zero eigenvalue, respectively. Here, all the platonic graphs are connected simple graphs, so $\lambda_2 = \tilde{\lambda}_2$.}
        \begin{tabular}{cccccccc}
        \hline
         Platonic Solid     & $V$  & $E$   & $\kappa$    & $\kappa^{\prime}$          & $\beta_0$  & $\tilde{\lambda}_2$\\ \hline
         Tetrahedron        & $4$  & $6$   & $4$         & $3$                        & $1$    & $4.00$\\
         Octahedron         & $6$  & $12$  & $4$         & $4$                        & $1$    & $4.00$\\
         Cube               & $8$  & $12$  & $3$         & $3$                        & $1$    & $2.00$\\
         Dodecahedron       & $20$ & $30$  & $3$         & $3$                        & $1$    & $0.7639$\\
         Icosahedron        & $12$ & $30$  & $5$         & $5$                        & $1$    & $2.76$\\ \hline
        \end{tabular}
        \label{table:platonic solid}
    \end{table}

    \subsection{Simplicial complex}
    A simplicial complex is a powerful algebraic topology tool that has wide applications in graph theory, topological data analysis  \cite{edelsbrunner2008persistent}, and many physical fields \cite{cang2017topologynet}. We briefly review simplicial complexes to generate notation and provide  essential preparation  for  introducing persistent spectral graph.

  \subsubsection{Simplex}
        Let $\{v_0, v_1, \cdots, v_q\}$ be a set of points in $\mathbb{R}^n$. A point $v = \displaystyle{\sum_{i=0}^{q}}\lambda_iv_i, \lambda_i \in \mathbb{R}$ is an affine combination of  $v_i$ if  $\displaystyle{\sum_{i=0}^{q}}\lambda_i = 1$.  An affine hull is the set of affine combinations. Here, $q+1$ points $v_0, v_1,\cdots, v_q$ are affinely independent if $v_1-v_0, v_2-v_0, \cdots, v_q-v_0$ are linearly independent. A $q$-plane is well-defined if the $q+1$ points are affinely independent. In $\mathbb{R}^n$, one can have at most $n$ linearly independent vectors. Therefore, there are  at most $n+1$ affinely independent points. An affine combination $v = \displaystyle{\sum_{i=0}^{q}}\lambda_iv_i$ is a convex combination if all $\lambda_i$ are non-negative. The convex hull is the set of convex combinations.

\begin{figure}[H]
            \centering
            \includegraphics[width=13cm]{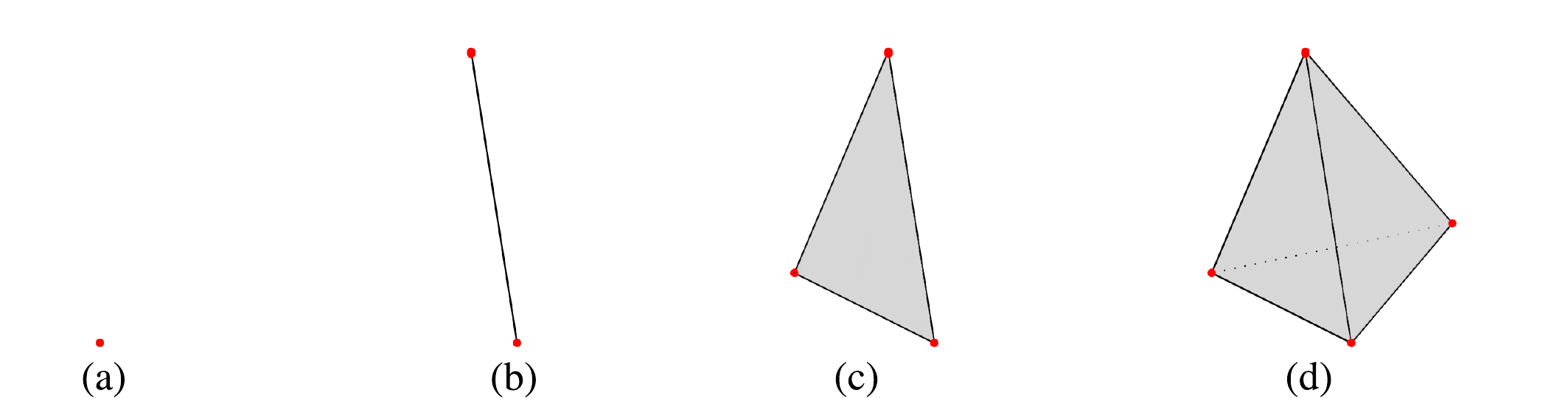}
            \captionsetup{margin=0.9cm}  
            \caption{Illustration of simplices. (a) $0$-simplex (a vertex), (b) $1$-simplex (an edge),   $(c)$ $2$-simplex (a triangle), and  (d) $3$-simplex (a tetrahedron)}
            \label{fig:simplex}
        \end{figure}

        A (geometric) $q$-simplex denoted as $\sigma_q$ is the convex hull of $q+1$ affinely independent points in $\mathbb{R}^n (n \ge q)$
        with dimension $\text{dim}(\sigma_q) = q$. A $0$-simplex is a vertex, a $1$-simplex is an edge, a $2$-simplex is a triangle, and a $3$-simplex is a tetrahedron, as shown in \autoref{fig:simplex}. The convex hull of each nonempty subset of $q+1$ points forms a subsimplex and is regraded as a face of $\sigma_q$ denoted $\tau$.  The $p$-face of a $q$-simplex is the subset $\{v_{i1}, \cdots, v_{ip}\}$ of the $q$-simplex.

        \subsubsection{Simplicial complex}
\begin{figure}[H]
            \centering
            \includegraphics[width=15cm]{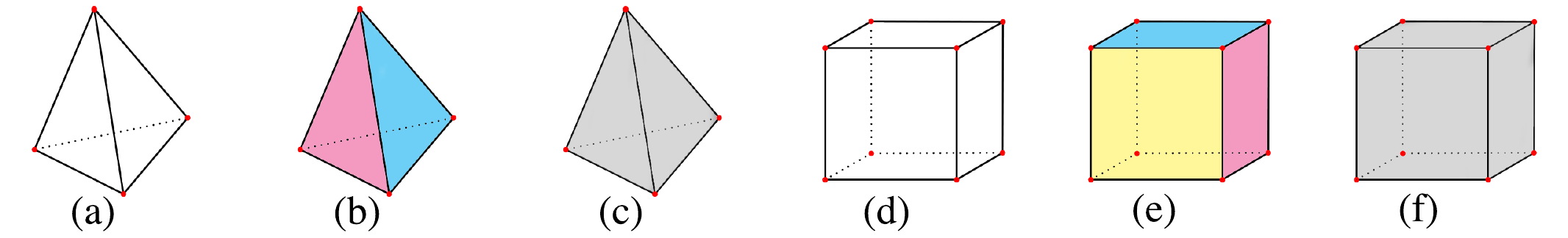}
            \captionsetup{margin=0.9cm}  
            \caption{Illustrations of simplicial complexes}
            \label{fig:simplical complex}
        \end{figure}

        A (finite) simplicial complex $K$  is a (finite) collection of simplices   in $\mathbb{R}^n$ satisfying the following conditions
        \begin{itemize}
            \item[(1)] If $\sigma_q \in K$ and $\sigma_p$ is a face of $\sigma_q$, then $\sigma_p \in K$.
            \item[(2)] The non-empty intersection of any two simplices $\sigma_q, \sigma_p \in K$ is a face of both  of $\sigma_q $ and $ \sigma_p$.
        \end{itemize}
        Each element $\sigma_q \in K$ is a $q$-simplex. The dimension of $K$ is defined as $\text{dim}(K) = \max\{\text{dim}(\sigma_q): \sigma_q \in K\}$.
				To distinguish topological spaces based on the connectivity of simplicial complexes, one uses  Betti numbers. The $k$-th Betti number, $\beta_k$, counts the number of $k$-dimensional holes on a topological surface. The geometric meaning of Betti numbers in $\mathbb{R}^3$ is the following: $\beta_0$ represents the number of connected components, $\beta_1$ counts the number of  one-dimensional loops or circles, and $\beta_2$ describes the number of two-dimensional voids or holes. In a nutshell, the Betti number sequence $\{\beta_0, \beta_1, \beta_2, \cdots\}$ reveals the intrinsic topological property of the system. To illustrate the simplicial complex and its corresponding Betti number, we have designed two simple models as is shown in \autoref{fig:simplical complex}. \footnote{These examples show in an intuitive way to count Betti numbers. However, In \autoref{sec: PSA}, it is impossible to generate structures (b), (e), and (f).}
        \\

        \begin{table}[H]
            \centering
            \setlength\tabcolsep{8pt}
            \captionsetup{margin=0.9cm}
            \caption{The Betti number of simplicial complexes in \autoref{fig:simplical complex}. Each color represents different faces.  The tetrahedron-shaped simplicial complexes are demonstrated in  (a)-(c), and the cube-shaped simplicial complexes are depicted in  (d) - (f).  (a)  and (d) only has $0$-simplices and $1$-simplices,  (b) has four $2$-simplices, and (c) has one more $3$-simplex.  (e) and (f) do not have any $2$-simplex.  }
            \begin{tabular}{ccccccc}
            \hline
             Betti number    & Fig. $3$ (a)    & Fig. $3$ (b)  & Fig. $3$ (c)  & Fig. $3$ (d)  & Fig. $3$ (e)  & Fig. $3$ (f) \\ \hline
             $\beta_0$       & $1$    & $1$    & $1$    & $1$    & $1$    & $1$  \\
             $\beta_1$       & $3$    & $0$    & $0$    & $5$    & $0$    & $0$  \\
             $\beta_2$       & $0$    & $1$    & $0$    & $0$    & $1$    & $0$  \\ \hline
            \end{tabular}
            \label{table:betti}
        \end{table}

        Recall that in graph theory, the degree of a vertex ($0$-simplex) $v$ is the number of edges that are adjacent to the vertex, denoted as deg$(v)$. However, once we generalize this notion to $q$-simplex, problem aroused since $q$-simplex can have $(q-1)$-simplices and  $(q+1)$-simplices adjacent to it at the same time. Therefore, the upper adjacency  and lower adjacency are required to define the degree of a $q$-simplex for $q>0$  \cite{serrano2019centrality,Maletic2014}.
        \begin{Def}
            Two $q$-simplices $\sigma^i_{q}$ and $\sigma^j_{q}$ of a simplicial complex $K$ are lower adjacent if they share a common $(q-1)$-face, denoted $\sigma^i_q \stackrel{L}\sim \sigma^j_q$. The lower degree of $q$-simplex, denoted deg$_L(\sigma_q)$, is the number of nonempty $(q-1)$-simplices in $K$ that are faces of $\sigma_q$, which is always $q+1$.
        \end{Def}
        \begin{Def}
            Two $q$-simplices $\sigma^i_{q}$ and $\sigma^j_{q}$ of a simplicial complex $K$ are upper adjacent if they share a common $(q+1)$-face, denoted $\sigma^i_q \stackrel{U}\sim \sigma^j_q$. The upper degree of $q$-simplex, denoted deg$_U(\sigma_q)$, is the number of $(q+1)$-simplices in $K$ of which $\sigma_q$ is a face.
        \end{Def}
        Then, the degree of a $q$-simplex ($q>0$) is defined as:
        \begin{equation}
            \text{deg}(\sigma_q) = \text{deg}_L(\sigma_q) + \text{deg}_U(\sigma_q) = \text{deg}_U(\sigma_q) + q + 1.
        \end{equation}

      \subsubsection{Chain complex}
        Chain complex is an important concept in topology, geometry, and algebra.
				Let $K$ be a simplicial complex of dimension $q$. A  $q$-chain is a formal sum of $q$-simplices in $K$ with $\mathbb{Z}_2$ field of the coefficients for the sum. A $q$-chain is called $q$-cycle if its boundary is zero. Under the addition operation of $\mathbb{Z}_2$, a set of all $q$-chains is called a chain group and denoted $C_q(K)$. To relate these chain groups, we denote boundary operator by $\partial_q: C_q(K) \longrightarrow C_{q-1}(K)$. The boundary operator maps a $q$-chain which is a linear combination of $q$-simplices to the same linear combination of the boundaries of the $q$-simplices. Denoting $\sigma_q = [v_0, v_1,\cdots,v_q]$ for the $q$-simplex spanned by its vertices, its boundary operator can be defined as:
        \begin{equation}
            \partial_q \sigma_q = \sum_{i=0}^{q}(-1)^i\sigma^{i}_{q-1},
        \end{equation}
        with $\sigma_q = [v_0, \cdots, v_q]$ being the $q$-simplex. Here, $\sigma^{i}_{q-1} = [v_0, \cdots, \hat{v_i},\cdots,v_q]$ is the $(q-1)$-simplex with $v_i$ being omitted. A chain complex is the sequence of chain groups connected by boundary operators
        \begin{equation}
            \cdots \stackrel{\partial_{q+2}}\longrightarrow C_{q+1}(K) \stackrel{\partial_{q+1}}\longrightarrow C_{q}(K) \stackrel{\partial_{q}}\longrightarrow C_{q-1}(K)\stackrel{\partial_{q-1}} \longrightarrow \cdots
        \end{equation}

\subsection{Persistent spectral analysis}\label{sec: PSA}

In this section, we introduce persistent spectral theory (PST) to extract rich topological and spectral information of simplicial complexes via a  filtration process. We briefly review preliminary concepts about the oriented simplicial complex and $q$-combinatorial Laplacian, while more detail information can be found  elsewhere \cite{hernandez2019higher, Maletic2014, goldberg2002combinatorial, horak2013spectra}. Then, we  discuss the properties of the $q$-combinatorial Laplacian matrix together with its spectrum. Moreover, we employ the $q$-combinatorial Laplacian to establish the PST. Finally, we discuss some variants of the persistent $q$-combinatorial Laplacian matrix and illustrate their formulation on simple geometry, i.e., a benzene molecule.

  \subsubsection{Oriented simplicial complex and $q$-combinatorial Laplacian}

    An oriented simplicial complex is the one in which all of the simplices in the simplicial complex, except for vertices and $\emptyset$, are oriented. A $q$-combinatorial Laplacian is defined based on oriented simplicial complexes, and its lower- and higher-dimensional simplexes can be employed to study a specifically oriented simplicial complex.

    We first introduce oriented simplex complexes. Let $\sigma_q$ be a $q$-simplex, we can define the ordering of its vertex set. If two orderings defined on  $\sigma_q$  differ from each other by an even permutation, we say that they are equivalent, and each of them  is called an orientation of $\sigma_q$. An oriented $q$-simplex is a simplex $\sigma_q$ with the orientation of $\sigma_q$. An oriented simplicial complex $K$ is defined if all of its simplices are oriented. Suppose $^i\sigma_q$ and $^j\sigma_q \in K$ with $K$ being an oriented simplicial complex. If $^i\sigma_q$ and $^j\sigma_q$ are upper adjacent with a common upper $(q+1)$-simplex $\tau_{q+1}$, we say they are similarly oriented if both have the same sign in  $\partial_{q+1}(\tau_{q+1})$ and dissimilarly oriented if the signs are opposite. Additionally, if $^i\sigma_q$ and $^j\sigma_q$ are lower adjacent with a common lower $(q-1)$-simplex $\eta_{q-1}$, we say they are similarly oriented if $\eta_{q-1}$ has the same sign in  $\partial_q(^i\sigma_q)$ and $\partial_q(^j\sigma_q)$, and dissimilarly oriented if the signs are opposite.

Similarly, we can define $q$-chains based on an oriented simplicial complex $K$. The $q$-chain $C_q(K)$ is also defined as the linear combinations of the basis, with the basis being the set of oriented $q$-simplices of $K$. The $q$-boundary operator $\partial_q: C_q(K) \longrightarrow C_{q-1}(K)$ is
    \begin{equation}
        \partial_q \sigma_q = \sum_{i=0}^{q}(-1)^i \sigma^{i}_{q-1},
    \end{equation}
    with $\sigma_q = [v_0, \cdots, v_q]$ to be the oriented $q$-simplex, and $\sigma^{i}_{q-1} = [v_0, \cdots, \hat{v_i} ,\cdots,v_q]$ the oriented $(q-1)$-simplex with its vertex $v_i$ being removed. Let $\mathcal{B}_q$
    be the matrix representation of a $q$-boundary operator with respect to the standard basis for $C_q(K)$ and $C_{q-1}(K)$ with some assigned orderings. Then, the number of rows in $\mathcal{B}_q$ corresponds to the number of $(q-1)$-simplices and the number of columns shows the number of $q$-simplices in $K$, respectively. Associated with the $q$-boundary operator is the adjoint operator denoted $q$-adjoint boundary operator,  defined as
    \begin{equation}
        \partial_q^{\ast}: C_{q-1}(K) \longrightarrow C_q(K),
    \end{equation}
    and the transpose of $\mathcal{B}_q$, denoted $\mathcal{B}_q^T$, is the matrix representation of $\partial_q^{\ast}$ relative to the same ordered orthonormal basis as $\partial_q$  \cite{spence2000elementary}.

    \begin{figure}[H]
        \centering
        \includegraphics[width=12cm]{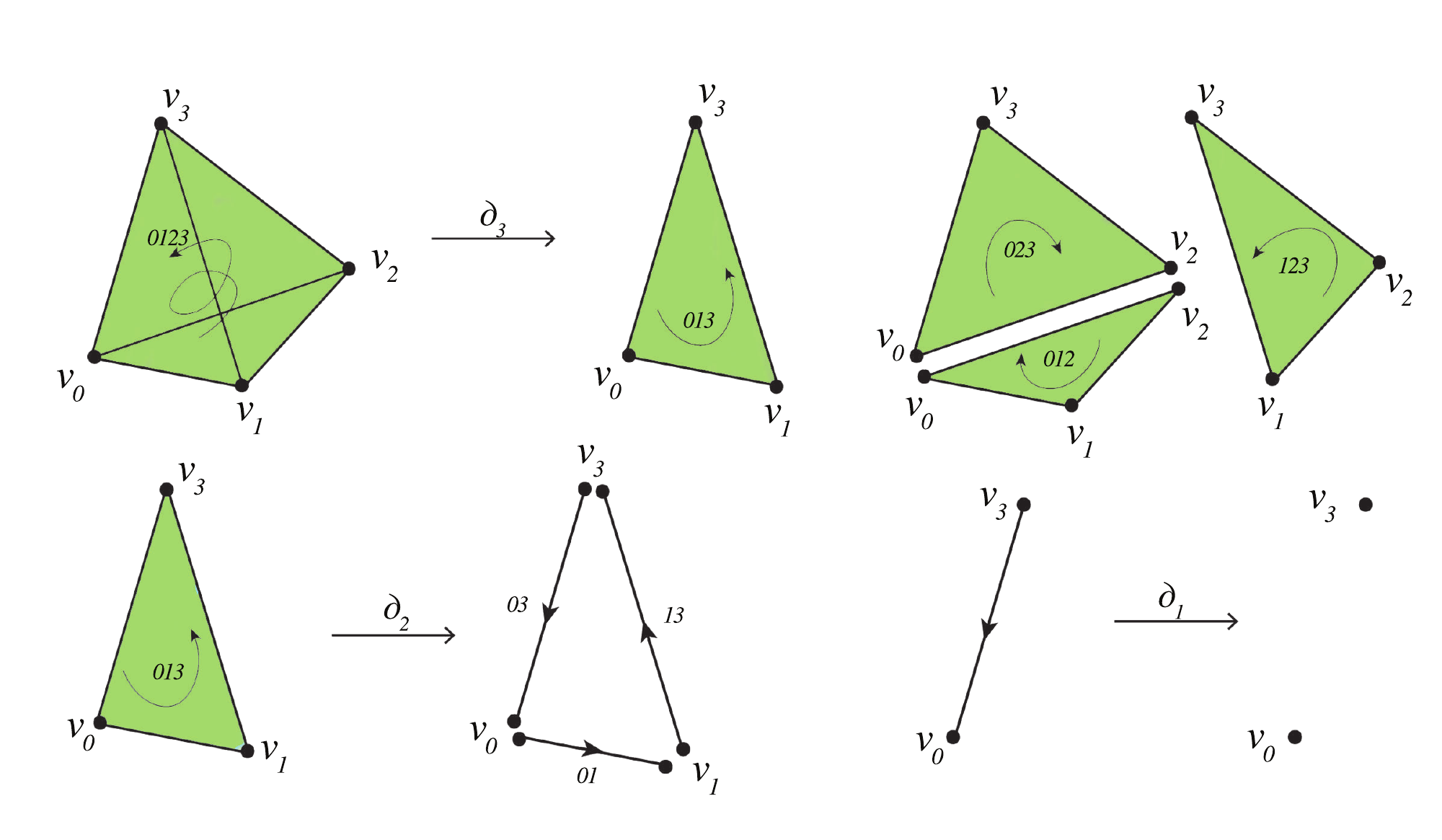}
        \captionsetup{margin=0.9cm}  
        \caption{Examples of $1,2,3$-boundary operators. For sake of brevity, we use $0123$ to represent a $3$-simplex. $013,023,012,$ and $123$ to represent $2$-simplices, and $01,13,$ and $03$ to represent $1$-simplices. The $3$-simplex $0123$ has a right-handed orientation. After the boundary map $\partial_3$, we have $\partial_3(0123) = 123 - 023 + 013 - 012$. Two $2$-simplices, $023$ and $012$, have opposite orientations as $123$ and $013$. Similarly, $\partial_2(013) = 13 - 03 + 01$, and $\partial_1(03) = 3 - 0$.}
        \label{fig:oriented}
    \end{figure}

    Let $K$ be an oriented simplicial complex, for integer $q\ge 0$, the $q$-combinatorial Laplacian is a linear operator $\Delta_q: C_q(K) \longrightarrow C_q(K)$
    \begin{equation}\label{equ:laplacian operator}
        \Delta_q := \partial_{q+1} \partial_{q+1}^{\ast} + \partial_{q}^{\ast} \partial_{q}
    \end{equation}
    with $\partial_q \partial_{q+1} = 0$, which implies $\text{Im}(\partial_{q+1}) \subset \text{ker}(\partial_q)$.
    The $q$-combinatorial Laplacian matrix, denoted $\mathcal{L}_q$, is the matrix representation\footnote{If $q=0$, $\partial_0$ is a zero map, and we denote $\mathcal{B}_0$ a zero matrix with dimension $1\times N$, where $N$ is the number of $0$-simplices. Note that this term is needed to attain the correct dimension of the null space. If $q > \text{dim}(K)$, we will not discuss   $\partial_q$ since there is no $\sigma_q$ in $K$.}.
    \begin{equation}\label{equ:combinatorial Laplacian}
        \mathcal{L}_q = \mathcal{B}_{q+1}\mathcal{B}_{q+1}^{T} + \mathcal{B}_q^T \mathcal{B}_q
    \end{equation}
    of operator $\Delta_q$, with $\mathcal{B}_q$ and $\mathcal{B}_{q+1}$ being the matrices of dimension $q$ and $q+1$. Additionally, we denote upper and lower $q$-combinatorial Laplacian matrices by $\mathcal{L}_q^U = \mathcal{B}_{q+1}\mathcal{B}_{q+1}^{T}$ and $\mathcal{L}_q^L = \mathcal{B}_q^T \mathcal{B}_q$, respectively. Note that $\partial_0$ is the zero map which leads to $\mathcal{B}_0$ being a zero matrix. Therefore, $\mathcal{L}_0(K) = \mathcal{B}_{1} \mathcal{B}_1^T + \mathcal{B}_0^T \mathcal{B}_0$, with $K$ the (oriented) simplicial complex of dimension $1$, which is actually a simple graph. Especially, $0$-combinatorial Laplacian matrix $\mathcal{L}_0(K)$ is actually the Laplacian matrix defined in the spectral graph theory. In fact, Eq. \eqref{equ:lap} is exactly the same as Eq. \eqref{equ:L0} given below.

    Given an oriented simplicial complex $K$ with $0 \le q \le \text{dim}(K)$, one can obtain the entries of its corresponding upper and lower $q$-combinatorial Laplacian matrices explicitly \cite{goldberg2002combinatorial}
    \begin{align}
        (\mathcal{L}_q^U)_{ij} &=
        \begin{cases}\label{equ:upper}
            \text{deg}_U(\sigma^i_{q}), & \mbox{if $i=j$,} \\
            1,                          & \mbox{if $i\neq j$, $\sigma^i_q \stackrel{U}\sim \sigma^j_q$ with similar orientation,}   \\
            -1,                         & \mbox{if $i\neq j$, $\sigma^i_q \stackrel{U}\sim \sigma^j_q$ with dissimilar orientation,} \\
            0,                          & \mbox{otherwise.}
        \end{cases} \\
        (\mathcal{L}_q^L)_{ij} &=
        \begin{cases}\label{equ:lower}
            \text{deg}_L(\sigma^i_{q})=q+1, & \mbox{if $i=j$,} \\
            1,                              & \mbox{if $i\neq j$, $\sigma^i_q \stackrel{L}\sim \sigma^j_q$ with similar orientation,}   \\
            -1,                             & \mbox{if $i\neq j$, $\sigma^i_q \stackrel{L}\sim \sigma^j_q$ with dissimilar orientation,} \\
            0,                              & \mbox{otherwise.}
        \end{cases}
    \end{align}

    The entries of $q$-combinatorial Laplacian matrices are
    \begin{align}
        q>0,\ (\mathcal{L}_q)_{ij} &=
        \begin{cases}\label{equ:combine}
            \text{deg}(\sigma^i_{q}) + q + 1, & \mbox{if $i=j$.} \\
            1,                 & \mbox{if $i\neq j$, $\sigma^i_q \stackrel{U}\nsim \sigma^j_q$ and $\sigma^i_q \stackrel{L}\sim \sigma^j_q$ with similar orientation.}   \\
            -1,                & \mbox{if $i\neq j$, $\sigma^i_q \stackrel{U}\nsim \sigma^j_q$ and $\sigma^i_q \stackrel{L}\sim \sigma^j_q$ with dissimilar orientation.} \\
            0,                 & \mbox{if $i\neq j$ and either , $\sigma^i_q \stackrel{U}\sim \sigma^j_q$ or $\sigma^i_q \stackrel{L}\nsim \sigma^j_q$.}
        \end{cases} \\
      q=0,\ (\mathcal{L}_q)_{ij}  &=
      \begin{cases}\label{equ:L0}
            \text{deg}(\sigma^i_{0}),   & \mbox{if $i=j$.} \\
            -1,                         & \mbox{if $\sigma^i_0 \stackrel{U}\sim \sigma^j_0$.} \\
            0,                          & \mbox{otherwise.}
      \end{cases}
    \end{align}

    \subsubsection{Spectral analysis of  $q$-combinatorial Laplacian matrices}

     A $q$-combinatorial Laplacian matrix for oriented simplicial complexes is a generalization of the Laplacian matrix in graph theory. The spectra of a Laplacian matrix play an essential role in understanding the connectivity and robustness of simple graphs (simplicial complexes of dimension $1$). They can also distinguish different topological structures. Inspired by the capability of the Laplacian spectra of analyzing topological structures, we study the spectral properties of $q$-combinatorial Laplacian matrices to reveal topological and spectral information of simplicial complexes with dimension $0 \le q \le \text{dim}(K)$.

\begin{figure}[H]
        \centering
        \includegraphics[width=12cm]{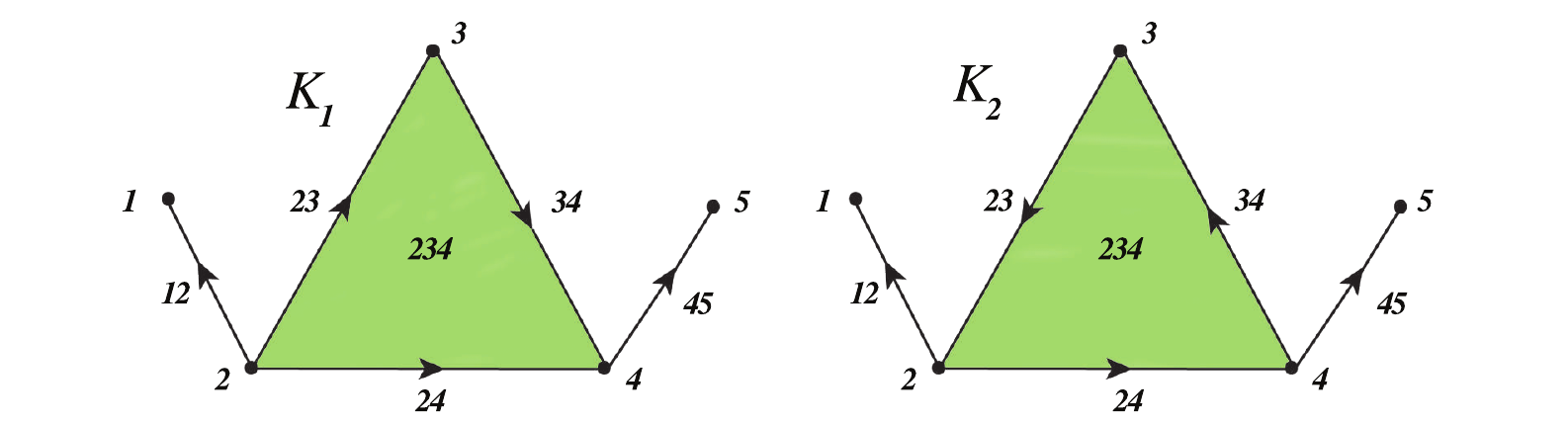}
        \captionsetup{margin=0.9cm}  
        \caption{Illustration of two oriented simplicial complexes with the same geometric structure but having different orientations. Here, we denote the vertices by $1,2,3,4,$ and $5$,  edges by $12,23,34,24,$ and $45$, and the triangle by $234$.}
        \label{fig:different orient}
    \end{figure}

      We clarify that for a given finite simplicial complex $K$, the spectra of its $q$-combinatorial Laplacian matrix is independent of the choice of the orientation for the $q$-simplices of $K$. The proof can be found in Ref.  \cite{goldberg2002combinatorial}. \autoref{fig:different orient} provides a simple example to illustrate this property. In \autoref{fig:different orient}, we have two oriented simplicial complexes, $K_1$ and $K_2$, with the same geometric structure but different orientations. For the sake of brevity, we use $1,2,3,4,$ and $5$ to represent $0$-simplices (vertices),  $12,23,34,24, $ and  $ 45$ to describe $1$-simplices (edges), and $234$ to stand for the $2$-simplex (triangle). Then the $0$-combinatorial Laplacian matrix of $K_1$ and $K_2$ is
    \[
    \mathcal{L}_0(K_1) = \mathcal{L}_0(K_2) =
    \left[\begin{array}{ccccc}
         1 & -1  &  0  &  0  &  0  \\
        -1 &  3  & -1  & -1  &  0  \\
         0 & -1  &  2  & -1  &  0  \\
         0 & -1  & -1  &  3  &  -1 \\
         0 &  0  &  0  & -1  &  1
    \end{array}\right].
    \]
    Obviously, $\mathcal{L}_0(K_1)$ and $\mathcal{L}_0(K_2)$ have the same spectra. For $q=1$, there are five $1$-simplices in $K_1$ and $K_2$, while $1$-combinatorial Laplacian matrices have  dimension $5 \times 5$.  Using $K_1$ as an example,  since $12$ and $23$ are lower adjacent with similar orientation, the element of $\mathcal{L}_1(K_1)$ addressed at first row and second column is $1$ according to  Eq. \eqref{equ:combine}. Since $34$ and $45$ are lower adjacent with dissimilar orientation, $(\mathcal{L}_1(K_1))_{35} = -1$. Moreover, $23$ and $34$ are upper adjacent which results in $(\mathcal{L}_1(K_1))_{23} = 0$. For the diagonal parts, $12$ and $45$ are not the faces of any $2$-simplex, while $23,34, $ and $24$ are the faces of $2$-simplex $234$. Therefore, $\text{deg}_U(12) = \text{deg}_U(12) =0$, and $\text{deg}_U(23) = \text{deg}_U(34) = \text{deg}_U(34) = 1$, so the diagonal terms of $\mathcal{L}_0(K_1)$ are $2,3,3,3,$ and $2$.

    \[
    \mathcal{L}_1(K_1) =
    \left[\begin{array}{ccccc}
         2 &  1  &  0  &  1  &  0  \\
         1 &  3  &  0  &  0  &  0  \\
         0 &  0  &  3  &  0  & -1  \\
         1 &  0  &  0  &  3  & -1 \\
         0 &  0  & -1  & -1  &  2
    \end{array}\right], \quad
    \mathcal{L}_1(K_2) =
    \left[\begin{array}{ccccc}
         2 & -1  &  0  &  1  &  0  \\
        -1 &  3  &  0  &  0  &  0  \\
         0 &  0  &  3  &  0  &  1  \\
         1 &  0  &  0  &  3  & -1 \\
         0 &  0  &  1  & -1  &  2
    \end{array}\right].
    \]
    The spectra of $\mathcal{L}_1(K_1)$ and $\mathcal{L}_1(K_2)$ have the same eigenvalues: $\{3,\dfrac{5 \pm \sqrt{5}}{2}, \dfrac{5 \pm \sqrt{13}}{2}\}$. Since $K_1$ and $K_2$ do not have $3$-simplices, $\mathcal{L}_q$ is a zero matrix when $q\ge 3$.

    A $q$-combinatorial Laplacian matrix is symmetric and positive semi-definite. Therefore, its eigenvalues are all real and non-negative. An analogy to the property that the number of zero eigenvalues of $\mathcal{L}_0$ represents the number of connected components ($\beta_0$) in the simple graph (simplicial complex with dimension $1$), the number of zero eigenvalues of $\mathcal{L}_q$ can also reveal the topological information. More specifically, for a given finite oriented simplicial complex, the  Betti number ${\beta_q}$ of $K$ satisfy
    \begin{equation}\label{equ:betti}
        \beta_q = \text{dim}(\mathcal{L}_q(K)) - \text{rank}(\mathcal{L}_q(K)) = \text{nullity}(\mathcal{L}_q(K)) = \# ~{\rm of ~ zero ~eigenvalues~ of}~   \mathcal{L}_q(K)
    \end{equation}
    We   consider tetrahedron-shaped structures in \autoref{fig:tetra_spectra} to illustrate the connection between Betti number and the dimension and the rank of $q$-combinatorial Laplacian matrix. For the sake of brevity, we will use $i$ to represents $0$-simplex $[v_i]$, $ij$ to represents $1$-simplex $[v_i, v_j]$, and $ijk$ to represents $[v_i, v_j, v_k]$. Then, $1$- and $2$-boundary operators map:
    \[
    \begin{split}
        \partial_1(ij)  &= j - i, \\
        \partial_2(ijk) &= jk - ik + ij,
    \end{split}
    \]
    Since  different orientations result in the same spectrum, there is no need  to label the orientation in \autoref{fig:tetra_spectra}. In the following, we analyze three tetrahedron-shaped simplicial complexes:

    \begin{figure}[H]
        \centering
        \includegraphics[width=12cm]{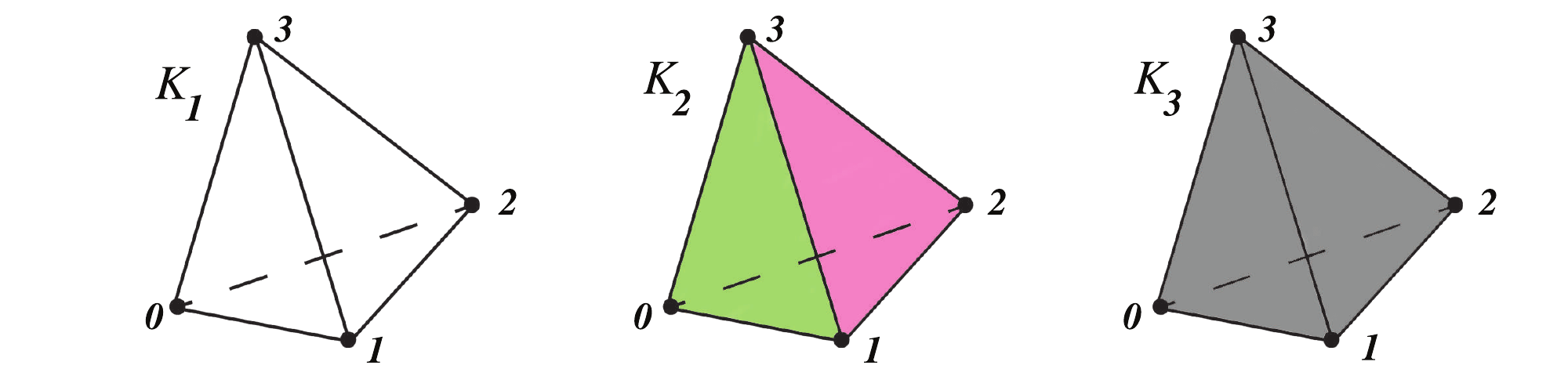}
        \captionsetup{margin=0.9cm}  
        \caption{Illustration of three different tetrahedron-shaped simplicial complexes. There are four $0$-simplices and six $1$-simplices  in $K_1$.  Here, $K_2$ has four more $2$-simplices than $K_1$ does, while $K_3$ owns one more $3$-simplex than $K_2$ does.}
        \label{fig:tetra_spectra}
    \end{figure}

    \begin{itemize}
        \item[$\boldsymbol{K_1}$. ] The left most chart in \autoref{fig:tetra_spectra} has four $0$-simplices: $0,1,2,$ and $3,$ and six $1$-simplices: $01, 02, 03, 12, 13, $ and $23$. It is clear that $C_{q}(K_1)$ is an empty set and $\partial_q$ is an zero map when $q\ge 2$. Then, its Laplacian operators are
        \[
        \Delta_1 = \partial_1^{\ast} \partial_1, \ \Delta_0 = \partial_1 \partial_1^{\ast} + \partial_0^{\ast} \partial_0.
        \]
        The combinatorial Laplacian matrices are:
        \[
        \mathcal{L}_1 = \mathcal{B}_1^{T} \mathcal{B}_1, \ \mathcal{L}_0 = \mathcal{B}_1 \mathcal{B}_1^{T} + \mathcal{B}_{0}^T \mathcal{B}_{0}.
        \]
        The matrix representation $\mathcal{B}_1$ for $\partial_1: C_1(K_1) \longrightarrow C_0(K_1)$ is:
        \begin{equation}\label{equ:K1}
            \mathcal{B}_1 =
                \begin{array}{@{}r@{}c@{}c@{}c@{}c@{}c@{}c@{}l@{}}
                    & 01 & 02 & 03 & 12 & 13 & 23  \\
                   \left.\begin{array}{c}
                    0 \\
                    1 \\
                    2 \\
                    3
                    \end{array}\right[
                    & \begin{array}{c} -1 \\  1  \\  0  \\  0    \end{array}
                    & \begin{array}{c} -1 \\  0  \\  1  \\  0    \end{array}
                    & \begin{array}{c} -1 \\  0  \\  0  \\  1    \end{array}
                    & \begin{array}{c}  0 \\ -1  \\  1  \\  0    \end{array}
                    & \begin{array}{c}  0 \\ -1  \\  0  \\  1    \end{array}
                    & \begin{array}{c}  0 \\  0  \\ -1  \\  1    \end{array}
                    & \left]. \begin{array}{c} \\ \\  \\ \\
                    \end{array}\right.
            \end{array}
        \end{equation}
        and $B_0$ is 
        \[
        B_0 = 
        \begin{array}{@{}r@{}c@{}c@{}c@{}c@{}l@{}}
                    & 0 & 1 & 2 & 3   \\
                   \left.\begin{array}{c}

                    \end{array}\right[
                    & \begin{array}{c}  0  \end{array}
                    & \begin{array}{c}  0  \end{array}
                    & \begin{array}{c}  0  \end{array}
                    & \begin{array}{c}  0  \end{array}
                    & \left]\begin{array}{c}
                    \end{array}\right.
                \end{array}
        \]
        Therefore, the associated combinatorial Laplacian matrices are
        \[
        \mathcal{L}_1(K_1) =
        \left[\begin{array}{cccccc}
             2  &   1  &   1  &  -1  &  -1  &   0 \\
             1  &   2  &   1  &   1  &   0  &  -1 \\
             1  &   1  &   2  &   0  &   1  &   1 \\
            -1  &   1  &   0  &   2  &   1  &  -1 \\
            -1  &   0  &   1  &   1  &   2  &   1 \\
             0  &  -1  &   1  &  -1  &   1  &   2
        \end{array}\right], \
        \mathcal{L}_0(K_1) =
        \left[\begin{array}{cccc}
             3  &  -1  &  -1  &  -1  \\
            -1  &   3  &  -1  &  -1  \\
            -1  &  -1  &   3  &  -1  \\
            -1  &  -1  &  -1  &   3
        \end{array}\right].
        \]
As shown in        \autoref{table:K1}, we can calculate the spectra and ranks from  combinatorial Laplacian matrices. We have $\beta_0 = 1$, $\beta_1 = 3$, which reveal that one connected component and three $1$-cycles are exist in $K_1$.
        \begin{table}[H]
            \centering
            \setlength\tabcolsep{18pt}
            \captionsetup{margin=0.9cm}
            \caption{Table of dimensions, ranks, nullity, spectra and Betti numbers of combinatorial Laplacian matrices $\mathcal{L}_0, $ and $\mathcal{L}_1$ for simplicial complex $K_1$.}
            \begin{tabular}{ccccccc}
            \hline
                             & $\mathcal{L}_1(K_1)$ & $\mathcal{L}_0(K_1)$  \\ \hline
             Betti number    & $\beta_1 = 3$        & $\beta_0 = 1$       \\
             dim             & $6$                  & $4$                 \\
             rank            & $3$                  & $3$                 \\
             nullity         & $3$                  & $1$                 \\
             Spectra         & $\{0,0,0,4,4,4\}$      & $\{0,4,4,4\}$         \\ \hline
            \end{tabular}
            \label{table:K1}
        \end{table}

        \item[$\boldsymbol{K_2}$.] We analyze the middle chart in \autoref{fig:tetra_spectra} in a similar manner. As one can see, $K_2$ has four $0$-simplices: $0,1,2,$ and $3$, six $1$-simplices: $01, 02, 03, 12, 13, $ and $ 23$, and four $2$-simplices: $012, 013, 023, $ and $123$. The associated Laplacian operators are
        \[
        \Delta_2 = \partial_2^{\ast} \partial_2,\ \Delta_1 = \partial_2 \partial_2^{\ast} + \partial_1^{\ast} \partial_1, \ \Delta_0 = \partial_1 \partial_1^{\ast} + \partial_0^{\ast} \partial_0.
        \]
        The resulting combinatorial Laplacian matrices are
        \[
        \mathcal{L}_2 = \mathcal{B}_2^T \mathcal{B}_2, \ \mathcal{L}_1 = \mathcal{B}_2 \mathcal{B}_2^T + \mathcal{B}_1^T \mathcal{B}_1, \  \mathcal{L}_0 = \mathcal{B}_1 \mathcal{B}_1^{T} + \mathcal{B}_{0}^T \mathcal{B}_{0}.
        \]
        The corresponding matrix representations for $\mathcal{B}_2$ and $\mathcal{B}_1$ are respectively
        \begin{equation}\label{equ:K2}
            \mathcal{B}_2 =
            \begin{array}{@{}r@{}c@{}c@{}c@{}c@{}l@{}}
                & 012 & 013 & 023 & 123   \\
               \left.\begin{array}{c}
                01 \\
                02 \\
                03 \\
                12 \\
                13 \\
                23
                \end{array}\right[
                & \begin{array}{c}  1 \\ -1  \\  0  \\  1  \\   0  \\  0 \end{array}
                & \begin{array}{c}  1 \\  0  \\ -1  \\  0  \\   1  \\  0 \end{array}
                & \begin{array}{c}  0 \\  1  \\ -1  \\  0  \\   0  \\  1 \end{array}
                & \begin{array}{c}  0 \\  0  \\  0  \\  1  \\  -1  \\  1 \end{array}
                & \left],\begin{array}{c} \\ \\  \\ \\ \\ \\
                \end{array}\right.
            \end{array}
            \mathcal{B}_1 =
            \begin{array}{@{}r@{}c@{}c@{}c@{}c@{}c@{}c@{}l@{}}
                & 01 & 02 & 03 & 12 & 13 & 23  \\
               \left.\begin{array}{c}
                0 \\
                1 \\
                2 \\
                3
                \end{array}\right[
                & \begin{array}{c} -1 \\  1  \\  0  \\  0    \end{array}
                & \begin{array}{c} -1 \\  0  \\  1  \\  0    \end{array}
                & \begin{array}{c} -1 \\  0  \\  0  \\  1    \end{array}
                & \begin{array}{c}  0 \\ -1  \\  1  \\  0    \end{array}
                & \begin{array}{c}  0 \\ -1  \\  0  \\  1    \end{array}
                & \begin{array}{c}  0 \\  0  \\ -1  \\  1    \end{array}
                & \left].\begin{array}{c} \\ \\  \\ \\
                \end{array}\right.
            \end{array}
        \end{equation}
        Then,  associated combinatorial Laplacian matrices are
        \[
        \mathcal{L}_2(K_2) =
        \left[\begin{array}{cccc}
             3  &   1  &  -1  &   1  \\
             1  &   3  &   1  &  -1  \\
            -1  &   1  &   3  &   1  \\
             1  &  -1  &   1  &   3
        \end{array}\right], \
        \mathcal{L}_1(K_2) =
        \left[\begin{array}{cccccc}
             4  &   0  &   0  &   0  &   0  &   0 \\
             0  &   4  &   0  &   0  &   0  &   0 \\
             0  &   0  &   4  &   0  &   0  &   0 \\
             0  &   0  &   0  &   4  &   0  &   0 \\
             0  &   0  &   0  &   0  &   4  &   0 \\
             0  &   0  &   0  &   0  &   0  &   4
        \end{array}\right],
        \]
				and $\mathcal{L}_0(K_2)=\mathcal{L}_0(K_1)$.
        Similarly, from \autoref{table:K2}, we see that there are one connected component and one $2$-cycle (void)  in $K_2$.

        \begin{table}[H]
            \centering
            \setlength\tabcolsep{11pt}
            \captionsetup{margin=0.9cm}
            \caption{Table of dimensions, ranks, nullity, spectra and Betti numbers of combinatorial Laplacian matrices $\mathcal{L}_0, \mathcal{L}_1,$ and $ \mathcal{L}_2$ for simplicial complex $K_2$.}
            \begin{tabular}{ccccccc}
            \hline
                             & $\mathcal{L}_2(K_2)$ & $\mathcal{L}_1(K_2)$ &  $\mathcal{L}_0(K_2)$ \\ \hline
             Betti number    & $\beta_2 = 1$        & $\beta_1 = 0$        & $\beta_1 = 1$         \\
             dim             & $4$                  & $6$                  & $4$                   \\
             rank            & $3$                  & $6$                  & $3$                   \\
             nullity         & $1$                  & $0$                  & $1$                   \\
             Spectra         & $\{0,4,4,4\}$          & $\{4,4,4,4,4,4\}$      & $\{0,4,4,4\}$           \\ \hline
            \end{tabular}
            \label{table:K2}
        \end{table}

        \item[$\boldsymbol{K_3}$.]  $K_3$ in the right most chart in \autoref{fig:tetra_spectra}  has four $0$-simplices: $0,1,2,$ and $3$, six $1$-simplices: $01, 02, 03, 12, 13, $ and $23$, four $4$-simplices: $012, 013, 023, $ and $123$, and one $3$-simplex $01234$. The associated Laplacian operators are
        \[
        \Delta_2 = \partial_3 \partial_3^{\ast} + \partial_2^{\ast} \partial_2, \ \Delta_1 = \partial_2 \partial_2^{\ast} + \partial_1^{\ast} \partial_1, \ \Delta_0 = \partial_1 \partial_1^{\ast} + \partial_0^{\ast} \partial_0.
        \]
        The corresponding combinatorial Laplacian matrices are:
        \[
        \mathcal{L}_2 = \mathcal{B}_3 \mathcal{B}_3^T + \mathcal{B}_2^T \mathcal{B}_2, \ \mathcal{L}_1 = \mathcal{B}_2 \mathcal{B}_2^T + \mathcal{B}_1^T \mathcal{B}_1, \  \mathcal{L}_0 = \mathcal{B}_1 \mathcal{B}_1^{T} + \mathcal{B}_{0}^T \mathcal{B}_{0}.
        \]
        The resulting matrix representations $\mathcal{B}_1 $ and $ \mathcal{B}_2$ are exactly the same as those in Eq. \eqref{equ:K2}. While $\mathcal{B}_3$ is given by
        \begin{equation}\label{equ:K3}
            \mathcal{B}_3 =
            \begin{array}{@{}r@{}c@{}l@{}}
                & 0123   \\
               \left.\begin{array}{c}
                012 \\
                013 \\
                023 \\
                123
                \end{array}\right[
                & \begin{array}{c} -1 \\ 1 \\ -1 \\ 1 \end{array}
                & \left].
								\begin{array}{c} \\ \\ \\ \\
                \end{array}\right.
            \end{array}
        \end{equation}
        Therefore, the $2$-combinatorial Laplacian matrices are
        \[
        \mathcal{L}_2(K_3) =
        \left[\begin{array}{cccc}
             4  &   0  &   0  &   0  \\
             0  &   4  &   0  &   0  \\
             0  &   0  &   4  &   0  \\
             0  &   0  &   0  &   4
        \end{array}\right],
        \]
				  $\mathcal{L}_1(K_3)=\mathcal{L}_1(K_2)$, and $\mathcal{L}_0(K_3)=\mathcal{L}_0(K_2)$.
        In this case,   \autoref{table:K3} shows that there is only one connected component in $K_3$, which matches the right most chart in \autoref{fig:tetra_spectra}.
        \begin{table}[H]
            \centering
            \setlength\tabcolsep{11pt}
            \captionsetup{margin=0.9cm}
            \caption{Table of dimensions, ranks, nullity, spectra and Betti numbers of combinatorial Laplacian matrices $\mathcal{L}_0, \mathcal{L}_1,$ and  $\mathcal{L}_2$ for simplicial complex $K_3$.}
            \begin{tabular}{ccccccc}
            \hline
                             & $\mathcal{L}_2(K_3)$ & $\mathcal{L}_1(K_3)$ &  $\mathcal{L}_0(K_3)$ \\ \hline
             Betti number    & $\beta_2 = 0$        & $\beta_1 = 0$        & $\beta_0 = 1$         \\
             dim             & $4$                  & $6$                  & $4$                   \\
             rank            & $4$                  & $6$                  & $3$                   \\
             nullity         & $0$                  & $0$                  & $1$                   \\
             Spectra         & $\{4,4,4,4\}$          & $\{4,4,4,4,4,4\}$      & $\{0,4,4,4\}$           \\ \hline
            \end{tabular}
            \label{table:K3}
        \end{table}
    \end{itemize}

    \subsubsection{Persistent spectral theory}
    Instead of using the aforementioned spectral analysis for $q$-combinatorial Laplacian matrix to describe a single configuration, we propose a persistent spectral theory to create a sequence of simplicial complexes induced by varying a filtration parameter, which is inspired by  persistent homology  and our earlier work in multiscale graphs \cite{xia2015multiscale,bramer2018multiscale}.

    A filtration of an oriented simplicial complex $K$ is a sequence of sub-complexes $(K_t)_{t=0}^m$ of $K$
    \begin{equation}
        \emptyset = K_0 \subseteq K_1 \subseteq K_2 \subseteq \cdots \subseteq K_m = K.
    \end{equation}
		It induces a sequence of  chain complexes
    \begin{equation}
        \left.\begin{array}{cccccccccccccc}
            \cdots & C_{q+1}^1 &
            \xrightleftharpoons[\partial_{q+1}^{1^\ast}]{\partial_{q+1}^1} & C_q^1 &
            \xrightleftharpoons[\partial_q^{1^\ast}]{\partial_q^1} & \cdots & \xrightleftharpoons[\partial_3^{1^\ast}]{\partial_3^1} & C_2^1 & \xrightleftharpoons[\partial_2^{1^\ast}]{\partial_2^1} & C_1^1 & \xrightleftharpoons[\partial_1^{1^\ast}]{\partial_1^1} & C_0^1 & \xrightleftharpoons[\partial_0^{1^\ast}]{\partial_0^1} & C_{-1}^1 \\
            & \rotatebox{-90}{$\subseteq$} &  & \rotatebox{-90}{$\subseteq$} &  &  &  & \rotatebox{-90}{$\subseteq$} &  & \rotatebox{-90}{$\subseteq$} &  & \rotatebox{-90}{$\subseteq$} &  &  \\
            \cdots & C_{q+1}^2 &
            \xrightleftharpoons[\partial_{q+1}^{2^\ast}]{\partial_{q+1}^2} & C_q^2 &
            \xrightleftharpoons[\partial_q^{2^\ast}]{\partial_q^2} & \cdots &
            \xrightleftharpoons[\partial_3^{2^\ast}]{\partial_3^2} & C_2^2 & \xrightleftharpoons[\partial_2^{2^\ast}]{\partial_2^2} & C_1^2 & \xrightleftharpoons[\partial_1^{2^\ast}]{\partial_1^2} & C_0^2 & \xrightleftharpoons[\partial_0^{2^\ast}]{\partial_0^2} & C_{-1}^1 \\
            & \rotatebox{-90}{$\subseteq$} &  & \rotatebox{-90}{$\subseteq$} &  &  &  & \rotatebox{-90}{$\subseteq$} &  & \rotatebox{-90}{$
            \subseteq$} &  & \rotatebox{-90}{$\subseteq$} &  &  \\
            \cdots & C_{q+1}^m &
            \xrightleftharpoons[\partial_{q+1}^{m^\ast}]{\partial_{q+1}^m} & C_q^m &
            \xrightleftharpoons[\partial_q^{m^\ast}]{\partial_q^m} & \cdots &
            \xrightleftharpoons[\partial_3^{m^\ast}]{\partial_3^m} & C_2^m & \xrightleftharpoons[\partial_2^{m^\ast}]{\partial_2^m} & C_1^m & \xrightleftharpoons[\partial_1^{m^\ast}]{\partial_1^m} & C_0^m & \xrightleftharpoons[\partial_0^{m^\ast}]{\partial_0^m} & C_{-1}^1
        \end{array}\right.
    \end{equation}
    where $C_q^t \coloneqq C_q(K_t)$ and $\partial_q^t \colon  C_q(K_t) \to C_{q-1}(K_t)$. Each $K_t$ itself is an oriented simplicial complex which has dimension denoted by $\text{dim}(K_t)$. If $q<0$, then $C_{q}(K_t) = \{ \emptyset \} $ and $\partial_q^t$ is actually a zero map.\footnote{We define the boundary matrix $\mathcal{B}_0^t$ for boundary map $\partial_0^t$ as a zero matrix. The number of columns of $\mathcal{B}_0^t$ is the number of $0$-simplices in $K_t$, the number of rows will be $1$.} For a general case of $0 < q \le \text{dim}(K_t)$, if $\sigma_q$ is an oriented $q$-simplex of $K_t$, then
    \[
    \partial_q^t(\sigma_q) = \sum_{i}^q(-1)^i \sigma^i_{q-1}, \sigma_q \in K_t,
    \]
    with $\sigma_q = [v_0, \cdots, v_q]$ being the oriented $q$-simplex, and $\sigma^{i}_{q-1} = [v_0, \cdots, \hat{v_i} ,\cdots,v_q]$ being the oriented $(q-1)$-simplex for which  its vertex $v_i$ is removed.

    Let $\mathbb{C}_{q}^{t+p}$ be the subset of $C_q^{t+p}$ whose boundary is in $C_{q-1}^t$:
    \begin{equation}
        \mathbb{C}_q^{t+p} \coloneqq \{ \alpha \in C_q^{t+p} \ | \ \partial_q^{t+p}(\alpha) \in C_{q-1}^{t}\}.
    \end{equation}
    We define
    \begin{equation}\label{equ: new boundary map}
        \eth_q^{t+p} : \mathbb{C}_q^{t+p} \to  C_{q-1}^{t}
    \end{equation}
    
    Based on the $q$-combinatorial Laplacian operator, the ${p}$-persistent ${q}$-combinatorial  Laplacian operator $\Delta_q^{t+p}: C_q(K_{t}) \to C_q(K_t)$ defined along the filtration can be expressed as
    \begin{equation}
        \Delta_q^{t+p} = \eth_{q+1}^{t+p} \left(\eth_{q+1}^{t+p}\right)^\ast + \partial_q^{t^\ast} \partial_q^t.
    \end{equation}
    We denote the matrix representations of boundary operator $\eth_{q+1}^{t+p}$ and $\partial_q^t$ by $\mathcal{B}_{q+1}^{t+p}$ and $\mathcal{B}_{q}^t$, respectively. It is clear that the number of rows in $\mathcal{B}_{q+1}^{t+p}$ is the number of oriented $q$-simplices in $K_t$, and the number of columns is the number of oriented $(q+1)$-simplices in $K_{t+p} \cap \mathbb{C}_{q+1}^{t+p}$. The transpose of the matrices $\mathcal{B}_{q+1}^{t+p}$ and $\mathcal{B}_{q}^t$, are the matrix representations of the adjoint boundary operator $\left(\partial_{q+1}^{t+p}\right)^\ast$ and $\partial_q^{t^\ast}$, respectively. Therefore, the  ${p}$-persistent ${q}$-combinatorial Laplacian matrix, $\mathcal{L}_q^{t+p}$, is
    \begin{equation}\label{equ:PLC}
        \mathcal{L}_q^{t+p} = \mathcal{B}_{q+1}^{t+p} (\mathcal{B}_{q+1}^{t+p})^T + (\mathcal{B}_{q}^t)^T \mathcal{B}_{q}^t.
    \end{equation}
    Intuitively, for a non-empty set $C_q^t$, the $p$-persistent $q$-combinatorial Laplacian matrix $\mathcal{L}_q^{t+p}$ is a square matrix with dimension to be the number of $q$-simplices in $K_t$. Moreover, $\mathcal{L}_q^{t+p}$ is symmetric and positive semi-defined and thus, all the spectra of $\mathcal{L}_q^{t+p}$ are real and non-negative. If $p=0$, then $\mathcal{L}_q^{t+0}$ is exactly the $q$-combinatorial Laplacian matrix defined in Eq. \eqref{equ:combinatorial Laplacian}.

    We are interested in the difference between $\mathcal{L}_q^{t+0}$ and $\mathcal{L}_q^{t+p}$. Suppose we have an oriented simplicial complex $K_t$, and also an oriented simplicial complex $K_{t+p}$ constructed by adding different dimensions simplices ("outer" topological structures) to $K_t$ with $\text{dim}(K_{t+p})=q+1$. Since $K_t \subset K_{t+p}$, we have
    \[
    \begin{split}
        \mathcal{L}_q^{t+0} &= \mathcal{B}_{q+1}^{t+0} (\mathcal{B}_{q+1}^{t+0})^T + (\mathcal{B}_{q}^t)^T \mathcal{B}_{q}^t \\
        \mathcal{L}_q^{t+p} &= \mathcal{B}_{q+1}^{t+p} (\mathcal{B}_{q+1}^{t+p})^T + (\mathcal{B}_{q}^t)^T \mathcal{B}_{q}^t.
    \end{split}
    \]
    
    \begin{enumerate}
        \item[Case 1]. If $(\text{Im}(\eth_{q+1}^{t+p}) \cap (K_{t+p}\backslash K_t)) \cap C_q^t = \emptyset$ for all possible $q$, then the "outer" topological structures are disconnected with $K_t$. Therefore, the boundary matrix $\mathcal{B}_{q+1}^{t+p}$ will add at least one additional zero column, which leads to at least one additional zero row for $(\mathcal{B}_{q+1}^{t+p})^T$. However, the matrix obtained by $\mathcal{B}_{q+1}^{t+p} (\mathcal{B}_{q+1}^{t+p})^T$ does not change and $\mathcal{L}_q^{t+p}$ is exactly the same as $\mathcal{L}_q^{t+0}$. In this situation, the spectra of $\mathcal{L}_q^{t+0}$ and $\mathcal{L}_q^{t+p}$ are exactly the same, which reveals the fact that the topological structure $K_t$ does not change  under the filtration process.

        \item[Case 2]. If $(\text{Im}(\eth_{q+1}^{t+p}) \cap (K_{t+p}\backslash K_t)) \cap C_q^t \neq \emptyset$ for at least one   $q$, then the boundary matrix $\mathcal{B}_{q+1}^{t+p}$ will be changed by adding additional non-zero columns. Therefore, $\mathcal{L}_q^{t+p}$ is no longer the same as $\mathcal{L}_q^{t+0}$, but the structure information in $K_t$ will still be preserved in $\mathcal{L}_q^{t+p}$. In this case, the topological structure $K_t$ builds connection with ``outer'' topological structures. By calculating the spectra of $\mathcal{L}_q^{t+p}$, the disappeared and preserved structure information of $K_t$ under the filtration process can be revealed.
    \end{enumerate}

    Based on the fact that the topological and spectral information of $K_t$ can also be analyzed from $\mathcal{L}_q(K_t)$ along with the filtration parameter by diagonalizing the  $q$-combinatorial Laplacian matrix, we  focus on the spectra information calculated from $\mathcal{L}_q^{t+p}$.   Denote the set of spectra of $\mathcal{L}_q^{t+p}$ by
    \[
    \text{Spectra}(\mathcal{L}_q^{t+p}) = \{(\lambda_1)_q^{t+p}, (\lambda_2)_q^{t+p}, \cdots, (\lambda_N)_q^{t+p}  \},
    \]
    where $\mathcal{L}_q^{t+p}$ has dimension $N\times N$ and spectra are arranged in ascending order. The smallest non-zero eigenvalue of $\mathcal{L}_q^{t+p}$ is defined as $(\tilde{\lambda}_2)_q^{t+p}$.  In the previous section, we have seen that Betti numbers (i.e $\#$ of zero eigenvalues) can reveal $q$-cycle information. Similarly, we  define the number of zero eigenvalues of $p$-persistent $q$-combinatorial Laplacian matrix $\mathcal{L}_q^{t+p}$ to be the  ${p}$-persistent ${q}$th Betti numbers
    \begin{equation}\label{equ:persistent betti}
        \beta_q^{t+p} = \text{dim}(\mathcal{L}_{q}^{t+p}) - \text{rank}(\mathcal{L}_{q}^{t+p}) = \text{nullity}(\mathcal{L}_{q}^{t+p}) = \boldsymbol{\#} ~{\rm of ~zero ~eigenvalues ~of ~}   \mathcal{L}_{q}^{t+p}.
    \end{equation}
 In fact, $\beta_q^{t+p}$ counts the number of $q$-cycles in $K_t$ that are still alive in $K_{t+p}$, which exactly provides the same topological information as persistent homology does. However, persistent spectral theory offers additional geometric information from the spectra of persistent combinatorial Laplacian matrix beyond topological  persistence. In general, the topological changes can be read off from  persistent   Betti numbers (harmonic persistent spectra) and the geometric changes can be derived from the non-harmonic persistent spectra.

\begin{figure}[H]
    \centering
    \includegraphics[width=16cm]{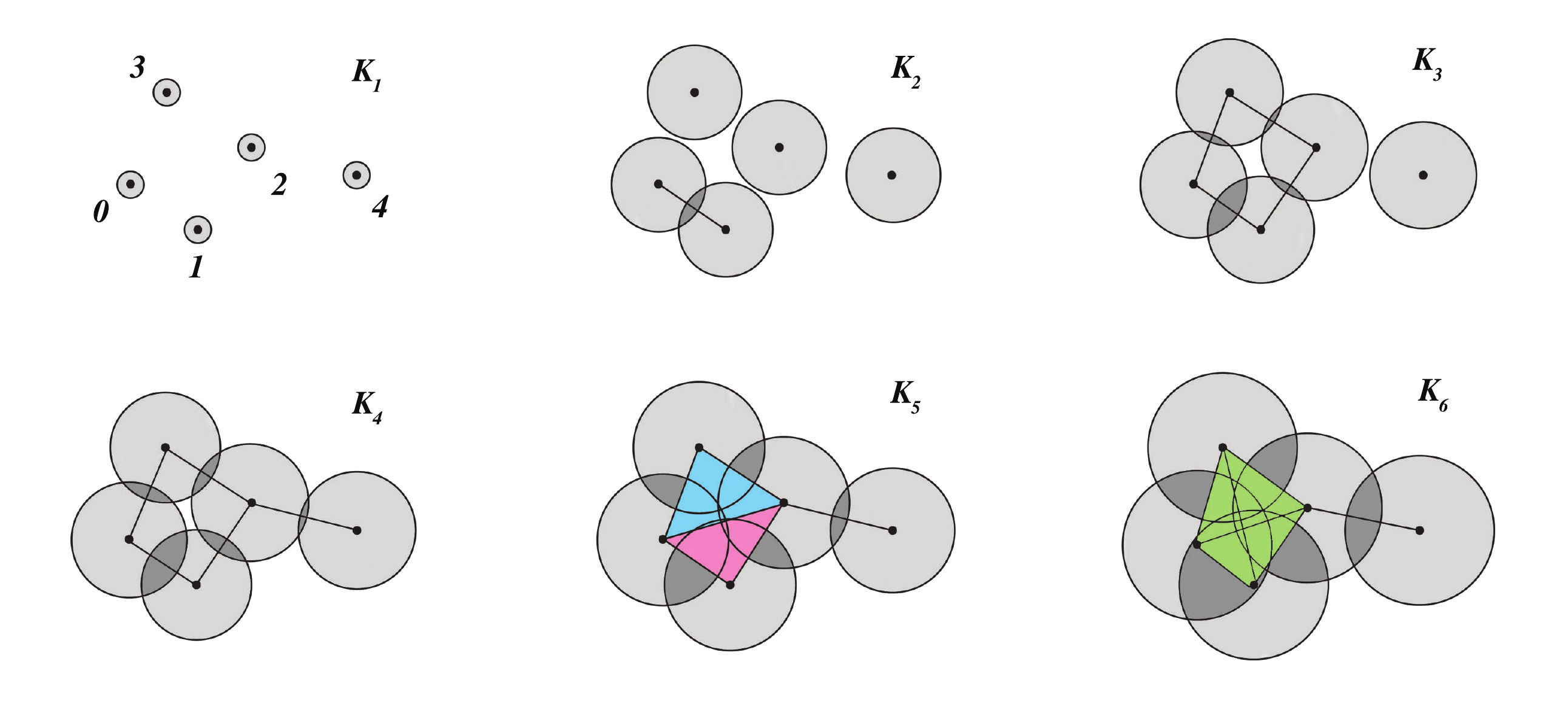}
    \captionsetup{margin=0.9cm}  
    \caption{Illustration of filtration. We use $0, 1, 2, 3, $ and $ 4$ to stand for $0$-simplices, $01, 12, 23,03, 24, 02, $ and $13$ for $1$-simplices, $012, 023, 013, $ and $123$ for $2$-simplices, and $0123$ for the $3$-simplex. Here, $K_1$ has five 0-cycles, $K_2$ has four 0-cycles, $K_3$ has two $0$-cycles and a $1$-cycle, $K_4$ has a $0$-cycle and a $1$-cycle, $K_5$ has one $0$-cycle, and $K_6$ has a $0$-cycle.}
    \label{fig:filtration}
    \end{figure}

    \autoref{fig:filtration} demonstrates an example of a standard filtration process. Here the initial setup $K_1$ consists of five $0$-simplices (vertices). We construct Vietoris–Rips complexes by using an ever-growing circle centered at each vertex with radius $r$.  Once two circles overlapped with each other, an $1$-simplex (edge) is formed. A $2$-simplex (triangle) will be created when $3$ circles contact with one another, and a $3$-simplex will be generated once $4$ circles get overlapped one another. As \autoref{fig:filtration} shows, we can attain a series of simplicial complexes from $K_1$ to $K_6$ with the radius of circles increasing. To fully illustrate how to construct  $p$-persistent $q$-combinatorial Laplacian matrices by the boundary operator and determine persistent Betti numbers, we  analyze $6$ $p$-persistent $q$-combinatorial Laplacian matrices and their corresponding harmonic persistent spectra (i.e., persistent Betti numbers) and non-harmonic persistent spectra. Additional matrices are analyzed in Appendix \autoref{app:examples}.

    \begin{table}[H]
            \centering
            \setlength\tabcolsep{11pt}
            \captionsetup{margin=0.9cm}
            \caption{The number of $q$-cycles of simplicial complexes demonstrated in \autoref{fig:filtration}.}
            \begin{tabular}{ccccccc}
            \hline
             $\#$ of $q$-cycles   & $K_1$  & $K_2$  & $K_3$  & $K_4$  & $K_5$  & $K_6$   \\ \hline
             $q=0$                & $5$    & $4$    & $2$    & $1$    & $1$    & $1$        \\
             $q=1$                & $0$    & $0$    & $1$    & $1$    & $0$    & $0$        \\
             $q=2$                & $0$    & $0$    & $0$    & $0$    & $0$    & $0$        \\ \hline
            \end{tabular}
            \label{table:cycle number}
        \end{table}

    \begin{enumerate}
        \item[Case $1$.] In this case, the initial setup is $K_1$ and the end status is $K_3$. Therefore, $t=1$ and $p=2$ in Eq. \eqref{equ:PLC}. We will calculate $\mathcal{L}_0^{1+2}, \mathcal{L}_1^{1+2}, $ and $\mathcal{L}_2^{1+2}$ first and find out their corresponding persistent spectra.

        The $2$-persistent $0,1,2$-combinatorial Laplacian operators are:
        \[
        \begin{split}
        \Delta_0^{1+2} &= \eth_{1}^{1+2} \left(\eth_{1}^{1+2}\right)^\ast + \partial_0^{1^\ast} \partial_0^1, \\
        \Delta_1^{1+2} &= \eth_{2}^{1+2} \left(\eth_{2}^{1+2}\right)^\ast + \partial_1^{1^\ast} \partial_1^1, \\
        \Delta_2^{1+2} &= \eth_{3}^{1+2} \left(\eth_{3}^{1+2}\right)^\ast + \partial_2^{1^\ast} \partial_2^1,
        \end{split}
        \]
        Since $2$-simplex and $3$-simplex do not exist in $K_1 $ and $ K_3$, $\eth_2^{1+2}, \partial_1^1, \eth_3^{1+2}, $ and $ \partial_2^{1}$ do not exist and  $\partial_0^1$ is a zero map.  Then, there is  only one persistent combinatorial Laplacian matrix
        \[
            \mathcal{L}_0^{1+2} = \mathcal{B}_{1}^{1+2} (\mathcal{B}_{1}^{1+2})^T + (\mathcal{B}_0^{1})^T \mathcal{B}_0^{1}.
        \]
        It can be seen in  \autoref{fig:filtration} that two $0$-cycles (connected components) in $K_1$ are still alive in $K_3$, while no $1$-cycle and $2$-cycle exist in the initial set up $K_1$, which perfectly match the calculations in \autoref{table:case 1}: $\beta_0^{1+2} = 2$.
        \begin{table}[H]
            \centering
            \setlength\tabcolsep{14pt}
            \captionsetup{margin=0.9cm}
            \caption{$K_1 \to K_3$}
            \begin{tabular}{c|ccc}
            \hline
            $q$  & $q=0$ & $q=1$ & $q=2$ \\ \hline \\
            $\mathcal{B}_{q+1}^{1+2}$  & $\begin{array}{@{}r@{}c@{}c@{}c@{}c@{}l@{}}
                    & 01 & 12 & 23 & 03   \\
                   \left.\begin{array}{c}
                    0 \\
                    1 \\
                    2 \\
                    3 \\
                    4
                    \end{array}\right[
                    & \begin{array}{c} -1 \\  1  \\  0  \\  0 \\ 0 \end{array}
                    & \begin{array}{c}  0 \\ -1  \\  1  \\  0 \\ 0 \end{array}
                    & \begin{array}{c}  0 \\  0  \\ -1  \\  1 \\ 0 \end{array}
                    & \begin{array}{c} -1 \\  0  \\  0  \\  1 \\ 0 \end{array}
                    & \left]\begin{array}{c} \\ \\  \\ \\ \\
                    \end{array}\right.
                \end{array}$ & / &  /    \\  \\
            $\mathcal{B}_{q}^{1}$  &  $\begin{array}{@{}r@{}c@{}c@{}c@{}c@{}c@{}l@{}}
                    & 0 & 1 & 2 & 3 & 4   \\
                   \left.\begin{array}{c}
                    /
                    \end{array}\right[
                    & \begin{array}{c}  0  \end{array}
                    & \begin{array}{c}  0  \end{array}
                    & \begin{array}{c}  0  \end{array}
                    & \begin{array}{c}  0  \end{array}
                    & \begin{array}{c}  0  \end{array}
                    & \left]\begin{array}{c} \\
                    \end{array}\right.
                \end{array}$ & / &  /    \\  \\
            $\mathcal{L}_{q}^{1+2}$  & $\left[\begin{array}{ccccc}
                 2  &  -1  &   0  &  -1  &  0  \\
                -1  &   2  &  -1  &   0  &  0  \\
                 0  &  -1  &   2  &  -1  &  0  \\
                -1  &   0  &  -1  &   2  &  0  \\
                 0  &   0  &   0  &   0  &  0
            \end{array}\right]$  & / &  /    \\  \\
            $\beta_{q}^{1+2}$                        & 2               & / & /    \\  \\
            $\text{dim}(\mathcal{L}_{q}^{1+2})$      & 5               & / & /    \\  \\
            $\text{rank}(\mathcal{L}_{q}^{1+2})$     & 3               & / & /    \\  \\
            $\text{nullity}(\mathcal{L}_{q}^{1+2})$  & 2               & / & /  \\  \\
            $\text{Spectrum}(\mathcal{L}_{q}^{1+2})$    & $\{0,0,2,2,4\}$   & / & /  \\ \hline
            \end{tabular}
            \\
            \label{table:case 1}
        \end{table}

        \item[Case $2$.] The initial setup is $K_3$ and the end status is $K_4$.
        The $1$-persistent $0,1,2$-combinatorial Laplacian operators are
        \[
        \begin{split}
            \Delta_0^{3+1} &= \eth_{1}^{3+1} \left(\eth_{1}^{3+1}\right)^\ast + \partial_0^{3^\ast} \partial_0^3, \\
            \Delta_1^{3+1} &= \eth_{2}^{3+1} \left(\eth_{2}^{3+1}\right)^\ast + \partial_1^{3^\ast} \partial 1^3, \\
            \Delta_2^{3+1} &= \eth_{3}^{3+1} \left(\eth_{3}^{3+1}\right)^\ast + \partial_2^{3^\ast} \partial_2^3, \\
        \end{split}
        \]
        Since  $2$-simplex and $3$-simplex do not exist in $K_4$,  $\partial_2^{3}, \partial_2^{3+1}, $ and $ \partial_2^{3}$ do not exist, then
        \[
        \begin{split}
            \mathcal{L}_0^{3+1} &= \mathcal{B}_{1}^{3+1} \left(\mathcal{B}_{1}^{3+1}\right)^T + (\mathcal{B}_0^{3})^T \mathcal{B}_0^{3}, \\
            \mathcal{L}_1^{3+1} &=  (\mathcal{B}_1^{3})^T \mathcal{B}_1^3.
        \end{split}
        \]
        From \autoref{table:case 2}, one can see that $\beta_0^{3+1}=0 $ and $ \beta_1^{3+1}=1$, which reveals only one $0$-cycle and one $1$-cycle in $K_3$ are still alive in $K_4$.

        \begin{table}[H]
            \centering
            \setlength\tabcolsep{8pt}
            \captionsetup{margin=0.9cm}
            \caption{$K_3 \to K_4$}
            \begin{tabular}{c|ccc}
            \hline
            $q$  & $q=0$ & $q=1$ & $q=2$ \\ \hline \\
            $\mathcal{B}_{q+1}^{3+1}$  & $\begin{array}{@{}r@{}c@{}c@{}c@{}c@{}c@{}l@{}}
                    & 01 & 12 & 23 & 03 & 24  \\
                   \left.\begin{array}{c}
                    0 \\
                    1 \\
                    2 \\
                    3 \\
                    4
                    \end{array}\right[
                    & \begin{array}{c} -1 \\  1  \\  0  \\  0 \\ 0 \end{array}
                    & \begin{array}{c}  0 \\ -1  \\  1  \\  0 \\ 0 \end{array}
                    & \begin{array}{c}  0 \\  0  \\ -1  \\  1 \\ 0 \end{array}
                    & \begin{array}{c} -1 \\  0  \\  0  \\  1 \\ 0 \end{array}
                    & \begin{array}{c}  0 \\  0  \\ -1  \\  0 \\ 1 \end{array}
                    & \left]\begin{array}{c} \\ \\  \\ \\ \\
                    \end{array}\right.
                \end{array}$ & / &  /    \\  \\
            $\mathcal{B}_{q}^{3}$  &  $\begin{array}{@{}r@{}c@{}c@{}c@{}c@{}c@{}l@{}}
                    & 0 & 1 & 2 & 3 & 4   \\
                   \left.\begin{array}{c}

                    \end{array}\right[
                    & \begin{array}{c}  0  \end{array}
                    & \begin{array}{c}  0  \end{array}
                    & \begin{array}{c}  0  \end{array}
                    & \begin{array}{c}  0  \end{array}
                    & \begin{array}{c}  0  \end{array}
                    & \left]\begin{array}{c}
                    \end{array}\right.
                \end{array}$ & $\begin{array}{@{}r@{}c@{}c@{}c@{}c@{}l@{}}
                    & 01 & 12 & 23 & 03   \\
                   \left.\begin{array}{c}
                    0 \\
                    1 \\
                    2 \\
                    3 \\
                    4
                    \end{array}\right[
                    & \begin{array}{c} -1 \\  1  \\  0  \\  0 \\ 0 \end{array}
                    & \begin{array}{c}  0 \\ -1  \\  1  \\  0 \\ 0 \end{array}
                    & \begin{array}{c}  0 \\  0  \\ -1  \\  1 \\ 0 \end{array}
                    & \begin{array}{c} -1 \\  0  \\  0  \\  1 \\ 0 \end{array}
                    & \left]\begin{array}{c} \\ \\  \\ \\ \\
                    \end{array}\right.
                \end{array}$ &  /    \\  \\
            $\mathcal{L}_{q}^{3+1}$  & $\left[\begin{array}{ccccc}
                 2  &  -1  &   0  &  -1  &  0  \\
                -1  &   2  &  -1  &   0  &  0  \\
                 0  &  -1  &   3  &  -1  &  -1  \\
                -1  &   0  &  -1  &   2  &  0  \\
                 0  &   0  &  -1  &   0  &  1
            \end{array}\right]$  & $\left[\begin{array}{cccc}
                 2  &  -1  &   0  &   1   \\
                -1  &   2  &  -1  &   0   \\
                 0  &  -1  &   2  &   1   \\
                 1  &   0  &   1  &   2   \\

            \end{array}\right]$ &  /    \\  \\
            $\beta_{q}^{3+1}$                        & 1               & 1 & /    \\  \\
            $\text{dim}(\mathcal{L}_{q}^{3+1})$      & 5               & 4 & /    \\  \\
            $\text{rank}(\mathcal{L}_{q}^{3+1})$     & 4               & 3 & /    \\  \\
            $\text{nullity}(\mathcal{L}_{q}^{3+1})$  & 1               & 1 & /  \\  \\
            $\text{Spectra}(\mathcal{L}_{q}^{3+1})$    & $\{0,0.8299,2,2.6889,4.4812\}$   & $\{0,2,2,4\}$ & /  \\ \hline
            \end{tabular}
            \\
            \label{table:case 2}
        \end{table}

        \item[Case $3$.] The initial setup is $K_4$ and the end status is $K_4$. Similarly,
        \[
        \begin{split}
            \mathcal{L}_0^{4+0} &= \mathcal{B}_{1}^{4+0} \left(\mathcal{B}_{1}^{4+0}\right)^T + (\mathcal{B}_0^{4})^T \mathcal{B}_0^{4}, \\
            \mathcal{L}_1^{4+0} &=  (\mathcal{B}_1^{4})^T \mathcal{B}_1^4,
        \end{split}
        \]
        and $\mathcal{L}_2^{4+0}$ does not exist.
        In this case, the $0$-persistent $q$-combinatorial Laplacian matrix is actually the $q$-combinatorial Laplacian matrix defined in Eq. \eqref{equ:combinatorial Laplacian}. Therefore, $\beta_0^{4+0}, \beta_1^{4+0}, $ and $ \beta_2^{4+0}$ actually represent the number of $0,1,2$-cycles in $K_4$. With the filtration parameter $r$ increasing, all the circles overlapped with at least another circle in $K_4$, which results in $\beta_0^{4+0} = 1$. Since  only one $1$-cycle formed in $K_4$, one has $\beta_1^{4+0} = 1$.

        \begin{table}[H]
            \centering
            \setlength\tabcolsep{6pt}
            \captionsetup{margin=0.9cm}
            \caption{$K_4 \to K_4$}
            \begin{tabular}{c|ccc}
            \hline
            $q$  & $q=0$ & $q=1$ & $q=2$ \\ \hline \\
            $\mathcal{B}_{q+1}^{4+0}$  & $\begin{array}{@{}r@{}c@{}c@{}c@{}c@{}c@{}l@{}}
                    & 01 & 12 & 23 & 03 & 24  \\
                   \left.\begin{array}{c}
                    0 \\
                    1 \\
                    2 \\
                    3 \\
                    4
                    \end{array}\right[
                    & \begin{array}{c} -1 \\  1  \\  0  \\  0 \\ 0 \end{array}
                    & \begin{array}{c}  0 \\ -1  \\  1  \\  0 \\ 0 \end{array}
                    & \begin{array}{c}  0 \\  0  \\ -1  \\  1 \\ 0 \end{array}
                    & \begin{array}{c} -1 \\  0  \\  0  \\  1 \\ 0 \end{array}
                    & \begin{array}{c}  0 \\  0  \\ -1  \\  0 \\ 1 \end{array}
                    & \left]\begin{array}{c} \\ \\  \\ \\ \\
                    \end{array}\right.
                \end{array}$ & / &  /    \\  \\
            $\mathcal{B}_{q}^{4}$  &  $\begin{array}{@{}r@{}c@{}c@{}c@{}c@{}c@{}l@{}}
                    & 0 & 1 & 2 & 3 & 4   \\
                   \left.\begin{array}{c}
                    /
                    \end{array}\right[
                    & \begin{array}{c}  0  \end{array}
                    & \begin{array}{c}  0  \end{array}
                    & \begin{array}{c}  0  \end{array}
                    & \begin{array}{c}  0  \end{array}
                    & \begin{array}{c}  0  \end{array}
                    & \left]\begin{array}{c} \\
                    \end{array}\right.
                \end{array}$ & $\begin{array}{@{}r@{}c@{}c@{}c@{}c@{}c@{}l@{}}
                    & 01 & 12 & 23 & 03 & 24  \\
                   \left.\begin{array}{c}
                    0 \\
                    1 \\
                    2 \\
                    3 \\
                    4
                    \end{array}\right[
                    & \begin{array}{c} -1 \\  1  \\  0  \\  0 \\ 0 \end{array}
                    & \begin{array}{c}  0 \\ -1  \\  1  \\  0 \\ 0 \end{array}
                    & \begin{array}{c}  0 \\  0  \\ -1  \\  1 \\ 0 \end{array}
                    & \begin{array}{c} -1 \\  0  \\  0  \\  1 \\ 0 \end{array}
                    & \begin{array}{c}  0 \\  0  \\ -1  \\  0 \\ 1 \end{array}
                    & \left]\begin{array}{c} \\ \\  \\ \\ \\
                    \end{array}\right.
                \end{array}$ &  /    \\  \\
            $\mathcal{L}_{q}^{4+0}$  & $\left[\begin{array}{ccccc}
                 2  &  -1  &   0  &  -1  &  0  \\
                -1  &   2  &  -1  &   0  &  0  \\
                 0  &  -1  &   3  &  -1  &  -1  \\
                -1  &   0  &  -1  &   2  &  0  \\
                 0  &   0  &  -1  &   0  &  1
            \end{array}\right]$  & $\left[\begin{array}{ccccc}
                 2  &  -1  &   0  &   1  &  0  \\
                -1  &   2  &  -1  &   0  & -1  \\
                 0  &  -1  &   2  &   1  &  1  \\
                 1  &   0  &   1  &   2  &  0  \\
                 0  &  -1  &   1  &   0  &  2
            \end{array}\right]$ &  /    \\  \\
            $\beta_{q}^{4+0}$                        & 1               & 1 & /    \\  \\
            $\text{dim}(\mathcal{L}_{q}^{4+0})$      & 5               & 5 & /    \\  \\
            $\text{rank}(\mathcal{L}_{q}^{4+0})$     & 4               & 4 & /    \\  \\
            $\text{nullity}(\mathcal{L}_{q}^{4+0})$  & 1               & 1 & /  \\  \\
            $\text{Spectra}(\mathcal{L}_{q}^{4+0})$    & $\{0,0.8299,2,2.6889,4.4812\}$   & $\{0,0.8299,2,2.6889,4.4812\}$ & /  \\ \hline
            \end{tabular}
            \\
            \label{table:case 3}
        \end{table}

        \item[Case $4$.] The initial setup is $K_4$ and the end status is $K_5$. Using similar analysis as in previous cases, we have
        \[
        \begin{split}
            \mathcal{L}_0^{4+1} &= \mathcal{B}_{1}^{4+1} \left(\mathcal{B}_{1}^{4+1}\right)^T + (\mathcal{B}_0^{4})^T \mathcal{B}_0^{4}, \\
            \mathcal{L}_1^{4+1} &= \mathcal{B}_{2}^{4+1} \left(\mathcal{B}_{2}^{4+1}\right)^T + (\mathcal{B}_1^{4})^T \mathcal{B}_1^4,
        \end{split}
        \]
        and $\mathcal{L}_2^{4+1}$ does not exist. Notice that two $2$-simplices $012 $ and $ 023$ are created under the filtration process. The appearance of these two newborns results in the $1$-cycle that was alive in $K_4$ being killed. Therefore $\beta_1^{4+1} = 0$ and $\beta_0^{4+1} = 1$ because only one connected component keeps alive until $K_{5}$.
        \begin{table}[H]
            \centering
            \setlength\tabcolsep{2pt}
            \captionsetup{margin=0.9cm}
            \caption{$K_4 \to K_5$}
            \begin{tabular}{c|ccc}
            \hline
            $q$  & $q=0$ & $q=1$ & $q=2$ \\ \hline \\
            $\mathcal{B}_{q+1}^{4+1}$  & $\begin{array}{@{}r@{}c@{}c@{}c@{}c@{}c@{}c@{}l@{}}
                    & 01 & 12 & 23 & 03 & 24 & 02 \\
                   \left.\begin{array}{c}
                    0 \\
                    1 \\
                    2 \\
                    3 \\
                    4
                    \end{array}\right[
                    & \begin{array}{c} -1 \\  1  \\  0  \\  0 \\ 0 \end{array}
                    & \begin{array}{c}  0 \\ -1  \\  1  \\  0 \\ 0 \end{array}
                    & \begin{array}{c}  0 \\  0  \\ -1  \\  1 \\ 0 \end{array}
                    & \begin{array}{c} -1 \\  0  \\  0  \\  1 \\ 0 \end{array}
                    & \begin{array}{c}  0 \\  0  \\ -1  \\  0 \\ 1 \end{array}
                    & \begin{array}{c} -1 \\  0  \\  1  \\  0 \\ 0 \end{array}
                    & \left]\begin{array}{c} \\ \\  \\ \\ \\
                    \end{array}\right.
                \end{array}$ & $\begin{array}{@{}r@{}c@{}c@{}l@{}}
                    & 012 & 023  \\
                   \left.\begin{array}{c}
                    01 \\
                    12 \\
                    23 \\
                    03 \\
                    24
                    \end{array}\right[
                    & \begin{array}{c}  1 \\  1  \\  0  \\  0 \\ 0  \end{array}
                    & \begin{array}{c}  0 \\  0  \\  1  \\ -1 \\ 0  \end{array}
                    & \left]\begin{array}{c} \\ \\  \\ \\ \\
                    \end{array}\right.
                \end{array}$ &  /    \\  \\
            $\mathcal{B}_{q}^{4}$  &  $\begin{array}{@{}r@{}c@{}c@{}c@{}c@{}c@{}l@{}}
                    & 0 & 1 & 2 & 3 & 4   \\
                   \left.\begin{array}{c}
                    /
                    \end{array}\right[
                    & \begin{array}{c}  0  \end{array}
                    & \begin{array}{c}  0  \end{array}
                    & \begin{array}{c}  0  \end{array}
                    & \begin{array}{c}  0  \end{array}
                    & \begin{array}{c}  0  \end{array}
                    & \left]\begin{array}{c} \\
                    \end{array}\right.
                \end{array}$ & $\begin{array}{@{}r@{}c@{}c@{}c@{}c@{}c@{}l@{}}
                    & 01 & 12 & 23 & 03 & 24 \\
                   \left.\begin{array}{c}
                    0 \\
                    1 \\
                    2 \\
                    3 \\
                    4
                    \end{array}\right[
                    & \begin{array}{c} -1 \\  1  \\  0  \\  0 \\ 0  \end{array}
                    & \begin{array}{c}  0 \\ -1  \\  1  \\  0 \\ 0  \end{array}
                    & \begin{array}{c}  0 \\  0  \\ -1  \\  1 \\ 0  \end{array}
                    & \begin{array}{c} -1 \\  0  \\  0  \\  1 \\ 0  \end{array}
                    & \begin{array}{c}  0 \\  0  \\ -1  \\  0 \\ 1  \end{array}
                    & \left]\begin{array}{c} \\ \\  \\ \\ \\
                    \end{array}\right.
                \end{array}$ &  /    \\  \\
            $\mathcal{L}_{q}^{4+1}$  & $\left[\begin{array}{ccccc}
                 3  &  -1  &  -1  &  -1  &  0  \\
                -1  &   2  &  -1  &   0  &  0  \\
                -1  &  -1  &   4  &  -1  &  -1  \\
                -1  &   0  &  -1  &   2  &  0  \\
                 0  &   0  &  -1  &   0  &  1
            \end{array}\right]$  & $\left[\begin{array}{ccccc}
                 3  &   0  &   0  &   1  &  0  \\
                 0  &   3  &  -1  &   0  & -1  \\
                 0  &  -1  &   3  &   0  &  1  \\
                 1  &   0  &   0  &   3  &  0  \\
                 0  &  -1  &   1  &   0  &  2
            \end{array}\right]$ &  /    \\  \\
            $\beta_{q}^{4+1}$                        & 1               & 0 & /    \\  \\
            $\text{dim}(\mathcal{L}_{q}^{4+1})$      & 5               & 5 & /    \\  \\
            $\text{rank}(\mathcal{L}_{q}^{4+1})$     & 4               & 5 & /    \\  \\
            $\text{nullity}(\mathcal{L}_{q}^{4+1})$  & 1               & 0 & /  \\  \\
            $\text{Spectra}(\mathcal{L}_{q}^{4+1})$    & $\{0,1,2,4,5\}$   & $\{1.2677,2,2,4,4.7321\}$ & /  \\ \hline
            \end{tabular}
            \\
            \label{table:case 4}
        \end{table}

        \item[Case $5$.] The initial setup is $K_5$ and the end status is $K_6$. The $1$-persistent $0,1,2$-combinatorial Laplacian matrices are
        \[
        \begin{split}
            \mathcal{L}_0^{5+1} &= \mathcal{B}_{1}^{5+1} \left(\mathcal{B}_{1}^{5+1}\right)^T+ (\mathcal{B}_0^{5})^T \mathcal{B}_0^{5}, \\
            \mathcal{L}_1^{5+1} &= \mathcal{B}_{2}^{5+1} \left(\mathcal{B}_{2}^{5+1}\right)^T + (\mathcal{B}_1^{5})^T \mathcal{B}_1^5,  \\
            \mathcal{L}_2^{5+1} &= \mathcal{B}_{3}^{5+1} \left(\mathcal{B}_{3}^{5+1}\right)^T + (\mathcal{B}_2^{5})^T \mathcal{B}_2^5.
        \end{split}
        \]
    In this situation, a new $3$-simplex is formed in $K_6$. From \autoref{table:case 5}, we can see that $\beta_2^{5+1} = 0$ because  $K_5$ does not own any $2$-cycle and thus, there is no $2$-cycle keeping alive up to $K_6$. $\beta_0^{5+1}$ implies only one $0$-cycle preserved along the filtration process.

        \begin{table}[H]
            \centering
            \setlength\tabcolsep{0pt}
            \captionsetup{margin=0.9cm}
            \caption{$K_5 \to K_6$}
            \begin{tabular}{c|ccc}
            \hline
            $q$  & $q=0$ & $q=1$ & $q=2$ \\ \hline \\
            $\mathcal{B}_{q+1}^{5+1}$  & $\begin{array}{@{}r@{}c@{}c@{}c@{}c@{}c@{}c@{}l@{}}
                    & 01 & 12 & 23 & 03 & 24 & 02 \\
                   \left.\begin{array}{c}
                    0 \\
                    1 \\
                    2 \\
                    3 \\
                    4
                    \end{array}\right[
                    & \begin{array}{c} -1 \\  1  \\  0  \\  0 \\ 0 \end{array}
                    & \begin{array}{c}  0 \\ -1  \\  1  \\  0 \\ 0 \end{array}
                    & \begin{array}{c}  0 \\  0  \\ -1  \\  1 \\ 0 \end{array}
                    & \begin{array}{c} -1 \\  0  \\  0  \\  1 \\ 0 \end{array}
                    & \begin{array}{c}  0 \\  0  \\ -1  \\  0 \\ 1 \end{array}
                    & \begin{array}{c} -1 \\  0  \\  1  \\  0 \\ 0 \end{array}
                    & \left]\begin{array}{c} \\ \\  \\ \\ \\
                    \end{array}\right.
                \end{array}$ & $\begin{array}{@{}r@{}c@{}c@{}c@{}c@{}l@{}}
                    & 012 & 023 & 013 & 123  \\
                   \left.\begin{array}{c}
                    01 \\
                    12 \\
                    23 \\
                    03 \\
                    24 \\
                    02
                    \end{array}\right[
                    & \begin{array}{c}  1 \\  1  \\  0  \\  0 \\ 0 \\ -1  \end{array}
                    & \begin{array}{c}  0 \\  0  \\  1  \\ -1 \\ 0 \\  1  \end{array}
                    & \begin{array}{c}  1 \\  0  \\  0  \\ -1 \\ 0 \\  0  \end{array}
                    & \begin{array}{c}  0 \\  1  \\  1  \\  0 \\ 0 \\  0  \end{array}
                    & \left]\begin{array}{c} \\ \\  \\ \\ \\ \\
                    \end{array}\right.
                \end{array}$ &  /   \\  \\
            $\mathcal{B}_{q}^{5}$  &  $\begin{array}{@{}r@{}c@{}c@{}c@{}c@{}c@{}l@{}}
                    & 0 & 1 & 2 & 3 & 4   \\
                   \left.\begin{array}{c}
                    /
                    \end{array}\right[
                    & \begin{array}{c}  0  \end{array}
                    & \begin{array}{c}  0  \end{array}
                    & \begin{array}{c}  0  \end{array}
                    & \begin{array}{c}  0  \end{array}
                    & \begin{array}{c}  0  \end{array}
                    & \left]\begin{array}{c} \\
                    \end{array}\right.
                \end{array}$ & $\begin{array}{@{}r@{}c@{}c@{}c@{}c@{}c@{}c@{}l@{}}
                    & 01 & 12 & 23 & 03 & 24 & 02  \\
                   \left.\begin{array}{c}
                    0 \\
                    1 \\
                    2 \\
                    3 \\
                    4
                    \end{array}\right[
                    & \begin{array}{c} -1 \\  1  \\  0  \\  0 \\ 0 \end{array}
                    & \begin{array}{c}  0 \\ -1  \\  1  \\  0 \\ 0 \end{array}
                    & \begin{array}{c}  0 \\  0  \\ -1  \\  1 \\ 0 \end{array}
                    & \begin{array}{c} -1 \\  0  \\  0  \\  1 \\ 0 \end{array}
                    & \begin{array}{c}  0 \\  0  \\ -1  \\  0 \\ 1 \end{array}
                    & \begin{array}{c} -1 \\  0  \\  1  \\  0 \\ 0 \end{array}
                    & \left]\begin{array}{c} \\ \\  \\ \\ \\
                    \end{array}\right.
                \end{array}$ &  $\begin{array}{@{}r@{}c@{}c@{}l@{}}
                    & 012 & 023   \\
                   \left.\begin{array}{c}
                    01 \\
                    12 \\
                    23 \\
                    03 \\
                    24 \\
                    02
                    \end{array}\right[
                    & \begin{array}{c}  1 \\  1  \\  0  \\  0 \\ 0 \\ -1 \end{array}
                    & \begin{array}{c}  0 \\  0  \\  1  \\ -1 \\ 0 \\  1 \end{array}
                    & \left]\begin{array}{c} \\ \\  \\ \\ \\ \\
                    \end{array}\right.
                \end{array}$    \\  \\
            $\mathcal{L}_{q}^{5+1}$  & $\left[\begin{array}{ccccc}
                 3  &  -1  &  -1  &  -1  &  0  \\
                -1  &   2  &  -1  &   0  &  0  \\
                -1  &  -1  &   4  &  -1  &  -1  \\
                -1  &   0  &  -1  &   2  &  0  \\
                 0  &   0  &  -1  &   0  &  1
            \end{array}\right]$  & $\left[\begin{array}{cccccc}
                 4  &   0  &   0  &   0  &  0  &  0 \\
                 0  &   4  &   0  &   0  & -1  &  0 \\
                 0  &   0  &   4  &   0  &  1  &  0 \\
                 0  &   0  &   0  &   4  &  0  &  0 \\
                 0  &  -1  &   1  &   0  &  2  & -1 \\
                 0  &   0  &   0  &   0  & -1  &  4
            \end{array}\right]$ &  $\left[\begin{array}{cc}
                 3  &  -1  \\
                -1  &   3
            \end{array}\right]$    \\  \\
            $\beta_{q}^{5+1}$                        & 1               & 0 & 0    \\  \\
            $\text{dim}(\mathcal{L}_{q}^{5+1})$      & 5               & 6 & 2    \\  \\
            $\text{rank}(\mathcal{L}_{q}^{5+1})$     & 4               & 6 & 2    \\  \\
            $\text{nullity}(\mathcal{L}_{q}^{5+1})$  & 1               & 0 & 0     \\  \\
            $\text{Spectra}(\mathcal{L}_{q}^{5+1})$    & $\{0,1,2,4,5\}$   & $\{1,4,4,4,4,5\}$ & $\{2,4\}$ \\ \hline
            \end{tabular}
            \\
            \label{table:case 5}
        \end{table}

        \item[Case $6$.] The initial setup is $K_6$ and the end status is $K_6$.
        The $0$-persistent $0,1,2$-combinatorial Laplacian operators are
        \[
        \begin{split}
            \mathcal{L}_0^{6+0} &= \mathcal{B}_{1}^{6+0} (\mathcal{B}_{1}^{6+0})^T + (\mathcal{B}_0^{6})^T \mathcal{B}_0^{6}, \\
            \mathcal{L}_1^{6+0} &= \mathcal{B}_{2}^{6+0} (\mathcal{B}_{2}^{6+0})^T + (\mathcal{B}_{1}^{6})^T \mathcal{B}_{1}^{6},\\
            \mathcal{L}_2^{6+0} &= \mathcal{B}_{3}^{6+0} (\mathcal{B}_{3}^{6+0})^T + (\mathcal{B}_{2}^{6})^T \mathcal{B}_{2}^{6},
        \end{split}
        \]
        $\beta_0^{6+0}=1, \beta_1^{6+0}=0, $ and $ \beta_2^{6+0}=0$ imply that only one $0$-cycle (connected component) exists in $K_6$.

        \begin{table}[H]
            \centering
            \setlength\tabcolsep{0pt}
            \captionsetup{margin=0.9cm}
            \caption{$K_6 \to K_6$}
            \begin{tabular}{c|ccc}
            \hline
            $q$  & $q=0$ & $q=1$ & $q=2$ \\ \hline \\
            $\mathcal{B}_{q+1}^{6+0}$  & $\begin{array}{@{}r@{}c@{}c@{}c@{}c@{}c@{}c@{}c@{}l@{}}
                    & 01 & 12 & 23 & 03 & 24 & 02 & 13 \\
                   \left.\begin{array}{c}
                    0 \\
                    1 \\
                    2 \\
                    3 \\
                    4
                    \end{array}\right[
                    & \begin{array}{c} -1 \\  1  \\  0  \\  0 \\ 0 \end{array}
                    & \begin{array}{c}  0 \\ -1  \\  1  \\  0 \\ 0 \end{array}
                    & \begin{array}{c}  0 \\  0  \\ -1  \\  1 \\ 0 \end{array}
                    & \begin{array}{c} -1 \\  0  \\  0  \\  1 \\ 0 \end{array}
                    & \begin{array}{c}  0 \\  0  \\ -1  \\  0 \\ 1 \end{array}
                    & \begin{array}{c} -1 \\  0  \\  1  \\  0 \\ 0 \end{array}
                    & \begin{array}{c}  0 \\ -1  \\  0  \\  1 \\ 0 \end{array}
                    & \left]\begin{array}{c} \\ \\  \\ \\  \\
                    \end{array}\right.
                \end{array}$ & $\begin{array}{@{}r@{}c@{}c@{}c@{}c@{}l@{}}
                    & 012 & 023 & 013 & 123  \\
                   \left.\begin{array}{c}
                    01 \\
                    12 \\
                    23 \\
                    03 \\
                    24 \\
                    02 \\
                    13
                    \end{array}\right[
                    & \begin{array}{c}  1 \\  1  \\  0  \\  0 \\ 0 \\ -1  \\  0  \end{array}
                    & \begin{array}{c}  0 \\  0  \\  1  \\ -1 \\ 0 \\  1  \\  0  \end{array}
                    & \begin{array}{c}  1 \\  0  \\  0  \\ -1 \\ 0 \\  0  \\  1  \end{array}
                    & \begin{array}{c}  0 \\  1  \\  1  \\  0 \\ 0 \\  0  \\ -1  \end{array}
                    & \left]\begin{array}{c} \\ \\  \\ \\ \\ \\ \\
                    \end{array}\right.
                \end{array}$ &  $\mathcal{B}_3^{6+0}$    \\  \\
            $\mathcal{B}_{q}^{6}$  &  $\begin{array}{@{}r@{}c@{}c@{}c@{}c@{}c@{}l@{}}
                    & 0 & 1 & 2 & 3 & 4   \\
                   \left.\begin{array}{c}
                    /
                    \end{array}\right[
                    & \begin{array}{c}  0  \end{array}
                    & \begin{array}{c}  0  \end{array}
                    & \begin{array}{c}  0  \end{array}
                    & \begin{array}{c}  0  \end{array}
                    & \begin{array}{c}  0  \end{array}
                    & \left]\begin{array}{c} \\
                    \end{array}\right.
                \end{array}$ & $\begin{array}{@{}r@{}c@{}c@{}c@{}c@{}c@{}c@{}c@{}l@{}}
                    & 01 & 12 & 23 & 03 & 24 & 02 & 13 \\
                   \left.\begin{array}{c}
                    0 \\
                    1 \\
                    2 \\
                    3 \\
                    4
                    \end{array}\right[
                    & \begin{array}{c} -1 \\  1  \\  0  \\  0 \\ 0 \end{array}
                    & \begin{array}{c}  0 \\ -1  \\  1  \\  0 \\ 0 \end{array}
                    & \begin{array}{c}  0 \\  0  \\ -1  \\  1 \\ 0 \end{array}
                    & \begin{array}{c} -1 \\  0  \\  0  \\  1 \\ 0 \end{array}
                    & \begin{array}{c}  0 \\  0  \\ -1  \\  0 \\ 1 \end{array}
                    & \begin{array}{c} -1 \\  0  \\  1  \\  0 \\ 0 \end{array}
                    & \begin{array}{c}  0 \\ -1  \\  0  \\  1 \\ 0 \end{array}
                    & \left]\begin{array}{c} \\ \\  \\ \\ \\
                    \end{array}\right.
                \end{array}$ &   $\mathcal{B}_2^{6}$   \\  \\
            $\mathcal{L}_{q}^{6+0}$  & $\left[\begin{array}{ccccc}
                 3  &  -1  &  -1  &  -1  &  0  \\
                -1  &   2  &  -1  &   0  &  0  \\
                -1  &  -1  &   4  &  -1  &  -1  \\
                -1  &   0  &  -1  &   2  &  0  \\
                 0  &   0  &  -1  &   0  &  1
            \end{array}\right]$  & $\left[\begin{array}{ccccccc}
                 4  &   0  &   0  &   0  &  0  &  0  &  0 \\
                 0  &   4  &   0  &   0  & -1  &  0  &  0 \\
                 0  &   0  &   4  &   0  &  1  &  0  &  0 \\
                 0  &   0  &   0  &   4  &  0  &  0  &  0 \\
                 0  &  -1  &   1  &   0  &  2  & -1  &  0 \\
                 0  &   0  &   0  &   0  & -1  &  4  &  0 \\
                 0  &   0  &   0  &   0  &  0  &  0  &  4
            \end{array}\right]$ &  $\mathcal{L}_3^{6+0}$    \\  \\
            $\beta_{q}^{6+0}$                        & 1               & 0 & 0    \\  \\
            $\text{dim}(\mathcal{L}_{q}^{6+0})$      & 5               & 7 & 4    \\  \\
            $\text{rank}(\mathcal{L}_{q}^{6+0})$     & 4               & 7 & 4    \\  \\
            $\text{nullity}(\mathcal{L}_{q}^{6+0})$  & 1               & 0 & 0     \\  \\
            $\text{Spectra}(\mathcal{L}_{q}^{6+0})$    & $\{0,1,4,4,5\}$   & $\{1,4,4,4,4,4,5\}$ & $\{4,4,4,4\}$ \\ \hline
            \end{tabular}
            \\
            \label{table:case 6}
        \end{table}
        with
        \[
        \mathcal{B}_3^{6+0} =
        \begin{array}{@{}r@{}c@{}l@{}}
                    & 0123  \\
                   \left.\begin{array}{c}
                    012 \\
                    023 \\
                    013 \\
                    123
                    \end{array}\right[
                    & \begin{array}{c}  -1 \\  -1 \\ 1 \\ 1  \end{array}
                    & \left]\begin{array}{c} \\ \\  \\ \\
                    \end{array}\right.
        \end{array},
        \mathcal{B}_2^6 =
        \begin{array}{@{}r@{}c@{}c@{}c@{}c@{}l@{}}
                    & 012 & 023 & 013 & 123  \\
                   \left.\begin{array}{c}
                    01 \\
                    12 \\
                    23 \\
                    03 \\
                    24 \\
                    02 \\
                    13
                    \end{array}\right[
                    & \begin{array}{c}  1 \\  1  \\  0  \\  0 \\ 0 \\ -1  \\  0  \end{array}
                    & \begin{array}{c}  0 \\  0  \\  1  \\ -1 \\ 0 \\  1  \\  0  \end{array}
                    & \begin{array}{c}  1 \\  0  \\  0  \\ -1 \\ 0 \\  0  \\  1  \end{array}
                    & \begin{array}{c}  0 \\  1  \\  1  \\  0 \\ 0 \\  0  \\ -1  \end{array}
                    & \left]\begin{array}{c} \\ \\  \\ \\ \\ \\ \\
                    \end{array}\right.
        \end{array}
        \]
        and
        \[
        \mathcal{L}_3^{6+0} =
        \left[\begin{array}{cccc}
                 4 & 0 & 0 & 0 \\
                 0 & 4 & 0 & 0 \\
                 0 & 0 & 4 & 0 \\
                 0 & 0 & 0 & 4
        \end{array}\right].
        \]
    \end{enumerate}

		We have constructed a family of persistent spectral graphs induced by a filtration parameter. For the sake of simplicity, we focus on the analysis of high-dimensional spectra with $p=0$ in the rest of this section.   As clarified before, the $0$-persistent $q$-combinatorial Laplacian matrix is the $q$-combinatorial Laplacian matrix.

A  graph structure encodes inter-dependencies among constituents and provides a convenient representation of the high-dimensional data. Naturally, the same idea can be applied to higher-dimensional spaces. For a set of points $V\subset \mathbb{R}^{n}$ without additional structures, we consider growing an $(n-1)$-sphere centered at each point with an ever-increasing radius $r$. Therefore, a family of $0$-persistent $q$-combinatorial Laplacian matrices as well as spectra can be generated as the radius $r$ increases, which provides topological and spectral features to distinguish individual entries of the dataset. 

\begin{figure}[H]
        \centering
        \includegraphics[width=16cm]{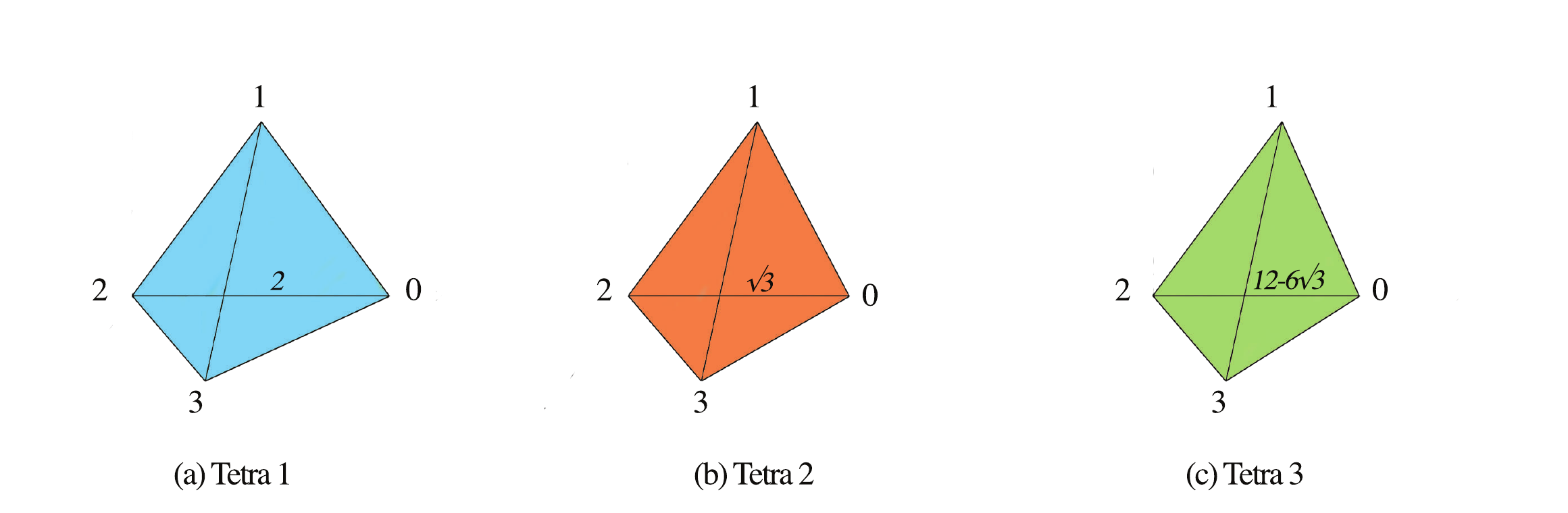}
        \captionsetup{margin=0.9cm}  
        \caption{Three tetrahedrons with different topological shape in ${\mathbb R}^3$. $(a)$ Regular tetrahedron with edge length $2$. $(b)$ Move $v_0$ along the edge $[v_0, v_1]$ and construct a new tetrahedron with the length of $[v_0, v_1]$ to be $\sqrt{3}$. $(c)$ Move $v_0$ along the edge $[v_0, v_1]$ and construct a new tetrahedron with the length of $[v_0, v_1]$ being $12-6\sqrt{3}$.}
        \label{fig: tetrahedron}
    \end{figure}

\begin{figure}[H]
        \centering
        \captionsetup{margin=0.9cm}
        \subfigure[]{
            \label{fig:subfig:tetra zero}
            \includegraphics[scale=0.45]{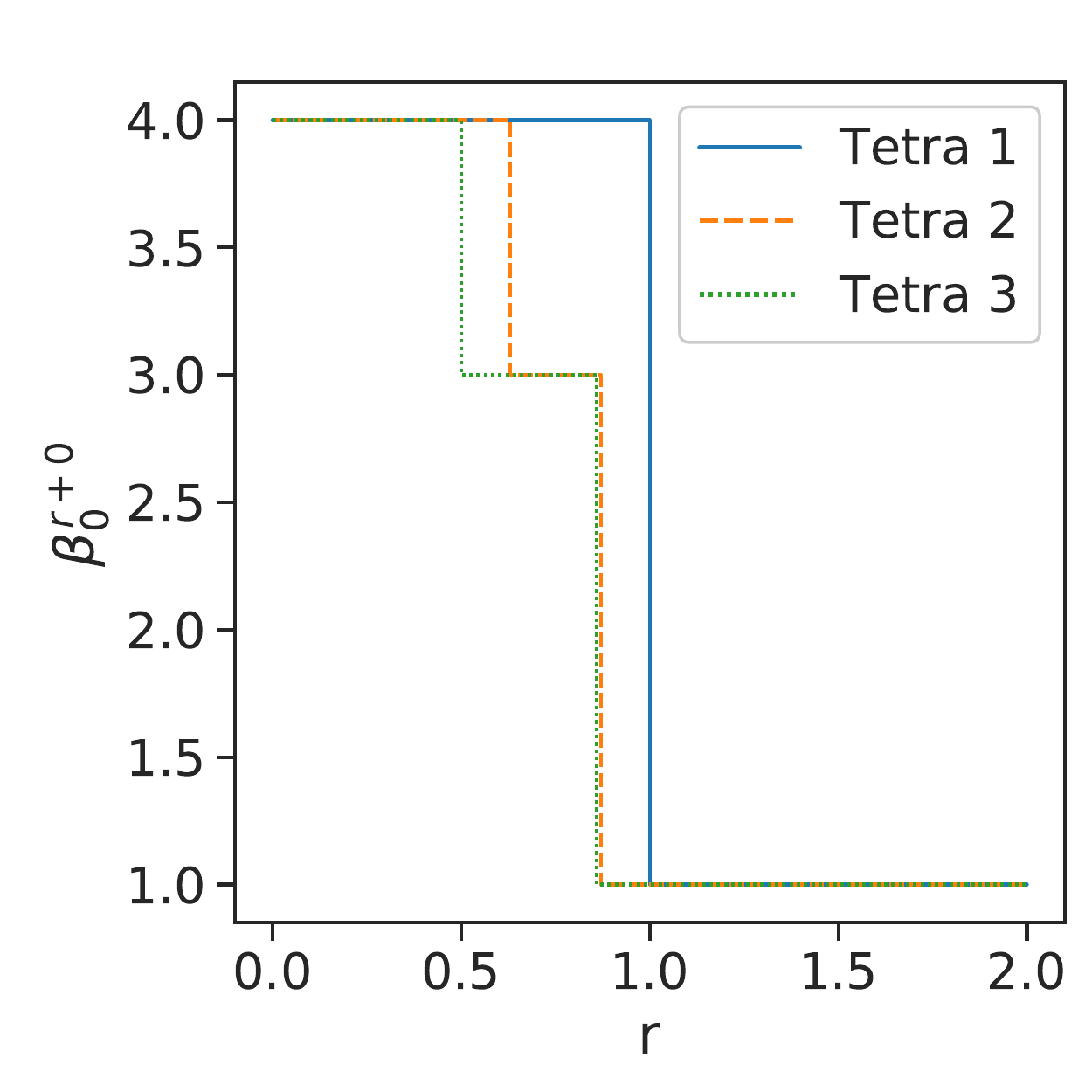}}
            \hspace{0.05in}
        \subfigure[]{
            \label{fig:subfig:tetra sec}
            \includegraphics[scale = 0.45]{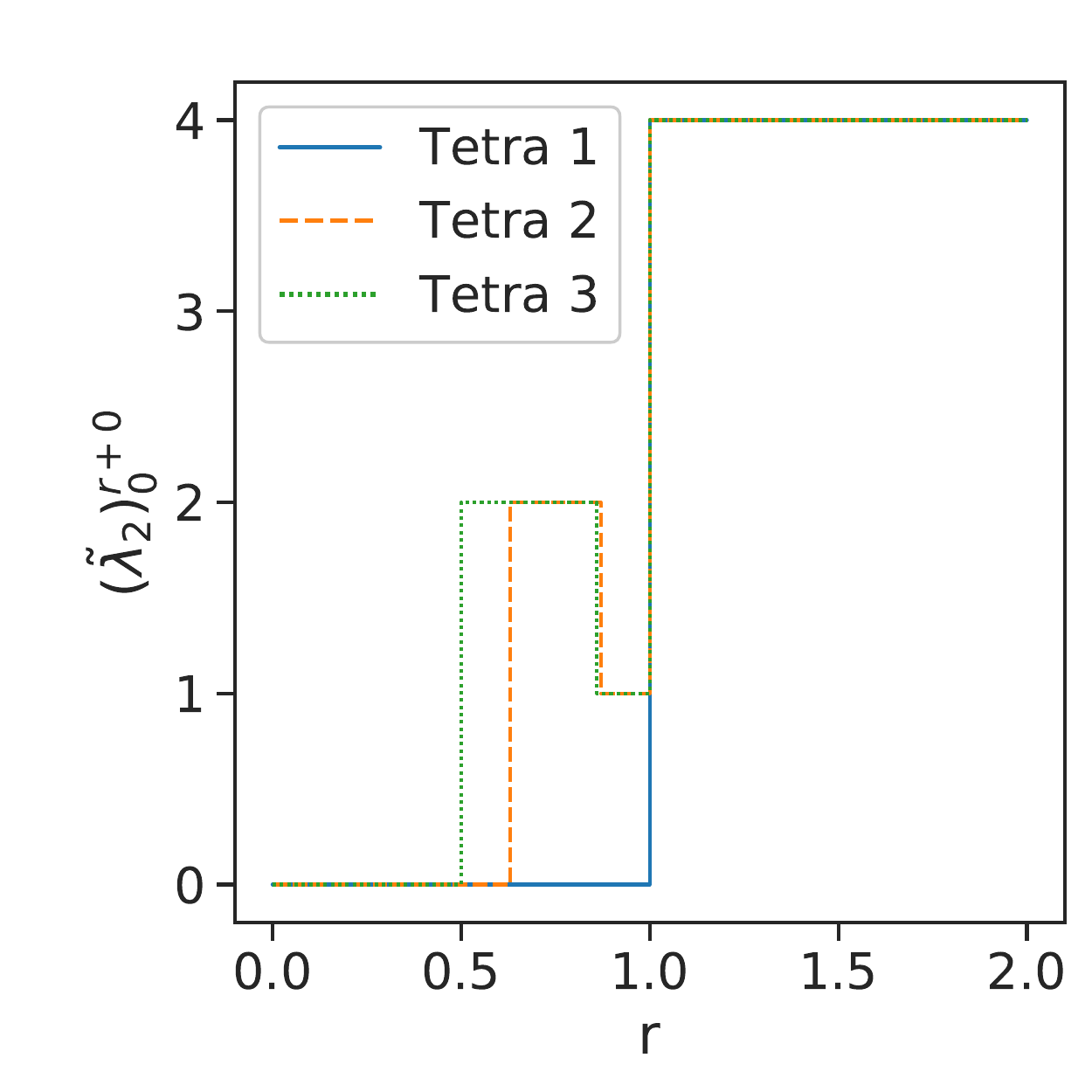}}
        \caption{(a) Plot of $\beta_0^{r+0}$ with radius filtration $r$ among $3$ different tetrahedrons. (b) Plot of $(\tilde{\lambda}_2)^{r+0}_0$ with radius filtration $r$ among $3$ different tetrahedrons.}
    \label{fig:Tetra filtration}
    \end{figure}

\autoref{fig: tetrahedron} and \autoref{fig:Tetra filtration} exemplify the capacity of persistent spectral theory to discriminate between different structures in $\mathbb{R}^{3}$. In Figure \autoref{fig:subfig:tetra zero}, we employ the persistent spectral analysis based on the $\beta_0^{r+0}$ tendency along the filtration to distinguish three tetrahedrons. As $r$ grows, isolated points ($0$-simplices) will gradually grow into solid 2-spheres, and a new isolated component will be created once two spheres corresponding to two isolated points overlap with each other. Since $\beta_0^{r+0}$ represents the number of isolated components, the value of $\beta_0^{r+0}$ will finally decrease to $1$. Take Tetra $2$ as an example. It is seen that at the initial setup ($r=0.0$), the number of isolated components is $4$, which represents the number of isolated points. When $r$ is around $0.63$, two spheres centered at $v_0$ and $v_2$ with radius $0.63$ overlapped with each other. Therefore, $\beta_0^{r+0}$ reduces to $3$ at this point. With $r$ keeping growing, the sphere centered at $v_0$ overlaps with spheres centered at $v_1, v_2,$ and $v_3$, which results $\beta_0^{r+0}=1$ after $r=0.87$. Similarly, the smallest non-zero eigenvalue $(\tilde{\lambda}_2)_0^{r+0}$ changes at radius $0.63$ and $0.87$ in Figure \autoref{fig:subfig:tetra sec}, which also affirms that the solid spheres get overlapped at these specific filtration parameters. It is clear  that Tetrahedron $1, 2,$ and $3$ have different $\beta_0^{r+0}$ and $(\tilde{\lambda}_2)_0^{r+0}$ values. Since $1$-cycle and $2$-cycle are not formed along with the filtration, analysis of $\beta_1^{r+0}$ and $\beta_2^{r+0}$ will not be mentioned in this case.

    \subsubsection{Variants of $p$-persistent $q$-combinatorial Laplacian matrices}

    The traditional approach in defining the $q$-boundary operator $\partial_q: C_q(K) \to C_{q-1}(K)$ can be expressed as:
    \begin{equation}\nonumber
        \partial_q \sigma_q = \sum_{i=0}^{q}(-1)^i \sigma^{i}_{q-1},
    \end{equation}
    which leads to the corresponding elements in the boundary matrices being either $1$ or $-1$. However, to encode more geometric information into the Laplacian operator, we  add volume information of $q$-simplex $\sigma_q$ to the expression of $q$-boundary operator.

    Given a vertex set $V = \{v_0, v_1, \cdots, v_q\}$ with $q+1$ isolated points ($0$-simplices) randomly arranged in the $n$-dimensional Euclidean space $\mathbb{R}^n$, often with $n \ge q$. Set $d_{ij}$ to be the distances between $v_i$ and $v_j$ with $0 \le i \le j \le q$ and obviously, $d_{ij} = d_{ji}$.  The Cayley-Menger determinant  can be expressed as  \cite{berger2009geometry}
    \begin{equation}
        {\rm Det}_{\rm CM}(v_0, v_1, \cdots, v_q) =
        \left|\begin{array}{cccccc}
                 0        & d_{01}^2 & d_{02}^2 & \cdots & d_{0q}^2 & 1 \\
                 d_{10}^2 & 0        & d_{12}^2 & \cdots & d_{1q}^2 & 1 \\
                 d_{20}^2 & d_{21}^2 & 0        & \cdots & d_{2q}^2 & 1 \\
                 \vdots   & \vdots   & \vdots   & \vdots & \ddots   & \vdots \\
                 d_{q0}^2 & d_{q1}^2 & d_{q2}^2 & \cdots & 0        & 1 \\
                 1        & 1        & 1        & 1      & 1        & 0
        \end{array}\right|
    \end{equation}
    The ${q}$-dimensional volume of $q$-simplex $\sigma_q$ with vertices $\{v_0, v_1, \cdots, v_q\}$ is defined by
    \begin{equation}
        \text{Vol}(\sigma_q) = \sqrt{\dfrac{(-1)^{q+1}}{(q!)^2 2^q}{\rm Det}_{\rm CM}(v_0, v_1, \cdots, v_q)}.
    \end{equation}
    In trivial cases, $\text{Vol}(\sigma_0) = 1$, meaning the $0$-dimensional volume of $0$-simplex is $1$, i.e., there is only $1$ vertex in a $0$-simplex. Also, the $1$-dimensional volume of $1$-simplex $\sigma_1 = [v_i,v_j]$ is the distance between $v_i$ and $v_j$, and the $2$-dimensional volume of $2$-simplex is the area of a triangle $[v_i, v_j, v_k]$.

    The weighted boundary operator  equipped with volume, denoted $\hat{\partial}_q$, is given by
    \begin{equation}
        \hat{\partial}_q \sigma_q = \sum_{i=0}^{q}(-1)^i \text{Vol}(\sigma^{i}_{q}) \sigma^{i}_{q-1}.
    \end{equation}
    Employed the same concept to the persistent spectral theory, we have the volume-weighted $p$-persistent $q$-combinatorial Laplacian operator. We  also define
    \begin{equation}\label{equ: dis weighted boundary map}
        \hat{\eth}_q^{t+p} : \mathbb{C}_q^{t+p} \to  C_{q-1}^{t}
    \end{equation}
    with
    \[
        \mathbb{C}_q^{t+p} \coloneqq \{ \alpha \in C_q^{t+p} \ | \ \hat{\partial}_q^{t+p}(\alpha) \in C_{q-1}^{t}\}.
    \]
    Similarly, an inverse-volume weighted  boundary operator, denoted $\check{\partial}_q$, is given by
    \begin{equation}
        \check{\partial}_q \sigma_q = \sum_{i=0}^{q}(-1)^i \frac{1}{\text{Vol}(\sigma^{i}_{q})} \sigma^{i}_{q-1}.
    \end{equation}
  To define an inverse-volume weighted $p$-persistent $q$-combinatorial Laplacian operator. We  define
  \begin{equation}\label{equ: inv weighted boundary map}
        \check{\eth}_q^{t+p} : \mathbb{C}_q^{t+p} \to  C_{q-1}^{t}
  \end{equation}
    with
    \[
        \mathbb{C}_q^{t+p} \coloneqq \{ \alpha \in C_q^{t+p} \ | \ \check{\partial}_q^{t+p}(\alpha) \in C_{q-1}^{t}\}.
    \]
       Then volume-weighted and inverse-volume-weighted  ${p}$-persistent $\boldsymbol{q}$-combinatorial Laplacian operators defined along the filtration can be expressed as
    \begin{equation}
        \begin{split}
         \hat{\Delta}_q^{t+p} &= \hat{\eth}_{q+1}^{t+p} \left( \hat{\eth}_{q+1}^{t+p}\right)^\ast + \hat{\partial}_q^{t^\ast} \hat{\partial}_q^t ,\\
         \check{\Delta}_q^{t+p} &= \check{\eth}_{q+1}^{t+p} \left( \check{\eth}_{q+1}^{t+p}\right)^\ast + \check{\partial}_q^{t^\ast} \check{\partial}_q^t.
        \end{split}
    \end{equation}
    The corresponding weighted matrix representations of boundary operators $\hat{\eth}_{q+1}^{t+p}$, $\hat{\eth}_q^t$, $\check{\eth}_{q+1}^{t+p}$, and $\check{\eth}_q^t$ are denoted $\hat{\mathcal{B}}_{q+1}^{t+p}$, $\hat{\mathcal{B}}_{q}^t$, $\check{\mathcal{B}}_{q+1}^{t+p}$, and $\check{\mathcal{B}}_{q}^t$, respectively. Therefore, volume-weighted and inverse-volume-weighted ${p}$-persistent $\boldsymbol{q}$-combinatorial Laplacian matrices  can be expressed as
    \begin{equation}
    \begin{split}
        \hat{\mathcal{L}}_q^{t+p} &= \hat{\mathcal{B}}_{q+1}^{t+p} (\hat{\mathcal{B}}_{q+1}^{t+p})^T + (\hat{\mathcal{B}}_{q}^t)^T (\hat{\mathcal{B}}_{q}^t), \\
        \check{\mathcal{L}}_q^{t+p} &= \check{\mathcal{B}}_{q+1}^{t+p} (\check{\mathcal{B}}_{q+1}^{t+p})^T + (\check{\mathcal{B}}_{q}^t)^T (\check{\mathcal{B}}_{q}^t).
    \end{split}
    \end{equation}
     Although the expressions of the weighted persistent Laplacian matrices are different from the original persistent Laplacian matrices, some properties of $\mathcal{L}_q^{t+p}$ are preserved. The weighted persistent Laplacian operators are still symmetric and positive semi-defined. Additionally, their ranks are the same as $\mathcal{L}_q^{t+p}$. With the embedded volume information, weighted PSGs can provide richer topological and geometric information through the associated persistent Betti numbers and non-harmonic spectra  (i.e., non-zero eigenvalues).

    In real applications, we are more interested in the $0,1,2$-combinatorial Laplacian matrices because its more intuitive to depict the relation among vertex, edges, and faces. Given a set of vertices  $V = \{v_0, v_2, \cdots, v_N\}$ with $N+1$ isolated points ($0$-simplices) randomly arranged in $\mathbb{R}^n$. By varying the radius $r$ of the $(n-1)$-sphere centered at each vertex, a variety of simplicial complexes is created. We denote the simplicial complex generated at radius $r$ to be $K_r$, then the $0$-persistent $q$-combinatorial Laplacian operator and matrix at initial set up $K_r$ is
    \begin{equation}
        \mathcal{L}_q^{r+0} = \mathcal{B}_{q+1}^{r+0} (\mathcal{B}_{q+1}^{r+0})^T + (\mathcal{B}_{q}^{r})^T \mathcal{B}_{q}^{r}.
    \end{equation}
    The volume of any $1$-simplex $\sigma_1 = [v_i,v_j]$ is $\text{Vol}(\sigma_1)$ is actually the distance between $v_i$ and $v_j$ denoted $d_{ij}$. Then the $0$-persistent $0$-combinatorial Laplacian matrix based on filtration $r$ can be expressed explicitly as
    \begin{equation}\label{equ:contact 1}
        (\mathcal{L}_0^{r+0})_{ij}=
        \begin{cases}
            -\displaystyle{\sum_{j}}(\mathcal{L}_0^{r+0})_{ij}, & \mbox{if $i=j$} \\
            -1,                    & \mbox{if $i\neq j$ and $d_{ij}-2r<0$} \\
            0,                     & \mbox{otherwise.}
        \end{cases}
    \end{equation}
    Correspondingly, we can denote the $0$-persistent $1$-combinatorial Laplacian matrix based on filtration $r$ by $\mathcal{L}_1^{r+0}$, and the $0$-persistent $2$-combinatorial Laplacian matrix based on filtration $r$ by $\mathcal{L}_2^{r+0}$.

    Alternatively, variants of persistent $0$-combinatorial Laplacian matrices can be defined by adding the Euclidean distance information. The distance-weight persistent $0$-combinatorial Laplacian matrix  based on filtration $r$ can be expressed explicitly as
    \begin{equation}\label{equ:dis}
        (\hat{\mathcal{L}}_0^{r+0})_{ij}=
        \begin{cases}
            -\displaystyle{\sum_{j}}(\hat{\mathcal{L}}_0^{r+0})_{ij}, & \mbox{if $i=j$} \\
            -d_{ij},                              & \mbox{if $i\neq j$ and $d_{ij}-2r<0$} \\
            0,                                    & \mbox{otherwise.}
        \end{cases}
    \end{equation}
    Moreover, the inverse-distance-weight persistent $0$-combinatorial Laplacian matrix   based on filtration $r$ can also be implemented:
    \begin{equation}\label{equ:inv dis}
        (\check{\mathcal{L}}_0^{r+0})_{ij}=
        \begin{cases}
            -\displaystyle{\sum_{j}}(\check{\mathcal{L}}_0^{r+0})_{ij}, & \mbox{if $i=j$} \\
            -\displaystyle{\frac{1}{d_{ij}}},                    & \mbox{if $i\neq j$ and $d_{ij}-2r<0$} \\
            0,                     & \mbox{otherwise.}
        \end{cases}
    \end{equation}

  The spectra of the aforementioned $0$-persistent $0$-combinatorial Laplacian matrices based on filtration are given by
    \[
    \begin{split}
        \text{Spectra}(\mathcal{L}_0^{r+0}) &= \{ (\lambda_1)_0^{r+0}, (\lambda_2)_0^{r+0}, \cdots, (\lambda_N)_0^{r+0} \}, \\
        \text{Spectra}(\hat{\mathcal{L}}_0^{r+0}) &= \{ (\hat{\lambda}_1)_0^{r+0}, (\hat{\lambda}_2)_0^{r+0}, \cdots, (\hat{\lambda}_N)_0^{r+0} \}, \\
        \text{Spectra}(\check{\mathcal{L}}_0^{r+0}) &= \{ (\check{\lambda}_1)_0^{r+0}, (\check{\lambda}_2)_0^{r+0}, \cdots, (\check{\lambda}_N)_0^{r+0} \},
    \end{split}
    \]
    where $N$ is the dimension of persistent Laplacian matrices, $(\hat{\lambda}_j)_0^{r+0}$ and
        $(\check{\lambda}_j)_0^{r+0}$ are the $j$-th eigenvalues of $\hat{\mathcal{L}}_0^{r+0}$ and  $\check{\mathcal{L}}_0^{r+0}$, respectively.
     We denote $\hat{\beta}_q^{r+0}$ and $\check{\beta}_q^{r+0}$ the $q$th Betti for $\hat{\mathcal{L}}_q^{r+0}$ and $\check{\mathcal{L}}_q^{r+0}$, respectively.

    The smallest non-zero eigenvalue  of $\mathcal{L}_0^{r+0}$, denoted $(\tilde{\lambda}_2)_0^{r+0}$, is particularly useful in many applications. Similarly, the smallest non-zero eigenvalues  of $\hat{\mathcal{L}}_0^{r+0}$ and $\check{\mathcal{L}}_0^{r+0}$ are denoted as         $(\tilde{\hat{\lambda}}_2)_0^{r+0}$ and
    $(\tilde{\check{\lambda}}_2)_0^{r+0}$, respectively.

    Finally, it is mentioned that using the present procedure, more general weights, such as the radial basis function of the Euclidean distance, can be employed to construct weighted boundary operators and associated persistent combinatorial Laplacian matrices.

    \subsubsection{Multiscale spectral analysis}

 In the past few years, we have developed a multiscale spectral graph method such as generalized GNM and generalized ANM \cite{opron2015communication,xia2015multiscale}, to create a family of spectral graphs with different characteristic length scales for a given dataset. Similarly, in our persistent spectral theory, we can construct a family of spectral graphs induced by a filtration parameter. Moreover, we can sum over all the multiscale spectral graphs as an accumulated spectral graph. Specifically, a family of $\mathcal{L}_0^{r+0}$ matrices, as well as the accumulated combinatorial Laplacian matrices, can be generated via the filtration. By analyzing the persistent spectra of these matrices, the topological invariants and geometric shapes can be revealed from the given input point-cloud data.
	 
\begin{figure}[H]
        \centering
        \captionsetup{margin=0.9cm}
        \includegraphics[scale=0.15]{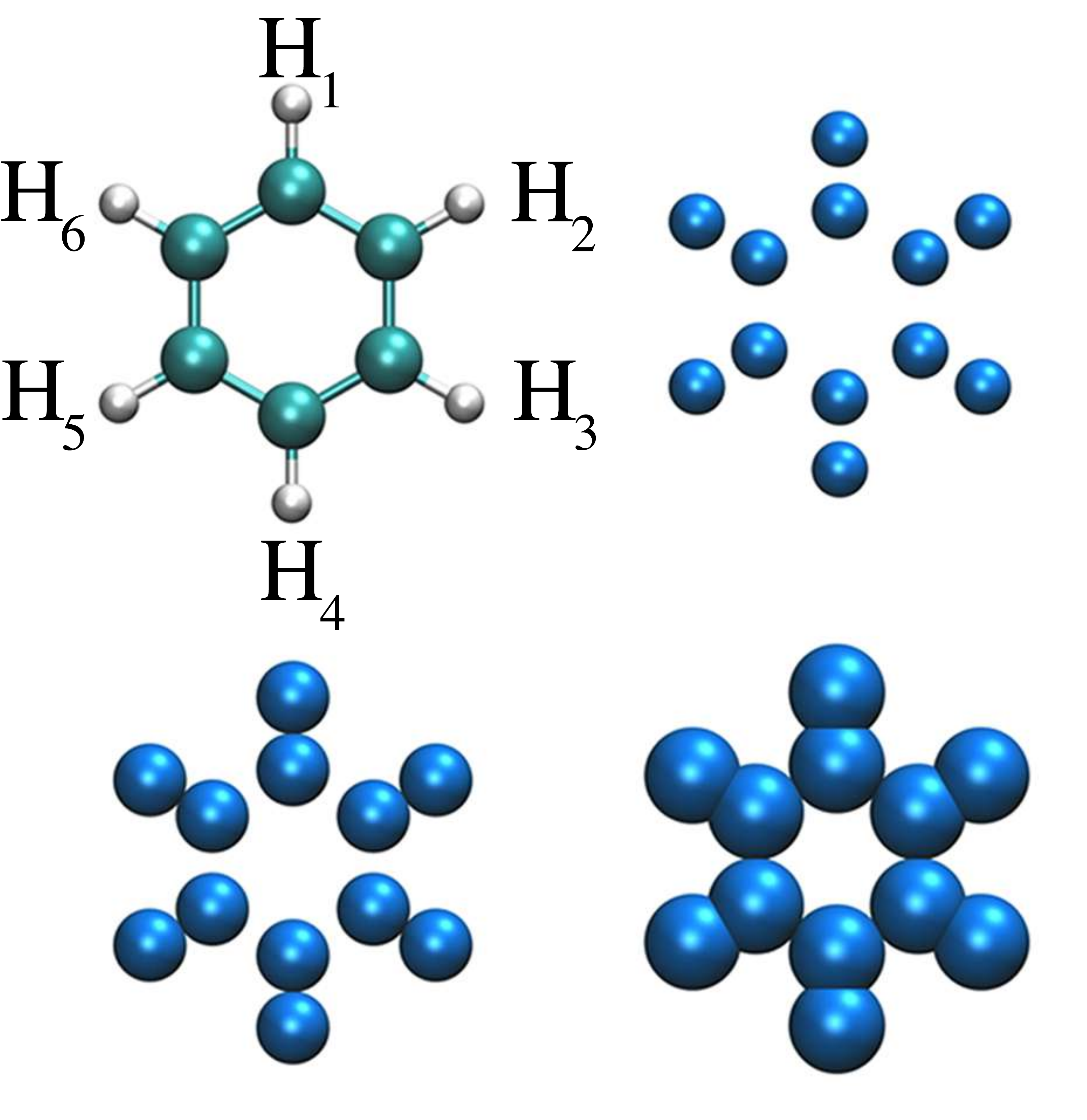}
        \caption{Benzene molecule and its topological changes during the filtration process. }
        \label{fig:Benzene Combine}
    \end{figure}

\begin{figure}[H]
        \centering
        \captionsetup{margin=0.9cm}
        \subfigure[]{
            \label{fig:subfig:Benzene Sec}
            \includegraphics[scale = 0.40]{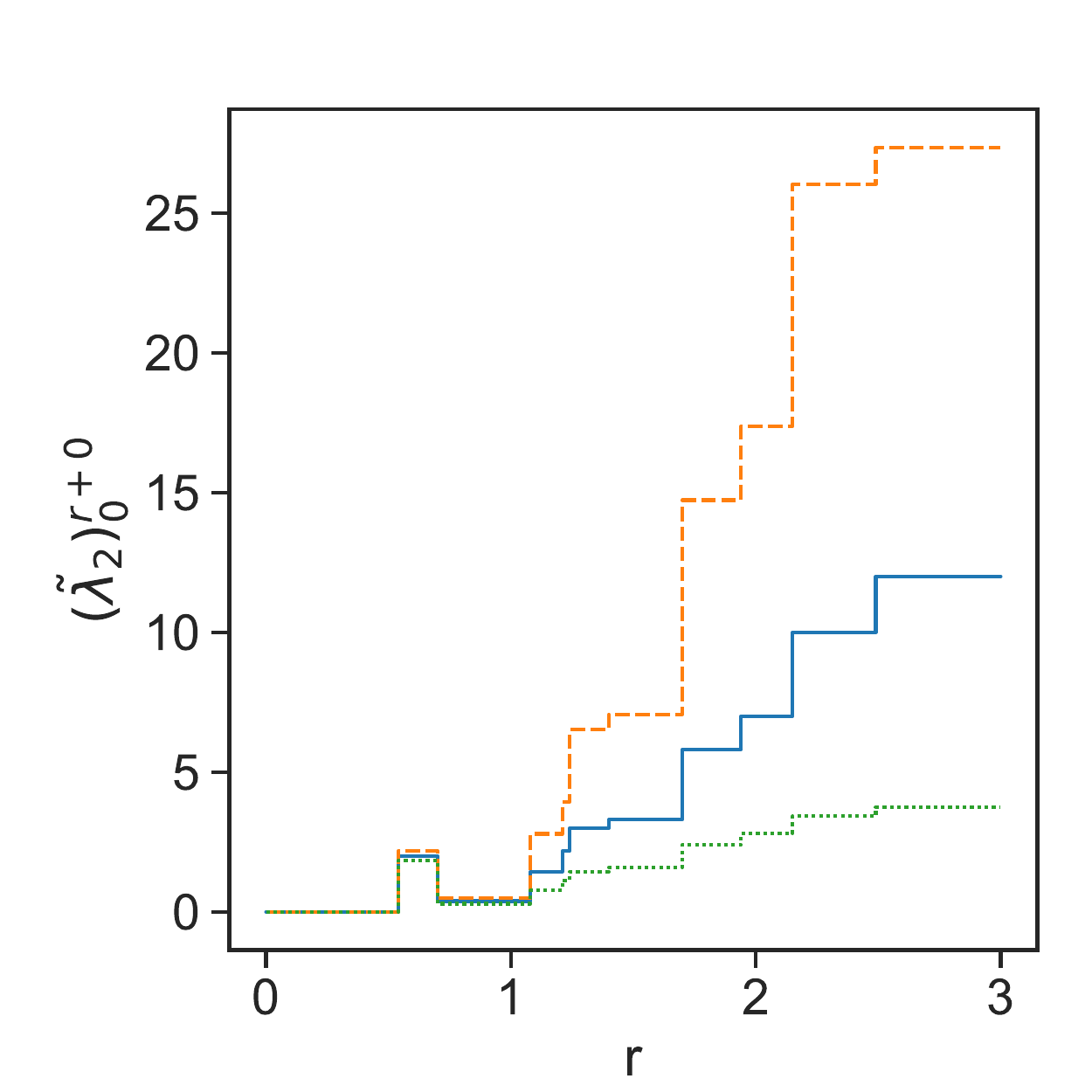}}
            \hspace{0.001\linewidth}
        \subfigure[]{
            \label{fig:subfig:Benzene Zero}
            \includegraphics[scale = 0.40]{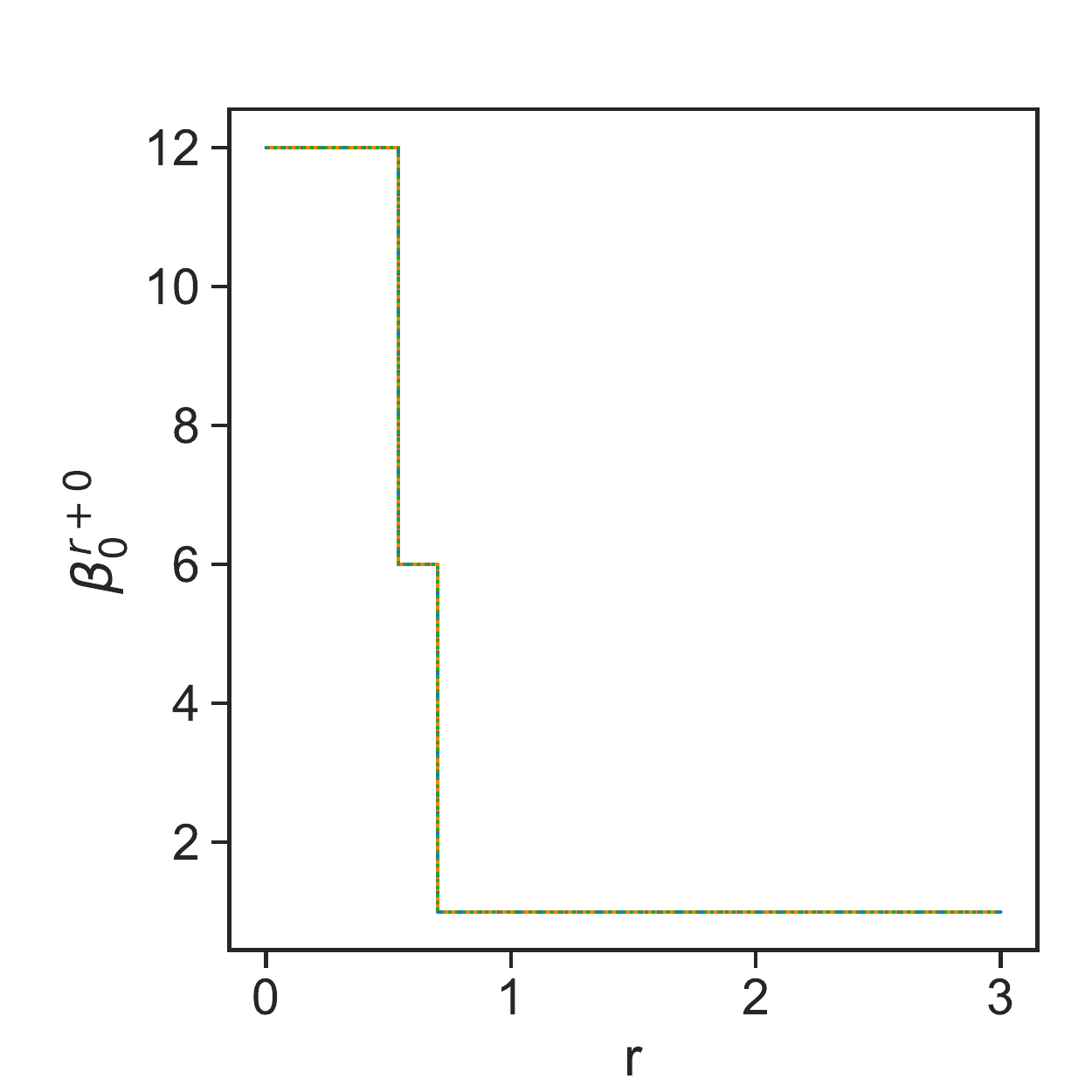}}
            \hspace{0.001\linewidth}
        \subfigure[]{
            \label{fig:subfig:Benzene betti1}
            \includegraphics[scale = 0.40]{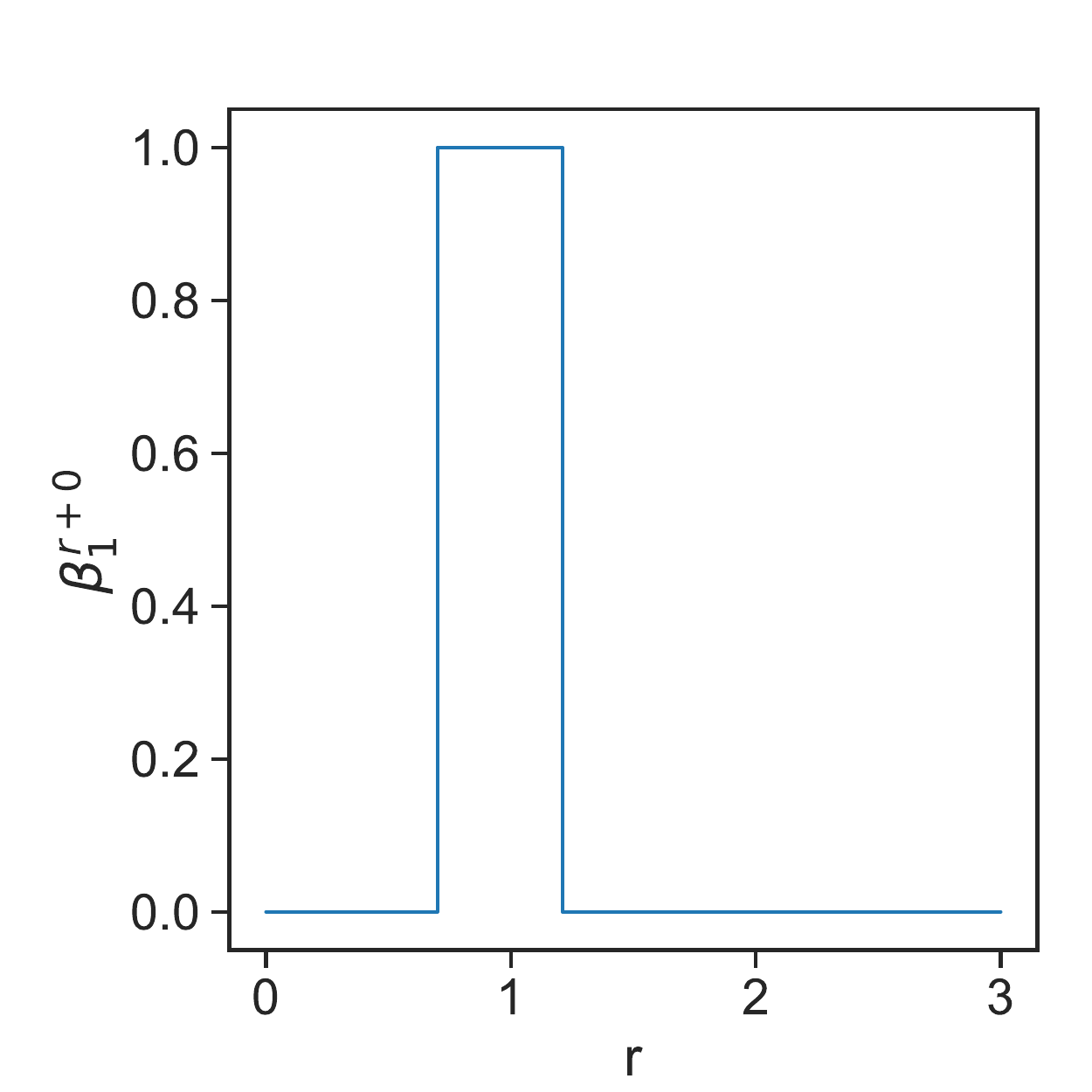}}
        \caption{Persistent spectral analysis of the benzene molecule induced by filtration parameter $r$. Blue line, orange line, and green line represent $\mathcal{L}_0^{r+0}$, $\hat{\mathcal{L}}_0^{r+0}$, and $\check{\mathcal{L}}_0^{r+0}$ respectively.
				(a) Plot of the smallest non-zero eigenvalues with radius filtration under $\mathcal{L}_0^{r+0}$ (blue line), $\hat{\mathcal{L}}_0^{r+0}$ (red line), and $\check{\mathcal{L}}_0^{r+0}$ (green line). Total $10$ jumps observed in this plot which represent $10$ possible distances between atoms.
				(b) Plot of  the number of zero eigenvalues ($\beta_0^{r+0}$) with radius filtration under $\mathcal{L}_0^{r+0}$, $\hat{\mathcal{L}}_0^{r+0}$, and $\check{\mathcal{L}}_0^{r+0}$ (three spectra are superimposed). When $r = \SI{0.00}{\angstrom}$, $12$ atoms are disconnected with each other. After $r = \SI{0.54}{\angstrom}$, H atoms and their adjacent C atoms are connected with one another resulting in $\beta_0^{r+0} = 6$. With $r$ keeps growing, all of the atoms are connected with one another and then $\beta_0^{r+0} = 1$.
				(c) Plot of the number of zero eigenvalues ($\beta_1^{r+0}$) with radius filtration under $\mathcal{L}_1^{r+0}$. When $r = \SI{0.70}{\angstrom}$, a $1$-cycle created since all of the C atoms are connected and form a hexagon, resulting in $\beta_1^{r+0} = 1$. After the radius reached $\SI{1.21}{\angstrom}$, the hexagon disappears and $\beta_1^{r+0} = 0$.}
        \label{fig:Benzene}
    \end{figure}

    The spectra of $\mathcal{L}_0^{r+0}$, $\hat{\mathcal{L}}_0^{r+0}$, and $\check{\mathcal{L}}_0^{r+0}$ mentioned above carry similar information on how the topological structures of a graph are changed during the filtration. Benzene molecule (C$_6$H$_6$), a typical aromatic hydrocarbon which is composed of six carbon atoms bonded in a planar regular hexagon ring with one hydrogen joined with each carbon atom. It provides a good example to demonstrate the proposed PST.  		\autoref{fig:Benzene Combine} illustrates the filtration of  the benzene molecule. Here, we label 6 hydrogen atoms by H$_1$, H$_2$, H$_3$, H$_4$, H$_5$, and H$_6$, and the carbon adjacent to the labeled hydrogen atoms are labeled by C$_1$, C$_2$, C$_3$, C$_4$, C$_5$, and C$_6$, respectively. Figure \autoref{fig:subfig:Benzene Zero} depicts that when the radius of the solid sphere reaches $\SI{0.54}{\angstrom}$, each carbon atom in the benzene ring is overlapped with its joined hydrogen atom, resulting in the reduction of $\beta_0^{r+0}$ to $6$. Moreover, once the radius of solid spheres is larger than $\SI{0.70}{\angstrom}$, all the atoms in the benzene molecule will connect and constitute a single component which gives rise  $\beta_0^{r+0}=1$. Furthermore, we can deduce that the C-C bond length of the benzene ring is about $\SI{1.40}{\angstrom}$, and the C-H bond length is around $\SI{1.08}{\angstrom}$, which are the real bond lengths in benzene molecule. Figure \autoref{fig:subfig:Benzene betti1} shows that a $1$-dimensional hole ($1$-cycle) is born when the filtration parameter $r$ increase to $\SI{0.70}{\angstrom}$ and dead when $r = \SI{1.21}{\angstrom}$. In Figures \autoref{fig:subfig:Benzene Zero} and   \autoref{fig:subfig:Benzene betti1}, it can be seen that variants of $0$-persistent $0$-combinatorial Laplacian and $1$ -combinatorial Laplacian matrices based on filtration give us the identical $\beta_0^{r+0}$ and $\beta_1^{r+0}$ information respectively.

     The C-C bond length of benzene is $\SI{1.39}{\angstrom}$, and the C-H bond length is $\SI{1.09}{\angstrom}$.    Due to the perfect hexagon structure of the benzene ring, we can calculate all of the distances between atoms.      The shortest and longest distances between carbons and the hydrogen atoms are $\SI{1.09}{\angstrom}$ and $\SI{3.87}{\angstrom}$. In Figure \autoref{fig:subfig:Benzene Sec}, a total of $10$ changes of $(\tilde{\lambda}_2)_0^{r+0}$ values is observed at various radii.   \autoref{table:Benzene distance} lists all the distances between atoms and the values of radii when the changes of $(\tilde{\lambda}_2)_0^{r+0}$ occur. It can be seen that the distance between atoms approximately equals twice of the radius value when a jump of $(\tilde{\lambda}_2)_0^{r+0}$ occurs. Therefore, we can detect all the possible distances between atoms with the nonzero spectral information. Moreover, in Figure \autoref{fig:subfig:Benzene Zero}, the values of the smallest nonzero eigenvalues of $\mathcal{L}_0^{r+0}$, $\hat{\mathcal{L}}_0^{r+0}$, and $\check{\mathcal{L}}_0^{r+0}$ change concurrently.

    \begin{table}[H]
        \centering
        \setlength\tabcolsep{3pt}
        \captionsetup{margin=0.9cm}
        \caption{Distances between atoms in the benzene molecule and the radii when the changes of $(\tilde{\lambda}_2)_0^{r+0}$ occur (Values increase from left to right).  
				}
        \begin{tabular}{ccccccccccc}
        \hline
         Type                           & C$_1$-H$_1$     & C$_1$-C$_2$     & C$_2$-H$_1$     & C$_1$-C$_3$     & H$_1$-H$_2$     & C$_1$-C$_4$     & C$_3$-H$_1$     & C$_4$-H$_1$     & H$_1$-H$_3$     & H$_1$-H$_4$       \\ \hline
         Distance ($\si{\angstrom}$)     & $1.09$  & $1.39$  & $2.15$  & $2.41$  & $2.48$  & $2.78$  & $3.39$  & $3.87$  & $4.30$  & $4.96$    \\
         r ($\si{\angstrom}$)       & $0.54$  & $0.70$  & $1.08$  & $1.21$  & $1.24$  & $1.40$  & $1.70$  & $1.94$  & $2.15$  & $2.49$    \\ \hline
        \end{tabular}
        \label{table:Benzene distance}
    \end{table}

\section{Applications}\label{sec: Application}

In this section, we apply the proposed persistent spectral theory to the study of two important systems, fullerenes and proteins.  All three different types of persistent combinatorial Laplacian matrices are employed in our investigation. The resulting persistent spectra contain not only the full set of topological persistence from the harmonic spectra, which is identical to that from a persistent homology analysis, but also non-harmonic eigenvalues and eigenvectors. Since the power of topological persistence has been fully explored and exploited in the past decade \cite{de2007coverage, YaoY:2009,bubenik2014categorification,dey2014computing},  to   demonstrate the additional utility of our persistent spectral analysis, we  mainly emphasize  the non-harmonic spectra in  the present application. In practical applications, such as drug design, the use of the full spectra from the  persistent spectral theory would certainly further enhance our models.     

\subsection{Fullerene analysis and prediction}
    In 1985 Kroto et all discovered the first structure of C$_{60}$  \cite{Kroto1985}, which was confirmed by Kratschmer et al in 1990  \cite{Kratschmer1990}. Since then, the quantitative analysis of fullerene molecules has become an interesting research topic. The understanding of the fullerene structure-function relationship is important for nanoscience and nanotechnology. Fullerene molecules are only made of carbon atoms that have various topological shapes, such as the hollow spheres, ellipsoids, tubes, or rings. Due to the monotony of the atom type and the variety of geometric shapes, the minor heterogeneity of fullerene structures can be ignored. The fullerene system offers a moderately  large dataset with relatively simple structures. Therefore, it is suitable for validating new computational methods because every single change in the spectra is interpretable. The proposed persistent spectral theory, i.e., persistent spectral analysis, is applied to characterize fullerene structures and predict their stability. 
        
    All the structural data can be downloaded from \href{http://www.ccl.net/cca/data/fullerenes/index.shtml}{CCL.NET Webpage}. This dataset gives the coordinates of fullerene carbon atoms. In this section, we will analyze fullerene structures and predict the heat of formation energy. 

 \subsubsection{Fullerene structure analysis}
         The smallest member of the fullerene family is C$_{20}$ molecule with a dodecahedral cage structure. Note that $12$ pentagons are required to form a closed fullerene structure. Following the Euler's formula, the number of vertices, edges, and faces on a polygon have the relationship $V-E+F=2$. Therefore, the $20$ carbon atoms in the dodecahedral cage form $30$ bonds with the same bond length. The C$_{20}$ is the only fullerene smaller than C$_{60}$ that has the molecular symmetry of the full icosahedral point group $I_h$. C$_{60}$ is a molecule that consists of $60$ carbon atoms arranged as $12$ pentagon rings and $20$ hexagon rings. Unlike C$_{20}$, C$_{60}$ has two types of bonds:  $6:6$ bonds and $6:5$ bonds. The $6:6$ bonds are shorter than $6:5$ bonds, which can also be considered as  ``double bond" \cite{Yadav2008}.  C$_{60}$ is the most well-know fullerene with geometric symmetry $I_h$. Since C$_{20}$ and C$_{60}$ are highly symmetrical,  they are ideal systems for illustrating the persistent spectral analysis.

        \begin{figure}[H]
            \centering
            \captionsetup{margin=0.9cm}
            \includegraphics[scale=0.22]{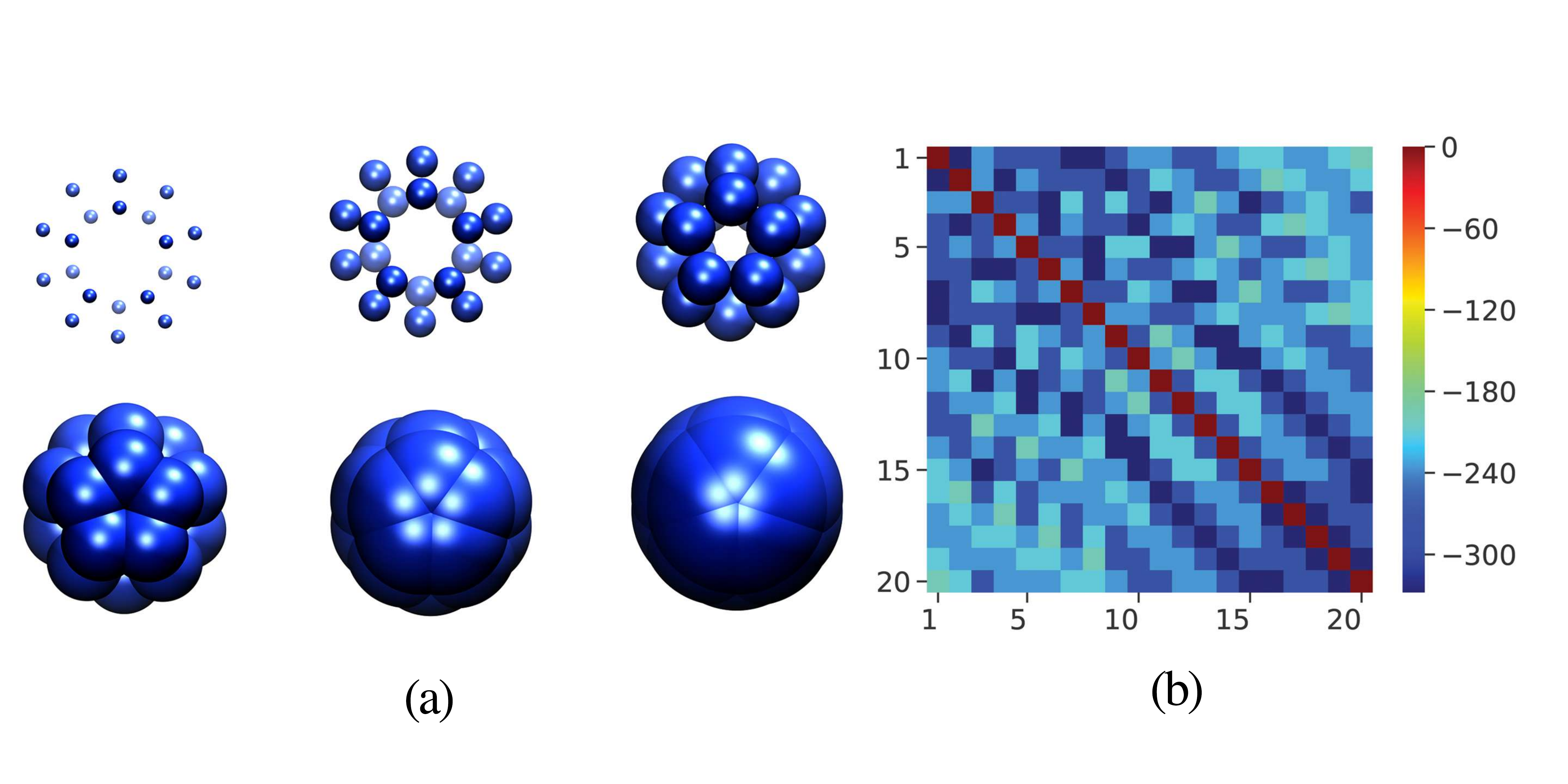}
            \caption{(a) Illustration of filtration built on fullerene C$_{20}$. Each carbon atom of C$_{20}$ is plotted by its given coordinates, which are associated with an ever-increasing radius $r$. The solid balls centered at given coordinates keep growing along with the radius filtration parameter. (b) The accumulated $\mathcal{L}_0^{r+0}$ matrix for C$_{20}$. For clarity, the diagonal terms are set to 0.}
            \label{fig:C20}
        \end{figure}

\begin{figure}[H]
            \centering
            \captionsetup{margin=0.9cm}
            \includegraphics[scale=0.44]{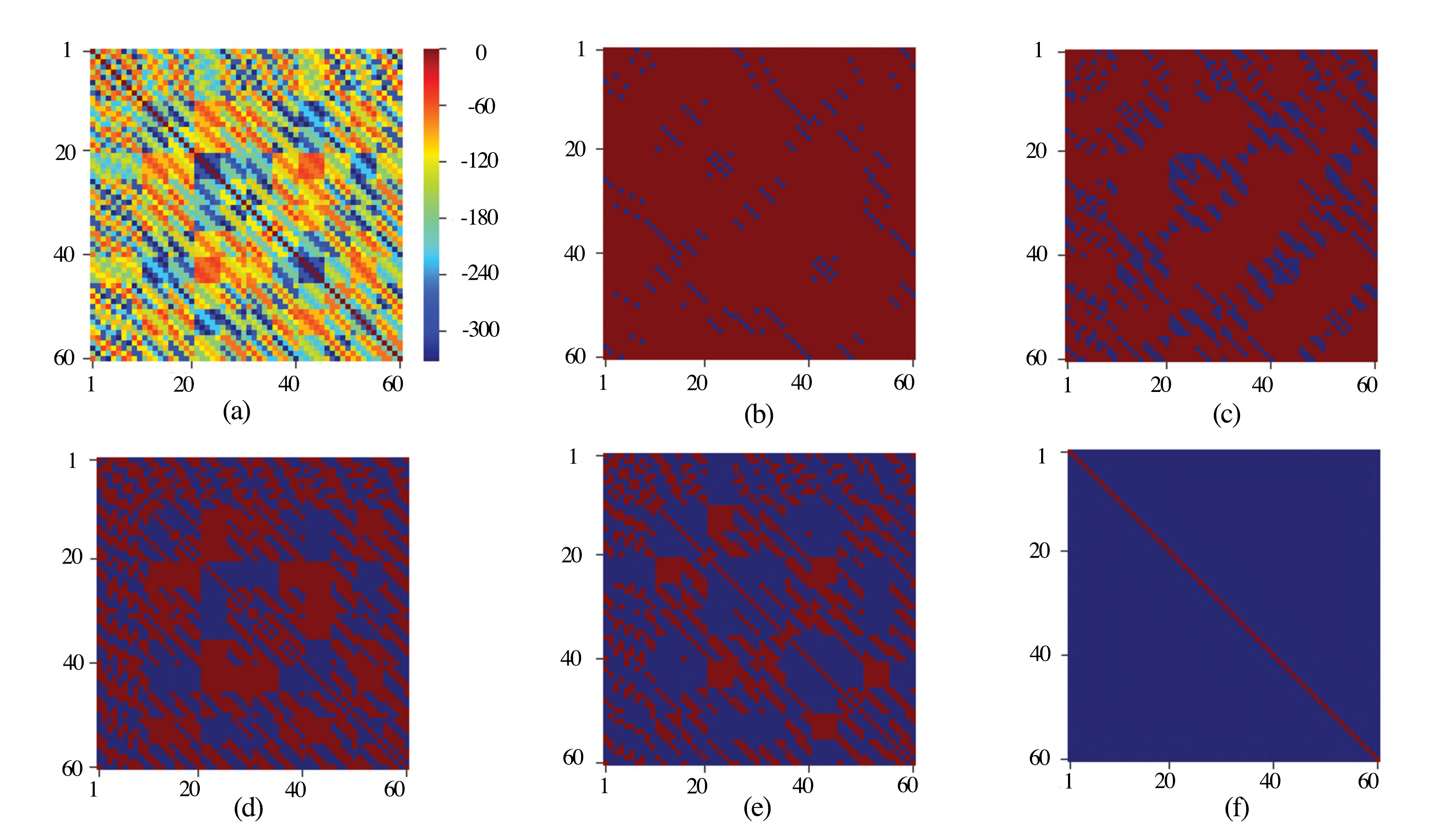}
            \caption{Illustration of persistent multiscale analysis of C$_{60}$ in terms of 0-combinatorial Laplacian matrices (b)-(f) and their accumulated matrix (a) induced by filtration.  
						As the value of filtration parameter $r$ increases, high-dimensional simplicial complex forms and grows accordingly.  (b), (c), (d), (e), and (d) demonstrate  the  0-combinatorial Laplacian matrices (i.e., the connectivity among C$_{60}$ atoms) at filtration $r =  \SI{1.0}{\angstrom},  \SI{1.5}{\angstrom},  \SI{2.5}{\angstrom}$, $ \SI{3.0}{\angstrom}$, and $\SI{3.6}{\angstrom}$, respectively. The blue cell located at the $i$th row and $j$th column represents the balls centered at atom $i$ and atom $j$ connected with each other. For clarity, the diagonal terms are set to 0 in all plots.}
            \label{fig:C60 LapFil}
        \end{figure}

      \autoref{fig:C20} (a) illustrates the radius filtration process built on C$_{20}$. As the radius increases, the solid balls corresponding to carbon atoms grow, and a sequence of $\mathcal{L}_0^{r+0}$ matrices can be defined through the overlap relations among the set of balls. At the initial state ($r = \SI{0.00}{\angstrom}$), all of the atoms are isolated from one another. Therefore,  $\mathcal{L}_0^{r+0}$ is a zero matrix with dimension $20\times 20$. Since the C$_{20}$ molecule has the same bond length which can be denoted as $l(\mbox{C}_{20})$, once the radius of solid balls is greater than $l(\mbox{C}_{20})$, all of the balls are overlapped, which makes the system a singly connected component.
        \autoref{fig:C20} (b) depicts the accumulated $\mathcal{L}_0^{r+0}$ for C$_{20}$. For C$_{60}$, the accumulated $\mathcal{L}_0^{r+0}$ is described in \autoref{fig:C60 LapFil} (a). \autoref{fig:C60 LapFil} (b)-(f) are the plots of $\mathcal{L}_0^{r+0}$ under different filtration $r$ values. The blue cell located at the $i$th row and $j$th column means the balls centered at atom $i$ and atom $j$ connected with each other, i.e., a $1$-simplex formed with its vertex to be $i$ and $j$. When the radius filtration increases, more and more bluer cells are created. In \autoref{fig:C60 LapFil} (f), the color of cells, except the cells located in the diagonal, turns to blue, which means all of the carbon atoms are connected with one another at $r = \SI{3.6}{\angstrom}$. For clarity, we set the diagonal terms to 0.
				
\begin{figure}[H]
            \centering
            \captionsetup{margin=0.9cm}
            \subfigure[$\beta_0^{r+0}$]{
                \label{fig:subfig:Fullerene zero}
                \includegraphics[scale=0.45]{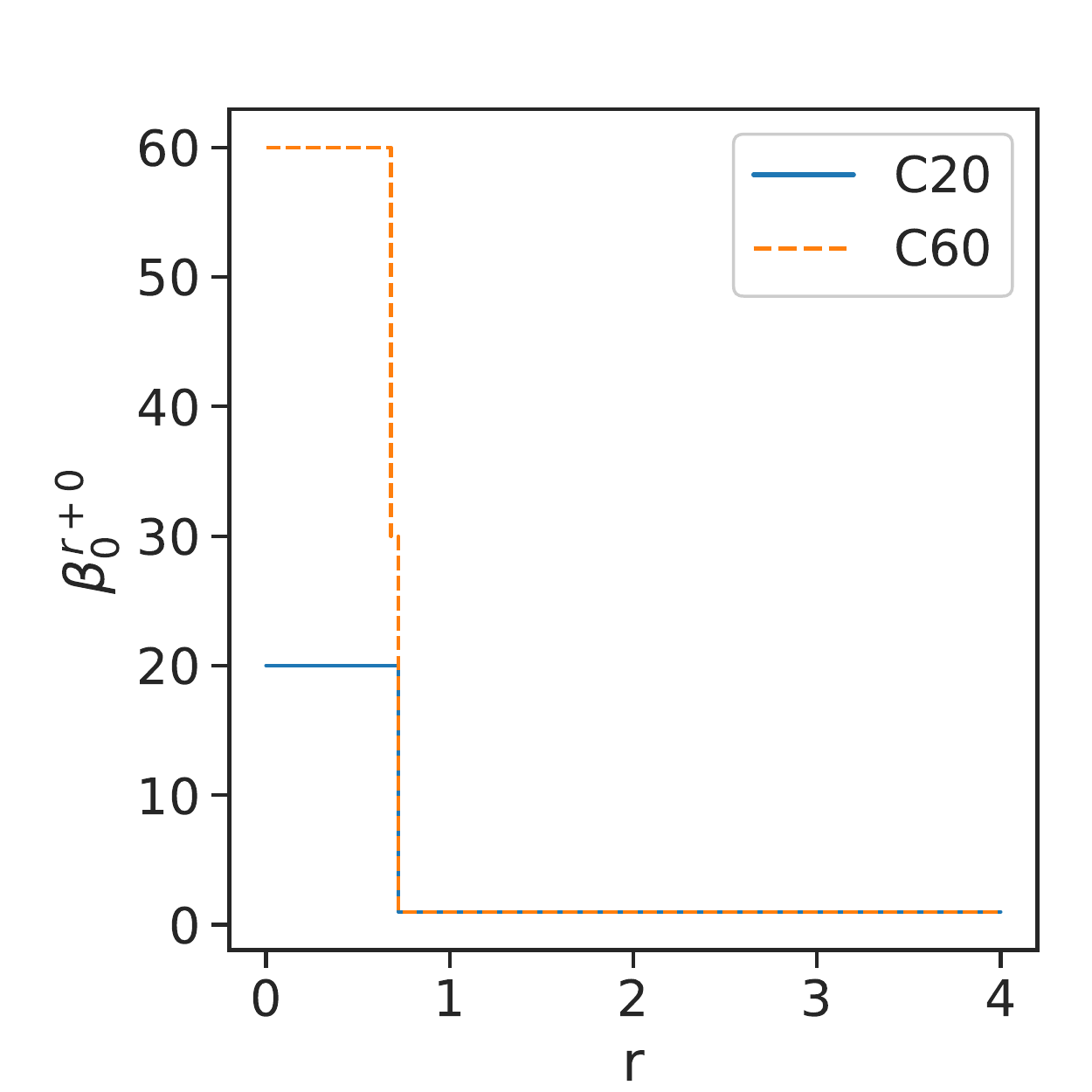}}
                \hspace{0.001in}
            \subfigure[$\beta_1^{r+0}$]{
                \label{fig:subfig:Fullerene betti1}
                \includegraphics[scale = 0.45]{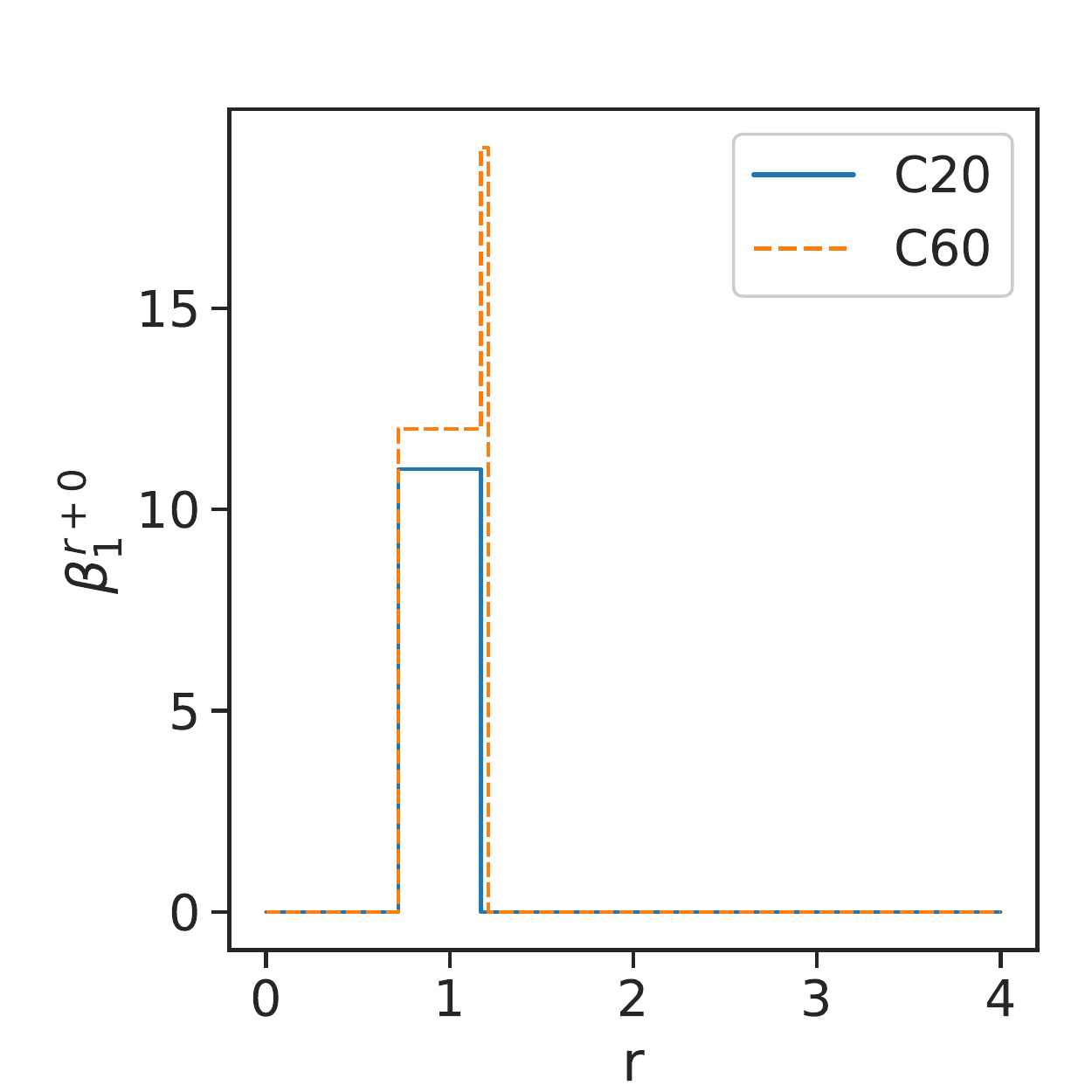}}
                \vfill
            \subfigure[$\beta_2^{r+0}$]{
                \label{fig:subfig:Fullerene betti2}
                \includegraphics[scale=0.45]{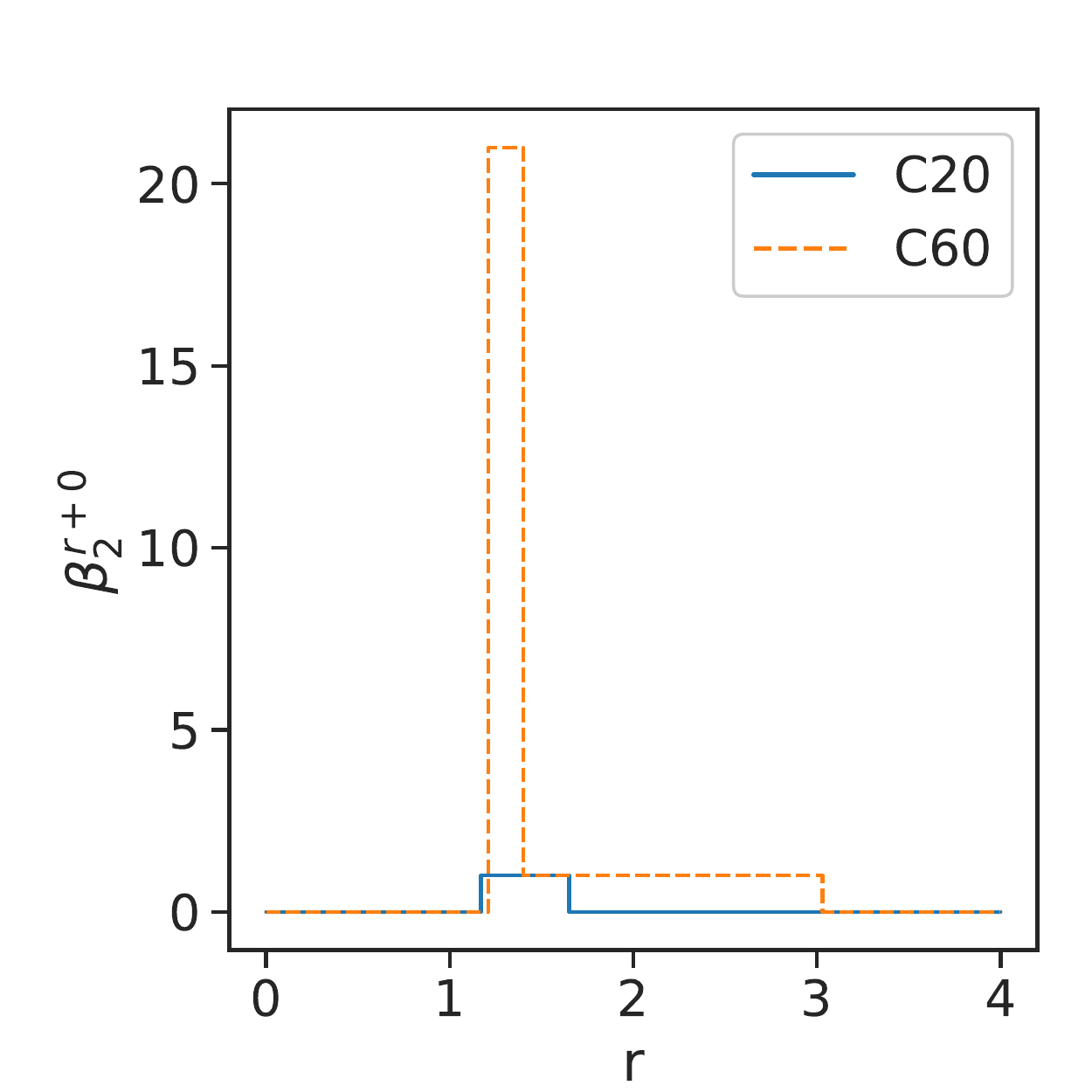}}
                \hspace{0.001in}
            \subfigure[$(\tilde{\lambda}_2)_0^{r+0}$]{
                \label{fig:subfig:Fullerene sec}
                \includegraphics[scale=0.45]{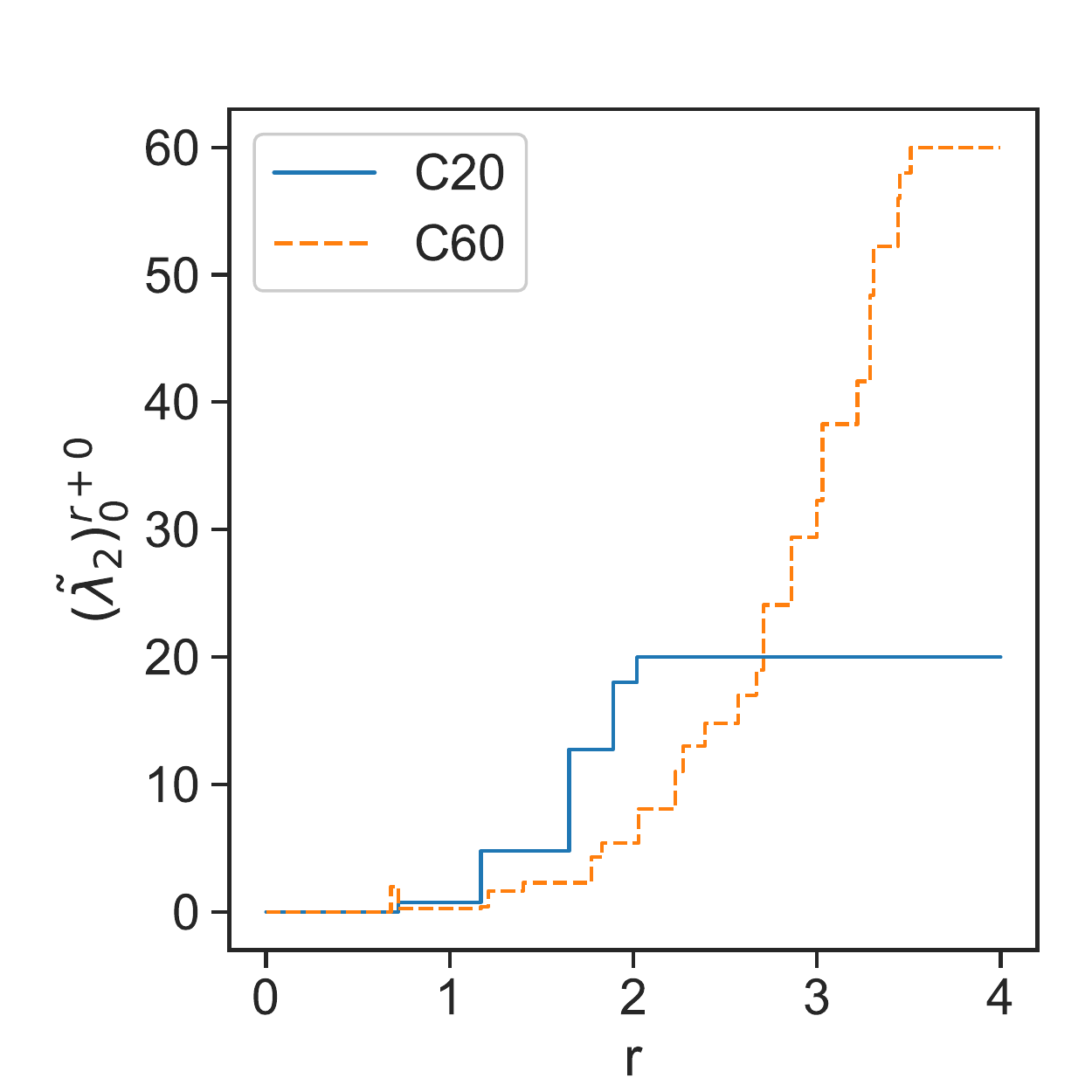}}
                \hspace{0.001in}
            \caption{ Illustration of persistent spectral analysis of C$_{20}$ and C$_{60}$ using the spectra of $\mathcal{L}_q^{r+0}$ ($q=1,2$ and 3).    
						(a) The number of zero eigenvalues of $\mathcal{L}_0^{r+0}$, i.e., $\beta_0^{r+0}$, under radius filtration. 
						(b) The number of zero eigenvalues of $\mathcal{L}_1^{r+0}$, i.e., $\beta_1^{r+0}$ under radius filtration. 
						(c) The number of zero eigenvalues of $\mathcal{L}_2^{r+0}$, i.e., $\beta_2^{r+0}$ under radius filtration. 
						(d) The smallest non-zero eigenvalue $(\tilde{\lambda}_2)_0^{r+0}$ under  radius filtration. The radius grid spacing is $\SI{0.01}{\angstrom}$.}
            \label{fig:C20C60}
        \end{figure}

        In \autoref{fig:C20C60}, the blue solid line represents C$_{20}$ properties and the dash orange line represents C$_{60}$ properties. For Figure \autoref{fig:subfig:Fullerene zero}, the blue line drops at $r=\SI{0.72}{\angstrom}$, which means the bond length of C$_{20}$ is around $\SI{1.44}{\angstrom}$. The orange line drops at $r=\SI{0.68}{\angstrom}$ and $\SI{0.72}{\angstrom}$, which means the ``double bond" length of C$_{60}$ is around $\SI{1.36}{\angstrom}$ and the $6:5$ bond length is around $\SI{1.44}{\angstrom}$. Moreover, the total number of ``double bond" is $30$, yielding $\beta_0^{r+0}=30$ when the radius of solid balls is over $\SI{0.68}{\angstrom}$. In conclusion, one can deduce the number of different types of bonds as well as the bond length information from the number of zero eigenvalues (i.e., $\beta_0^{r+0}$) under the radius filtration. Furthermore, the geometric information can also be derived from the plot of $(\tilde{\lambda}_2)_0^{r+0}$. Each jump in Figure \autoref{fig:subfig:Fullerene sec} at a specific radius represents the change of geometric and topological structure. The smallest non-zero eigenvalue $(\tilde{\lambda}_2)_0^{r+0}$ of $\mathcal{L}_0^{r+0}$ matrices for C$_{20}$ changes $5$ times in Figure \autoref{fig:subfig:Fullerene sec}, which means C$_{20}$ has $5$ different distances between carbon atoms. Furthermore, by Remark \autoref{rem1}, as $(\tilde{\lambda}_2)_0^{r+0}$ of C$_{20}$ keeps increasing,   the smallest vertex connectivity of the connected subgraph  continues growing and the topological structure  becomes steady. As can be seen in the right-corner chart of \autoref{fig:C20}, the carbon atoms will finally grow to a solid object with a steady topological structure. 
				
Figure \autoref{fig:subfig:Fullerene betti1} depicts the changes of Betti $1$ value $\beta_1^{r+0}$ (i.e., the number of zero eigenvalues for $\mathcal{L}_1^{r+0}$) under the filtration $r$. Since C$_{20}$ has $12$ pentagonal rings, $\beta_1^{r+0}$ jumps to $11$ when radius $r$ equals to the half of the bond length of $l(\text{C}_{20})$. These eleven  $1$-cycles disappear at $r = \SI{1.17}{\angstrom}$. There are $12$ pentagons and $20$ hexagons in C$_{60}$, which results in $\beta_1^{r+0} = 12$ at $r = \SI{0.72}{\angstrom}$, $\beta_1^{r+0} = 31$ at $r = \SI{1.17}{\angstrom}$. All of the pentagons and hexagons  disappear at $r = \SI{1.22}{\angstrom}$.

As the filtration process, even more structure information can be derived from the number of zero eigenvalues of $\mathcal{L}_2^{r+0}$ (i.e., $\beta_2^{r+0}$)  in Figure \autoref{fig:subfig:Fullerene betti2}. For C$_{20}$, $\beta_2^{r+0} = 1$ when $r = \SI{1.17}{\angstrom}$, which corresponds to the void structure  in the center of the dodecahedral cage. The void disappears at $r = \SI{1.65}{\angstrom}$ since a solid structure is generated at this point. For fullerene C$_{60}$, $20$ hexagonal cavities and a center void exist from $\SI{1.12}{\angstrom}$ to $\SI{1.40}{\angstrom}$ yielding $\beta_2^{r+0} = 21$. As the filtration goes, hexagonal cavities disappear which results $\beta_2^{r+0}$ decrease to $1$. The central void  keeps alive until  a solid block is formed at $r = \SI{3.03}{\angstrom}$. In a nutshell, we can deduce the number of different types of bonds, the bond length, and the topological invariants from the present persistent spectral analysis.

 \subsubsection{Fullerene stability  prediction}

\begin{figure}[H]
            \centering
            \captionsetup{margin=.9cm}
								 \includegraphics[scale=0.2]{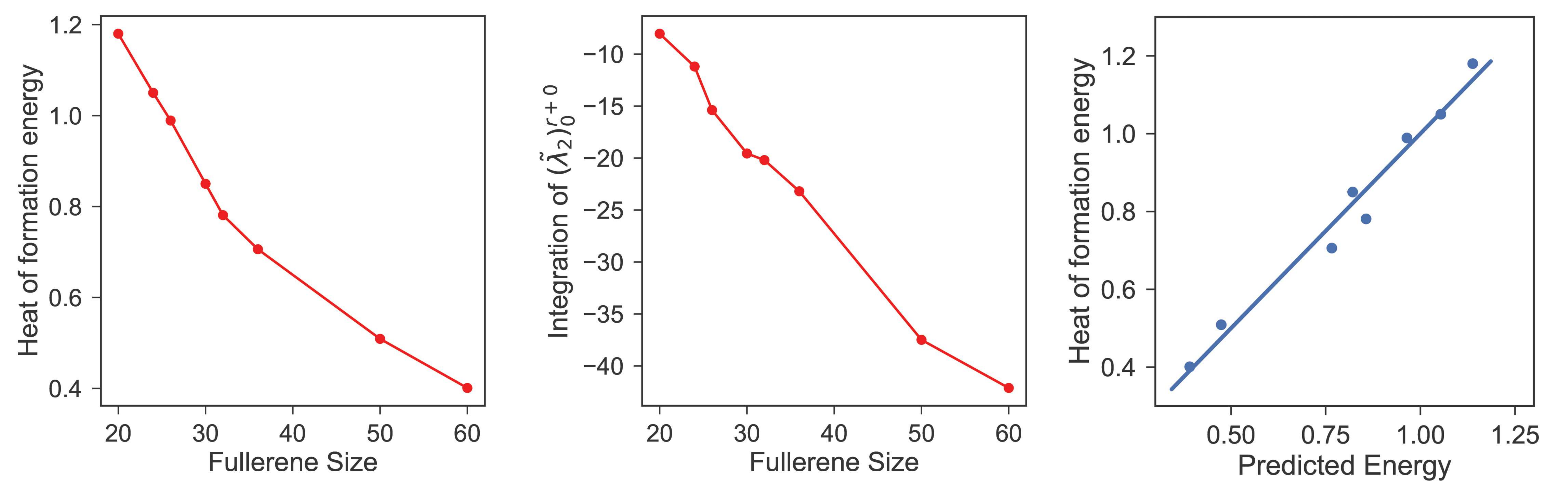}
            \caption{Persistent spectral analysis and prediction of fullerene heat formation energies.
						Left chart: the heat of formation energies of fullerenes obtained from quantum  calculations \cite{Zhang1992a}.
						Middle chart: PST model using the area under the plot of $(\tilde{\lambda}_2)_0^{r+0}$.
            Right chart: Correlation between the quantum calculation and the PST prediction. The highest correlation coefficient form the least-squares fitting is $0.986$ with the type index  of $\alpha=\mbox{Max}$.
						}
						\label{fig:Heat comparison}
        \end{figure}

        Having shown that the detailed fullerene structural information can be extracted into the spectra of $\mathcal{L}_q^{r+0}$, we further illustrate that fullerene functions can be predicted from their structures by using our persistent spectral theory in this section. Similar structure-function analysis  has been carried out by using other methods \cite{opron2015communication,xia2015persistent,Science2008}. For small fullerene molecule series C$_{20}$ to C$_{60}$, with the increase in the number of atoms, the ground-state heat of formation energies decrease \cite{Zhang1992,Zhang1992a}. The left chart in \autoref{fig:Heat comparison} describes this phenomenon. Similar patterns can also be found in the total energy (STO-$3$G/SCF at MM$3$) per atom and the average binding energy of C$_{2n}$. To analyze these patterns, many theories have been proposed. Isolated pentagon rule assumes that the most stable fullerene molecules are those in which all the pentagons are isolated. Zhang et al. \cite{Zhang1992a} stated that fullerene stability is related to the ratio between the number of pentagons and the number of carbon atoms. Xia and Wei  \cite{xia2015persistent} proposed that the stability of fullerene depends on the average number of hexagons per atom. However, these theories all focused on the pentagon and hexagon information. More specifically, they use topological information to reveal the stability of fullerene. In contrast, we believe that the non-harmonic persistent spectra can also model the structure-function relationship of fullerenes.  We hypothesize that the non-harmonic persistent spectra of $\mathcal{L}_0^{r+0}$ matrices are powerful enough to model the stability of fullerene molecules. To verify our hypothesis, we compute the summation, mean, maximal, standard deviation, variance of its eigenvalues, and $(\tilde{\lambda}_2)_0^{r+0}$ of the persistent spectra of $\mathcal{L}_0^{r+0}$ over various filtration radii $r$. We depict a plot with the horizontal axis represents radius $r$ and the vertical   axis represents the particular spectrum value, which is actually the same as  Figure  \autoref{fig:C20C60}. Then we  define the area under the plot of spectra with a negative sign as
        \begin{equation}\label{equ:intergration}
            A_\alpha = -\sum_{i=1} \Lambda_i^\alpha \delta r,
        \end{equation}
       where $\delta r$ is the radius grid spacing, in \autoref{fig:C20C60}, $\delta r = \SI{0.01}{\angstrom}$. Here, $\alpha=$ Sum, Avg, Max, Std, Var, Sec is the type index and thus, $\Lambda_i^\alpha$ represent the summation, mean, maximal, standard deviation, variance, and the smallest non-zero eigenvalue $(\tilde{\lambda}_2)_0^{r+0}$ of $\mathcal{L}_0^{r+0}$ at $i$-th radius step, respectively. The right chart in \autoref{fig:Heat comparison} describes the area under the plot of spectra and closely resembles that of the heat of formation energy. We can see that generally the left chart and the middle chart show the same pattern. The integration of $(\tilde{\lambda}_2)_0^{r+0}$  decreases as the number of carbon atoms increases. However, the structural data we used might not be the same ground-state data as in Ref. \cite{Zhang1992a}, which results in C$_{36}$ do not match the corresponding energy perfectly. Limited by the availability of the ground-state structural data, we are not able to analyze the full set of the fullerene family.

        To quantitatively validate our model, we apply one of the simplest machine learning algorithms,  linear least-squares method, to predict the heat of formation energy. The Pearson correlation coefficient is defined as
        \begin{equation}\label{equ:pcc}
            C_c^{\alpha} = \dfrac{\displaystyle{\sum_{i=1}^{N}}(A_{\alpha}^i-\bar{A}_{\alpha})(E_i-\bar{E})}{\left[\displaystyle{\sum_{i=1}^{N}}(A_{\alpha}^i-\bar{A}_{\alpha})^2 \displaystyle{\sum_{i=1}^{N}}(E_i-\bar{E})^2\right]^{\frac{1}{2}}}
        \end{equation}
         where $A_{\alpha}^i$ represents the theoretically predicted energy of the $i$-th fullerene molecule, $E_i$ represents the heat of formation energy of the $i$-th fullerene molecule, and $\bar{A}_{\alpha}$ and $\bar{E}$ are the corresponding mean values. When $\alpha=\mbox{Max}$, the Pearson correlation coefficient is around $0.986$.   The right chart of \autoref{fig:Heat comparison} plots the correlation between predicted energies and the heat of formation energy of the fullerene molecules computed from quantum mechanics \cite{Zhang1992a}. These results agree very well.
        \begin{table}[H]
            \centering
            \setlength\tabcolsep{5pt}
            \captionsetup{margin=0.9cm}
            \caption{The heat of formation energy of fullerenes \cite{Zhang1992a} and its corresponding predicted energies with $\alpha = \mbox{Max}$. The unit is EV/atom. }
            \begin{tabular}{lcccccccc}
            \hline
             Fullerene type              & C$_{20}$  & C$_{24}$  & C$_{26}$   & C$_{30}$  & C$_{32}$  & C$_{36}$  & C$_{50}$  & C$_{60}$  \\ \hline
             Heat of formation energy    & $1.180$   & $1.050$   & $0.989$      & $0.850$   & $0.781$   & $0.706$   & $0.509$   & $0.401$   \\
             Predicted energy     & $1.138$   & $1.050$   & $0.964$           & $0.821$   & $0.857$   & $0.766$   & $0.474$   & $0.391$    \\ \hline
            \end{tabular}
            \label{table:energy}
        \end{table}

     The right chart of  \autoref{fig:Heat comparison} illustrates the fitting results under different type index $\alpha$.   \autoref{table:pcc} lists the correlation coefficient under different type index $\alpha$. The highest correlation coefficient is close to unity ($0.986$) obtained with $\alpha = \mbox{Max}$. The lowest correlation coefficient is $0.942$ with $\alpha = \mbox{Sum}$. We can see that all the correlation coefficients are close to unity, which verifies our hypothesis that the non-harmonic spectra of $\mathcal{L}_0^{r+0}$ have the capacity of modeling the stability of fullerene molecules. Although we ignore the topological information (Betti numbers), our persistent spectral theory still works extremely well only with non-harmonic spectra, which means our persistent spectral theory is a powerful tool for quantitative data analysis and prediction.


        \begin{table}[H]
            \centering
            \setlength\tabcolsep{9.5pt}
            \captionsetup{margin=0.9cm}
            \caption{The correlation coefficients under different type index $\alpha$.}
            \begin{tabular}{lcccccc}
            \hline
             Type index                 & Sum       & Avg       & Max       & Std       & Var          & Sec   \\ \hline
             Correlation coefficient    & $0.942$   & $0.985$   & $0.986$   & $0.969$   & $0.977$      & $0.981$   \\ \hline
            \end{tabular}
            \label{table:pcc}
        \end{table}

    \subsection{Protein flexibility analysis}

 As clarified earlier, the number of zero eigenvalues of $p$-persistent $q$-Laplacian matrix ($p$-persistent $q$th Betti number) can also be derived from persistent homology. Persistent homology has been used to model fullerene stability \cite{xia2015persistent}. In this section, we further illustrate the applicability of present persistent spectral theory by a case that non-harmonic persistent spectra offer a unique theoretical model whereas it may be difficult to come up with a suitable persistent homology model for this problem. 

    The protein flexibility is known to correlate with a wide variety of protein functions. It can be modeled by the beta factors or B-factors, which are also called Debye-Waller factors. B-factors  are a measure of the atomic mean-square displacement or uncertainty in the X-ray scattering structure determination. Therefore, understanding the protein structure, flexibility, and function via the accurate protein B-factor prediction is a vital task in computational biophysics \cite{bramer2018blind}. Over the past few years, quite many methods are developed to predict protein B-factors, such as    GNM, \cite{Bahar1997}, ANM \cite{Atilgan2001},     FRI,   \cite{opron2014fast, xia2013multiscale} and MWCG \cite{xia2015multiscale, bramer2018blind}.  However, all of the aforementioned methods are based on a particular matrix derived from the graph network which is constructed using alpha carbon as nodes and connections between nodes as edges. In this section, we apply our persistent spectral theory to create richer geometric information in B-factor prediction.

     To illustrate our method, we consider protein $2$Y$7$L whose total number of residues is $N=319$. In this work, we employ the coarse-grained C$_\alpha$ representation of $2$Y$7$L. Therefore,  $319$ particles are taken into consideration in protein $2$Y$7$L. Similarly, like in the previous application of fullerene structure analysis, we treat each C$_\alpha$ atom as a $0$-simplex at the initial setup and assign it a solid ball with a radius of $r$. By varying the filtration parameter $r$, we can obtain a family of $\mathcal{L}_0^{r+0}$. For each matrix $\mathcal{L}_0^{r+0}$, its corresponding ordered spectrum is given by
    \[
    (\lambda_1)_0^{r+0}, (\lambda_2)_0^{r+0}, \cdots, (\lambda_N)_0^{r+0}.  
    \]
    Suppose the number of zero eigenvalues  is $m$,  then, we have  $\beta_0^{r+0} = m$. Since $\mathcal{L}_0^{r+0}$ is symmetric, then eigenvectors of $\mathcal{L}_0^{r+0}$ corresponding to different eigenvalues must be orthogonal to each other. The Moore-Penrose inverse of $\mathcal{L}_0^{r+0}$ can be calculated by the non-harmonic spectra of $\mathcal{L}_0^{r+0}$:
    \[
    (\mathcal{L}_0^{r+0})^{-1} = \sum_{k=m+1}^{N} \frac{1}{(\lambda_k)_0^{r+0}}[(u_k)_0^{r+0} ((u_k)_0^{r+0})^{T}],
    \]
    where $T$ is the transpose and $(u_k)_0^{r+0}$ is the $k$th eigenvector of $\mathcal{L}_0^{r+0}$. The modeling of $i$th B-factor of $2$Y$7$L at filtration parameter $r$ can be expressed as
    \[
    B_i^r = (\mathcal{L}_0^{r+0})^{-1}_{ii}, \forall i = 1,2,\cdots,N,
    \]
    and the final model of $i$th B-factor of $2$Y$7$L is given by 
    \[
    B^{\text{PST}}_i = \sum_{r}w_r B_i^r + w_0,  \forall i = 1,2,\cdots,N,
    \]
    where $w_r$ and $w_0$ are  fitting parameters which can be derived by linearly fitting B-factors from experimental data $B^{\text{Exp}}$. Consider the filtration radius from $2$ to $12$ with the grid spacing of $1$, then totally $11$ different $\mathcal{L}_0^{r+0}$ are created. By calculating all the non-harmonic spectra together with their eigenvectors, $11$ Moore-Penrose inverse matrices $(\mathcal{L}_0^{r+0})^{-1}$ can be constructed. Therefore, the predicted $i$th B-factor is 
    \[
    B^{\text{PST}}_i = \sum_{r=2}^{12}w_r B_i^r + w_0.
    \]
    The specific values of $w_r$ and $w_0$ can be found in \autoref{table:weights1} and \autoref{table:weights2} of Appendix \autoref{app:application}. \autoref{fig:2Y7L Combine} (c) shows that  the prediction B-factors are in an excellent agreement with   the experimental B-factors of protein $2$Y$7$L. The Pearson correlation coefficient is $0.925$ \footnote{We carry out feature scaling to make sure all $B_i^r$ are on a similar scale.}.

    \begin{figure}[H]
        \centering
        \captionsetup{margin=0.9cm}
        \includegraphics[scale=0.5]{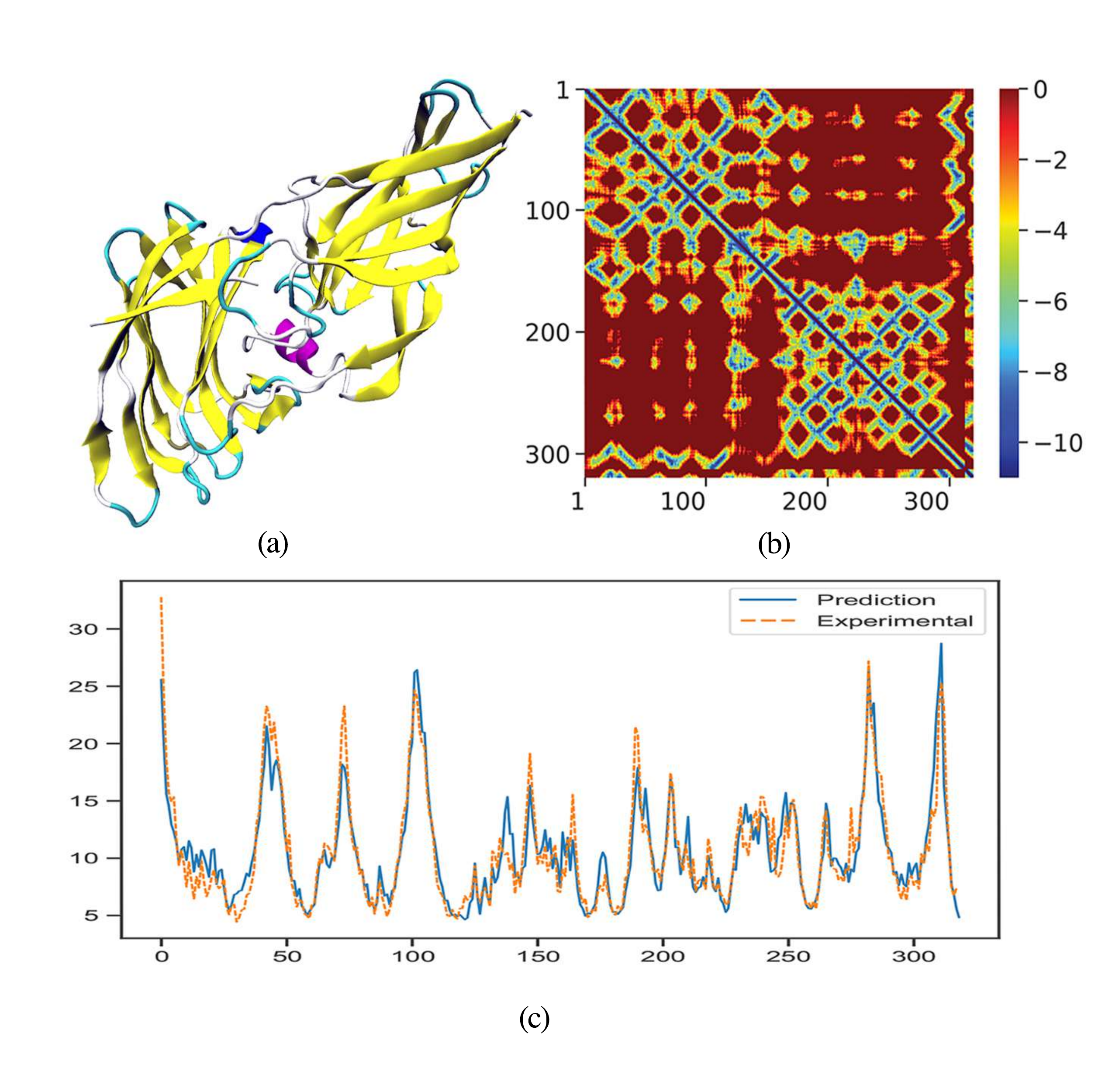}
        \caption{ Illustration of persistent spectral prediction of protein B-factors. (a) Plot of the secondary structure of protein $2$Y$7$L. (b) Accumulated persistent Laplacian matrix (For clarity, the diagonal terms are set to 0.). Note that the accumulated persistent Laplacian matrix maps out the detailed distance between each pair of residues. (c) Comparison of experimental B-factors and those predicted by PST for protein $2$Y$7$L.}
        \label{fig:2Y7L Combine}
    \end{figure}

   This example shows that our persistent spectral theory can be used beyond the persistent homology analysis. The number of zero eigenvalues of $0$-persistent $q$-combinatorial Laplacian matrices fully recover the persistent barcode or persistent diagram of persistent homology. Additional spectral information  from non-harmonic persistent spectra and persistent eigenvectors provides valuable information for data modeling, analysis, and prediction.

\section{Conclusion}
Spectral graph theory is a powerful tool for data analysis due to its ability to extract geometric  and topological information.  However, its performance can be quite limited for various reasons. One of them is that the current spectral graph theory does not provide a multiscale analysis. Motivated by persistent homology and multiscale graphs, we introduce persistent spectral theory as a unified paradigm to unveil both topological persistence and geometric shape from high-dimensional datasets. 

  For a point set $V\subset \mathbb{R}^n$ without additional structures, we construct a filtration using an $(n-1)$-sphere of a varying radius $r$ centered at  each point. A series of  persistent  combinatorial Laplacian matrices are induced by the filtration. It is noted that our harmonic persistent spectra (i.e., zero eigenvalues) fully recover the persistent barcode or persistent diagram of persistent homology. Specifically, the numbers of zero eigenvalues of  persistent $q$-combinatorial Laplacian matrices are the $q$-dimensional  persistent  Betti numbers for the same filtration given filtration. However, additional valuable spectral information is generated from the non-harmonic persistent spectra. In this work, in addition to persistent  Betti numbers and the smallest non-zero eigenvalues, five statistic values, namely, sum, mean, maximum, standard deviation, and variance, are also constructed for data analysis. We use a few simple two-dimensional (2D) and three-dimensional (3D) structures to carry out the proof of principle analysis of the persistent spectral theory. The detailed structural information can be incorporated into the persistent spectra of. For instant, for the benzene molecule, the approximate C-C bond and C-H bond length can be intuitively read from the plot of  the $0$-dimensional  persistent  Betti numbers. Moreover, persistent spectral theory also has the capacity to accurately predict the heat of formation energy of small fullerene molecules. We use the area under the plot of the persistent spectra to model fullerene stability and apply the linear least-squares method to fit our prediction with the heat of formation energy. The resulting correlation coefficient is close to $1$, which shows that our persistent spectral theory has an excellent performance on molecular data. Furthermore, we have applied our persistent spectral theory to the protein B-factor prediction. In this case, persistent homology does not give a straightforward model. This example shows that the additional non-harmonic persistent spectral information provides a powerful tool for dealing with molecular data.  

It is pointed out that the proposed persistent spectral analysis can be paired with advanced machine learning algorithms, including various deep learning methods, for a wide variety of applications in data science.  In particular, the further construction of element-specific persistent spectral theory and its application to protein-ligand binding affinity prediction and computer-aided drug design will be reported elsewhere.

\section{Acknowledgments}
This work was supported in part by NSF Grants DMS1721024, DMS1761320, IIS1900473, NIH grants GM126189 and GM129004,  Bristol-Myers Squibb,  and Pfizer.
RW thanks Dr. Jiahui Chen and Dr. Zixuan Cang for useful discussions.


\bibliographystyle{abbrv}
\bibliography{PST}

\newpage
\begin{appendices}
  \renewcommand\thetable{\thesection\arabic{table}}
  \renewcommand\thefigure{\thesection\arabic{figure}}
	
  \section{Persistence Homology} \label{app:PH}
  
		Persistence Homology is an algebraic topology-based  method for the multiscale analysis of the topological invariants of   functions and datasets.  It has been widely applied in the field of topological  data analysis. We provide a brief introduction to persistent homology and the interested readers are referred to the literature \cite{Fallis2013,edelsbrunner2010computational} for more detail.

        \subsection{Homology}
        For a topological space $X$,   a sequences of complexes $C_0(X), C_1(X), \cdots$   describes different dimensional information of the topological space $X$, which are connected by homomorphisms (or boundary operators)  $\partial_{k}: C_{k} \longrightarrow C_{k-1}$ such that $\im \partial_{k} \subseteq \ker \partial_{k-1}$, i.e.,  $ \partial_{k-1}   \partial_{k}=0$. With a $k$-simplex $\sigma_{k} = [v_0, \cdots, v_k]$ where $v_i$ are all the vertices of $\sigma_k$, $ \partial_k\sigma_k$ can be given by  a formal sum with coefficients in the $\mathbb{Z}_2$  field
        \begin{equation}
            \partial_k \sigma_k = \sum_{i=0}^{k}\sigma^i_{k-1},
        \end{equation}
        where $\sigma_i^{k-1}$ is the $(k-1)$-simplex with its $i$th vertex $v_i$ being omitted.
        The algebraic construction to connect a sequence of  complexes by boundary maps is called a chain complex
        \[
            \cdots \stackrel{\partial_{i+1}}\longrightarrow C_i(X) \stackrel{\partial_i}\longrightarrow C_{i-1}(X) \stackrel{\partial_{i-1}}\longrightarrow \cdots \stackrel{\partial_2} \longrightarrow C_{1}(X) \stackrel{\partial_{1}}\longrightarrow C_0(X) \stackrel{\partial_0} \longrightarrow 0
        \]
        and the $k$th homology group is the quotient group defined by
        \begin{equation}
            H_k = \ker \partial_k / \im \partial_{k+1}.
        \end{equation}
        By studying  homology groups, one can derive homological properties of the space $K$. The Betti numbers are defined by the ranks of $k$th homology group $H_k$ which counts  $k-$dimensional holes, especially, ${\rm rank}(H_0)$ reflects the number of connected components, ${\rm rank}(H_1)$ reflects the number of loops, and ${\rm rank}(H_2)$ reveals the number of voids or cavities. However, ${\rm rank}(H_k)$ only allows us to express the topological information for a specific setup. Persistent homology is devised to track the multiscale topological information over different scales along a filtration.

        \subsection{Persistent homology}
				
        A filtration of a topology space $K$ is a sequence of sub-spaces $(K_t)_{t=0}^m$ of $K$ such that
        \begin{equation}
            \emptyset = K_0 \subseteq K_1 \subseteq K_2 \subseteq \cdots \subseteq K_m = K.
        \end{equation}
        A sequence of chain complexes induced by the filtration is defined as
        \begin{equation}
        \left.\begin{array}{ccccccccc}
            \cdots & \stackrel{\partial_{3}}\longrightarrow & C_2^1 & \stackrel{\partial_{2}}\longrightarrow & C_1^1 & \stackrel{\partial_{1}}\longrightarrow & C_0^1 & \stackrel{\partial_{0}}\longrightarrow & 0 \\
             &  & \rotatebox{-90}{$\longrightarrow$} &  & \rotatebox{-90}{$\longrightarrow$} &  & \rotatebox{-90}{$\longrightarrow$} &  &  \\
            \cdots & \stackrel{\partial_{3}}\longrightarrow & C_2^2 & \stackrel{\partial_{2}}\longrightarrow & C_1^2 & \stackrel{\partial_{1}}\longrightarrow & C_0^2 & \stackrel{\partial_{0}}\longrightarrow & 0 \\
             &  & \rotatebox{-90}{$\dashrightarrow$} &  & \rotatebox{-90}{$\dashrightarrow$} &  & \rotatebox{-90}{$\dashrightarrow$} &  &  \\
            \cdots & \stackrel{\partial_{3}}\longrightarrow & C_2^m & \stackrel{\partial_{2}}\longrightarrow & C_1^m & \stackrel{\partial_{1}}\longrightarrow & C_0^m & \stackrel{\partial_{0}}\longrightarrow & 0
        \end{array}\right.
        \end{equation}
        with $C_k^t \coloneqq C_k(K^t)$ and $\downarrow$ denotes the inclusion \cite{hofer2017deep}. The $p$-persistent $k$th homology group of $K^t$ is defined as
        \begin{equation}
            H_k^p(K^t) = \ker \partial_{k}(K^t) / (\im \partial_{k+1}(K^{t+p}) \cap \ker \partial_{k}(K^t)),
        \end{equation}
         Intuitively, this homology group records the homology classes of $K^t$ that are persistent at least until $K^{t+p}$. When $k=0$, the rank of $H_k^p(K^t)$ reveals the number of connected components in $K^t$.

  \section{Additional Laplacian matrices and their properties }\label{app:examples}
	
  In this section, we give a further description of  additional  boundary and Laplacian matrices and their properties involved in  the filtration process  in \autoref{fig:filtration}.

        \begin{table}[H]
            \centering
            \setlength\tabcolsep{14pt}
            \captionsetup{margin=0.9cm}
            \caption{$K_1 \to K_1$}
            \begin{tabular}{c|ccc}
            \hline
            $q$  & $q=0$ & $q=1$ & $q=2$ \\ \hline \\
            $\mathcal{B}_{q+1}^{1+0}$  & / & / &  /    \\  \\
            $\mathcal{B}_{q}^{1}$  & $\begin{array}{@{}r@{}c@{}c@{}c@{}c@{}c@{}l@{}}
                    & 0 & 1 & 2 & 3 & 4   \\
                   \left.\begin{array}{c}

                    \end{array}\right[
                    & \begin{array}{c}  0  \end{array}
                    & \begin{array}{c}  0  \end{array}
                    & \begin{array}{c}  0  \end{array}
                    & \begin{array}{c}  0  \end{array}
                    & \begin{array}{c}  0  \end{array}
                    & \left]\begin{array}{c} \\
                    \end{array}\right.
                \end{array}$ & / &  /    \\  \\
            $\mathcal{L}_{q}^{1+0}$  & $\left[\begin{array}{ccccc}
                 0  &   0  &   0  &   0  &  0  \\
                 0  &   0  &   0  &   0  &  0  \\
                 0  &   0  &   0  &   0  &  0  \\
                 0  &   0  &   0  &   0  &  0  \\
                 0  &   0  &   0  &   0  &  0
            \end{array}\right]$  & / &  /    \\  \\
            $\beta_{q}^{1+0}$                        & 5               & / & /    \\  \\
            $\text{dim}(\mathcal{L}_{q}^{1+0})$      & 5               & / & /    \\  \\
            $\text{rank}(\mathcal{L}_{q}^{1+0})$     & 0               & / & /    \\  \\
            $\text{nullity}(\mathcal{L}_{q}^{1+0})$  & 5               & / & /  \\  \\
            $\text{Spectra}(\mathcal{L}_{q}^{1+0})$    & $\{0,0,0,0,0\}$   & / & /  \\ \hline
            \end{tabular}
            \\
            \label{table:ex 1}
        \end{table}


        \begin{table}[H]
            \centering
            \setlength\tabcolsep{14pt}
            \captionsetup{margin=0.9cm}
            \caption{$K_2 \to K_2$}
            \begin{tabular}{c|ccc}
            \hline
            $q$  & $q=0$ & $q=1$ & $q=2$ \\ \hline \\
            $\mathcal{B}_{q+1}^{2+0}$  & $\begin{array}{@{}r@{}c@{}l@{}}
                    & 01    \\
                   \left.\begin{array}{c}
                    0 \\
                    1 \\
                    2 \\
                    3 \\
                    4
                    \end{array}\right[
                    & \begin{array}{c} -1 \\  1  \\  0  \\  0 \\ 0 \end{array}
                    & \left]\begin{array}{c} \\ \\  \\ \\ \\
                    \end{array}\right.
                \end{array}$ & / &  /    \\  \\
            $\mathcal{B}_{q}^{2}$  & $\begin{array}{@{}r@{}c@{}c@{}c@{}c@{}c@{}l@{}}
                    & 0 & 1 & 2 & 3 & 4   \\
                   \left.\begin{array}{c}

                    \end{array}\right[
                    & \begin{array}{c}  0  \end{array}
                    & \begin{array}{c}  0  \end{array}
                    & \begin{array}{c}  0  \end{array}
                    & \begin{array}{c}  0  \end{array}
                    & \begin{array}{c}  0  \end{array}
                    & \left]\begin{array}{c}
                    \end{array}\right.
                \end{array}$ & $\begin{array}{@{}r@{}c@{}l@{}}
                    & 01    \\
                   \left.\begin{array}{c}
                    0 \\
                    1 \\
                    2 \\
                    3 \\
                    4
                    \end{array}\right[
                    & \begin{array}{c} -1 \\  1  \\  0  \\  0 \\ 0 \end{array}
                    & \left]\begin{array}{c} \\ \\  \\ \\ \\
                    \end{array}\right.
                \end{array}$ &  /    \\  \\
            $\mathcal{L}_{q}^{2+0}$  & $\left[\begin{array}{ccccc}
                 1  &  -1  &   0  &   0  &  0  \\
                -1  &   1  &   0  &   0  &  0  \\
                 0  &   0  &   0  &   0  &  0  \\
                 0  &   0  &   0  &   0  &  0  \\
                 0  &   0  &   0  &   0  &  0
            \end{array}\right]$  & [2] &  /    \\  \\
            $\beta_{q}^{2+0}$                        & 4               & 0 & /    \\  \\
            $\text{dim}(\mathcal{L}_{q}^{2+0})$      & 5               & 1 & /    \\  \\
            $\text{rank}(\mathcal{L}_{q}^{2+0})$     & 1               & 1 & /    \\  \\
            $\text{nullity}(\mathcal{L}_{q}^{2+0})$  & 4               & 0 & /  \\  \\
            $\text{Spectra}(\mathcal{L}_{q}^{2+0})$    & $\{0,0,0,0,2\}$   & 2 & /  \\ \hline
            \end{tabular}
            \\
            \label{table:ex 2}
        \end{table}


        \begin{table}[H]
            \centering
            \setlength\tabcolsep{8pt}
            \captionsetup{margin=0.9cm}
            \caption{$K_3 \to K_3$}
            \begin{tabular}{c|ccc}
            \hline
            $q$  & $q=0$ & $q=1$ & $q=2$ \\ \hline \\
            $\mathcal{B}_{q+1}^{3+0}$  & $\begin{array}{@{}r@{}c@{}c@{}c@{}c@{}l@{}}
                    & 01 & 12 & 23 & 03   \\
                   \left.\begin{array}{c}
                    0 \\
                    1 \\
                    2 \\
                    3 \\
                    4
                    \end{array}\right[
                    & \begin{array}{c} -1 \\  1  \\  0  \\  0 \\ 0 \end{array}
                    & \begin{array}{c}  0 \\ -1  \\  1  \\  0 \\ 0 \end{array}
                    & \begin{array}{c}  0 \\  0  \\ -1  \\  1 \\ 0 \end{array}
                    & \begin{array}{c} -1 \\  0  \\  0  \\  1 \\ 0 \end{array}
                    & \left]\begin{array}{c} \\ \\  \\ \\ \\
                    \end{array}\right.
                \end{array}$ & / &  /    \\  \\
            $\mathcal{B}_{q}^{3}$  & $\begin{array}{@{}r@{}c@{}c@{}c@{}c@{}c@{}l@{}}
                    & 0 & 1 & 2 & 3 & 4   \\
                   \left.\begin{array}{c}

                    \end{array}\right[
                    & \begin{array}{c}  0  \end{array}
                    & \begin{array}{c}  0  \end{array}
                    & \begin{array}{c}  0  \end{array}
                    & \begin{array}{c}  0  \end{array}
                    & \begin{array}{c}  0  \end{array}
                    & \left]\begin{array}{c}
                    \end{array}\right.
                \end{array}$ & $\begin{array}{@{}r@{}c@{}c@{}c@{}c@{}l@{}}
                    & 01 & 12 & 23 & 03   \\
                   \left.\begin{array}{c}
                    0 \\
                    1 \\
                    2 \\
                    3 \\
                    4
                    \end{array}\right[
                    & \begin{array}{c} -1 \\  1  \\  0  \\  0 \\ 0 \end{array}
                    & \begin{array}{c}  0 \\ -1  \\  1  \\  0 \\ 0 \end{array}
                    & \begin{array}{c}  0 \\  0  \\ -1  \\  1 \\ 0 \end{array}
                    & \begin{array}{c} -1 \\  0  \\  0  \\  1 \\ 0 \end{array}
                    & \left]\begin{array}{c} \\ \\  \\ \\ \\
                    \end{array}\right.
                \end{array}$ &  /    \\  \\
            $\mathcal{L}_{q}^{3+0}$  & $\left[\begin{array}{ccccc}
                 2  &  -1  &   0  &  -1  &  0  \\
                -1  &   2  &  -1  &   0  &  0  \\
                 0  &  -1  &   2  &  -1  &  0  \\
                -1  &   0  &  -1  &   2  &  0  \\
                 0  &   0  &   0  &   0  &  0
            \end{array}\right]$  & $\left[\begin{array}{cccc}
                 2  &  -1  &   0  &   1   \\
                -1  &   2  &  -1  &   0   \\
                 0  &  -1  &   2  &   1   \\
                 1  &   0  &   1  &   2   \\

            \end{array}\right]$ &  /    \\  \\
            $\beta_{q}^{3+0}$                        & 2               & 1 & /    \\  \\
            $\text{dim}(\mathcal{L}_{q}^{3+0})$      & 5               & 4 & /    \\  \\
            $\text{rank}(\mathcal{L}_{q}^{3+0})$     & 3               & 3 & /    \\  \\
            $\text{nullity}(\mathcal{L}_{q}^{3+0})$  & 2               & 1 & /  \\  \\
            $\text{Spectra}(\mathcal{L}_{q}^{3+0})$    & $\{0,0,2,2,4\}$   & $\{0,2,2,4\}$ & /  \\ \hline
            \end{tabular}
            \\
            \label{table:ex 3}
        \end{table}


        \begin{table}[H]
            \centering
            \setlength\tabcolsep{0pt}
            \captionsetup{margin=0.9cm}
            \caption{$K_5 \to K_5$}
            \begin{tabular}{c|ccc}
            \hline
            $q$  & $q=0$ & $q=1$ & $q=2$ \\ \hline \\
            $\mathcal{B}_{q+1}^{5+0}$  & $\begin{array}{@{}r@{}c@{}c@{}c@{}c@{}c@{}c@{}l@{}}
                    & 01 & 12 & 23 & 03 & 24 & 02  \\
                   \left.\begin{array}{c}
                    0 \\
                    1 \\
                    2 \\
                    3 \\
                    4
                    \end{array}\right[
                    & \begin{array}{c} -1 \\  1  \\  0  \\  0 \\ 0 \end{array}
                    & \begin{array}{c}  0 \\ -1  \\  1  \\  0 \\ 0 \end{array}
                    & \begin{array}{c}  0 \\  0  \\ -1  \\  1 \\ 0 \end{array}
                    & \begin{array}{c} -1 \\  0  \\  0  \\  1 \\ 0 \end{array}
                    & \begin{array}{c}  0 \\  0  \\ -1  \\  0 \\ 1 \end{array}
                    & \begin{array}{c} -1 \\  0  \\  1  \\  0 \\ 0 \end{array}
                    & \left]\begin{array}{c} \\ \\  \\ \\ \\
                    \end{array}\right.
                \end{array}$ & $\begin{array}{@{}r@{}c@{}c@{}l@{}}
                    & 012 & 023  \\
                   \left.\begin{array}{c}
                    01 \\
                    12 \\
                    23 \\
                    03 \\
                    24 \\
                    02
                    \end{array}\right[
                    & \begin{array}{c}  1 \\  1  \\  0  \\  0 \\ 0 \\ -1  \end{array}
                    & \begin{array}{c}  0 \\  0  \\  1  \\ -1 \\ 0 \\  1  \end{array}
                    & \left]\begin{array}{c} \\ \\  \\ \\ \\ \\
                    \end{array}\right.
                \end{array}$ &   $\begin{array}{@{}r@{}c@{}l@{}}
                    & 0123  \\
                   \left.\begin{array}{c}
                    012 \\
                    023
                    \end{array}\right[
                    & \begin{array}{c}  -1 \\  1   \end{array}
                    & \left]\begin{array}{c} \\ \\
                    \end{array}\right.
                \end{array}$    \\  \\
            $\mathcal{B}_{q}^{5}$  & $\begin{array}{@{}r@{}c@{}c@{}c@{}c@{}c@{}l@{}}
                    & 0 & 1 & 2 & 3 & 4   \\
                   \left.\begin{array}{c}

                    \end{array}\right[
                    & \begin{array}{c}  0  \end{array}
                    & \begin{array}{c}  0  \end{array}
                    & \begin{array}{c}  0  \end{array}
                    & \begin{array}{c}  0  \end{array}
                    & \begin{array}{c}  0  \end{array}
                    & \left]\begin{array}{c}
                    \end{array}\right.
                \end{array}$ & $\begin{array}{@{}r@{}c@{}c@{}c@{}c@{}c@{}c@{}l@{}}
                    & 01 & 12 & 23 & 03 & 24 & 02  \\
                   \left.\begin{array}{c}
                    0 \\
                    1 \\
                    2 \\
                    3 \\
                    4
                    \end{array}\right[
                    & \begin{array}{c} -1 \\  1  \\  0  \\  0 \\ 0 \end{array}
                    & \begin{array}{c}  0 \\ -1  \\  1  \\  0 \\ 0 \end{array}
                    & \begin{array}{c}  0 \\  0  \\ -1  \\  1 \\ 0 \end{array}
                    & \begin{array}{c} -1 \\  0  \\  0  \\  1 \\ 0 \end{array}
                    & \begin{array}{c}  0 \\  0  \\ -1  \\  0 \\ 1 \end{array}
                    & \begin{array}{c} -1 \\  0  \\  1  \\  0 \\ 0 \end{array}
                    & \left]\begin{array}{c} \\ \\  \\ \\ \\
                    \end{array}\right.
                \end{array}$ &   $\begin{array}{@{}r@{}c@{}c@{}l@{}}
                    & 012 & 023   \\
                   \left.\begin{array}{c}
                    01 \\
                    12 \\
                    23 \\
                    03 \\
                    24 \\
                    02
                    \end{array}\right[
                    & \begin{array}{c}  1 \\  1  \\  0  \\  0 \\ 0 \\ -1 \end{array}
                    & \begin{array}{c}  0 \\  0  \\  1  \\ -1 \\ 0 \\  1 \end{array}
                    & \left]\begin{array}{c} \\ \\  \\ \\ \\ \\
                    \end{array}\right.
                \end{array}$    \\  \\
            $\mathcal{L}_{q}^{5+0}$  & $\left[\begin{array}{ccccc}
                 3  &  -1  &  -1  &  -1  &  0  \\
                -1  &   2  &  -1  &   0  &  0  \\
                -1  &  -1  &   4  &  -1  & -1  \\
                -1  &   0  &  -1  &   2  &  0  \\
                 0  &   0  &  -1  &   0  &  1
            \end{array}\right]$  & $\left[\begin{array}{cccccc}
                 3  &   0  &   0  &   1  &  0  &  0   \\
                 0  &   3  &  -1  &   0  & -1  &  0   \\
                 0  &  -1  &   3  &   0  &  1  &  0   \\
                 1  &   0  &   0  &   3  &  0  &  0   \\
                 0  &  -1  &   1  &   0  &  2  & -1   \\
                 0  &   0  &   0  &   0  & -1  &  4
            \end{array}\right]$ &  $\left[\begin{array}{cc}
                 4  &   0  \\
                 0  &   4
            \end{array}\right]$   \\  \\
            $\beta_{q}^{5+0}$                        & 1               & 0 & 0    \\  \\
            $\text{dim}(\mathcal{L}_{q}^{5+0})$      & 5               & 6 & 2    \\  \\
            $\text{rank}(\mathcal{L}_{q}^{5+0})$     & 4               & 6 & 2    \\  \\
            $\text{nullity}(\mathcal{L}_{q}^{5+0})$  & 1               & 0 & 0  \\  \\
            $\text{Spectra}(\mathcal{L}_{q}^{5+0})$    & $\{0,1,2,4,5\}$   & $\{1,2,2,4,4,5\}$ & $\{4,4\}$  \\ \hline
            \end{tabular}
            \\
            \label{table:ex 4}
        \end{table}


        \begin{table}[H]
            \centering
            \setlength\tabcolsep{10pt}
            \captionsetup{margin=0.9cm}
            \caption{$K_1 \to K_2$}
            \begin{tabular}{c|ccc}
            \hline
            $q$  & $q=0$ & $q=1$ & $q=2$ \\ \hline \\
            $\mathcal{B}_{q+1}^{1+1}$  & $\begin{array}{@{}r@{}c@{}l@{}}
                    & 01    \\
                   \left.\begin{array}{c}
                    0 \\
                    1 \\
                    2 \\
                    3 \\
                    4
                    \end{array}\right[
                    & \begin{array}{c} -1 \\  1  \\  0  \\  0 \\ 0 \end{array}
                    & \left]\begin{array}{c} \\ \\  \\ \\ \\
                    \end{array}\right.
                \end{array}$ & / &  /    \\  \\
            $\mathcal{B}_{q}^{1}$  & $\begin{array}{@{}r@{}c@{}c@{}c@{}c@{}c@{}l@{}}
                    & 0 & 1 & 2 & 3 & 4   \\
                   \left.\begin{array}{c}

                    \end{array}\right[
                    & \begin{array}{c}  0  \end{array}
                    & \begin{array}{c}  0  \end{array}
                    & \begin{array}{c}  0  \end{array}
                    & \begin{array}{c}  0  \end{array}
                    & \begin{array}{c}  0  \end{array}
                    & \left]\begin{array}{c}
                    \end{array}\right.
                \end{array}$ & / &  /    \\  \\
            $\mathcal{L}_{q}^{1+1}$  & $\left[\begin{array}{ccccc}
                 1  &  -1  &   0  &   0  &  0  \\
                -1  &   1  &   0  &   0  &  0  \\
                 0  &   0  &   0  &   0  &  0  \\
                 0  &   0  &   0  &   0  &  0  \\
                 0  &   0  &   0  &   0  &  0
            \end{array}\right]$  & / &  /    \\  \\
            $\beta_{q}^{1+1}$                        & 4               & / & /    \\  \\
            $\text{dim}(\mathcal{L}_{q}^{1+1})$      & 5               & / & /    \\  \\
            $\text{rank}(\mathcal{L}_{q}^{1+1})$     & 1               & / & /    \\  \\
            $\text{nullity}(\mathcal{L}_{q}^{1+1})$  & 4               & / & /  \\  \\
            $\text{Spectra}(\mathcal{L}_{q}^{1+1})$    & $\{0,0,0,0,2\}$   & / & /  \\ \hline
            \end{tabular}
            \\
            \label{table:ex 5}
        \end{table}


        \begin{table}[H]
            \centering
            \setlength\tabcolsep{10pt}
            \captionsetup{margin=0.9cm}
            \caption{$K_1 \to K_4$}
            \begin{tabular}{c|ccc}
            \hline
            $q$  & $q=0$ & $q=1$ & $q=2$ \\ \hline \\
            $\mathcal{B}_{q+1}^{1+3}$  & $\begin{array}{@{}r@{}c@{}c@{}c@{}c@{}c@{}l@{}}
                    & 01 & 12 & 23 & 03 & 24   \\
                   \left.\begin{array}{c}
                    0 \\
                    1 \\
                    2 \\
                    3 \\
                    4
                    \end{array}\right[
                    & \begin{array}{c} -1 \\  1  \\  0  \\  0 \\ 0 \end{array}
                    & \begin{array}{c}  0 \\ -1  \\  1  \\  0 \\ 0 \end{array}
                    & \begin{array}{c}  0 \\  0  \\ -1  \\  1 \\ 0 \end{array}
                    & \begin{array}{c} -1 \\  0  \\  0  \\  1 \\ 0 \end{array}
                    & \begin{array}{c}  0 \\  0  \\ -1  \\  0 \\ 1 \end{array}
                    & \left]\begin{array}{c} \\ \\  \\ \\ \\
                    \end{array}\right.
                \end{array}$ & / &  /    \\  \\
            $\mathcal{B}_{q}^{1}$  &  $\begin{array}{@{}r@{}c@{}c@{}c@{}c@{}c@{}l@{}}
                    & 0 & 1 & 2 & 3 & 4   \\
                   \left.\begin{array}{c}

                    \end{array}\right[
                    & \begin{array}{c}  0  \end{array}
                    & \begin{array}{c}  0  \end{array}
                    & \begin{array}{c}  0  \end{array}
                    & \begin{array}{c}  0  \end{array}
                    & \begin{array}{c}  0  \end{array}
                    & \left]\begin{array}{c}
                    \end{array}\right.
                \end{array}$ & / &  /    \\  \\
            $\mathcal{L}_{q}^{1+3}$  & $\left[\begin{array}{ccccc}
                 2  &  -1  &   0  &  -1  &  0  \\
                -1  &   2  &  -1  &   0  &  0  \\
                 0  &  -1  &   3  &  -1  & -1  \\
                -1  &   0  &  -1  &   2  &  0  \\
                 0  &   0  &  -1  &   0  &  1
            \end{array}\right]$   & / &  /    \\  \\
            $\beta_{q}^{1+3}$                        & 1               & / & /    \\  \\
            $\text{dim}(\mathcal{L}_{q}^{1+3})$      & 5               & / & /    \\  \\
            $\text{rank}(\mathcal{L}_{q}^{1+3})$     & 4               & / & /    \\  \\
            $\text{nullity}(\mathcal{L}_{q}^{1+3})$  & 1               & / & /  \\  \\
            $\text{Spectra}(\mathcal{L}_{q}^{1+3})$    & $\{0,0.8299,2,2.6889,4.4812\}$   & / & /  \\ \hline
            \end{tabular}
            \\
            \label{table:ex 6}
        \end{table}


        \begin{table}[H]
            \centering
            \setlength\tabcolsep{10pt}
            \captionsetup{margin=0.9cm}
            \caption{$K_1 \to K_5$}
            \begin{tabular}{c|ccc}
            \hline
            $q$  & $q=0$ & $q=1$ & $q=2$ \\ \hline \\
            $\mathcal{B}_{q+1}^{1+4}$  & $\begin{array}{@{}r@{}c@{}c@{}c@{}c@{}c@{}c@{}l@{}}
                    & 01 & 12 & 23 & 03 & 24 & 02  \\
                   \left.\begin{array}{c}
                    0 \\
                    1 \\
                    2 \\
                    3 \\
                    4
                    \end{array}\right[
                    & \begin{array}{c} -1 \\  1  \\  0  \\  0 \\ 0 \end{array}
                    & \begin{array}{c}  0 \\ -1  \\  1  \\  0 \\ 0 \end{array}
                    & \begin{array}{c}  0 \\  0  \\ -1  \\  1 \\ 0 \end{array}
                    & \begin{array}{c} -1 \\  0  \\  0  \\  1 \\ 0 \end{array}
                    & \begin{array}{c}  0 \\  0  \\ -1  \\  0 \\ 1 \end{array}
                    & \begin{array}{c} -1 \\  0  \\  1  \\  0 \\ 0 \end{array}
                    & \left]\begin{array}{c} \\ \\  \\ \\ \\
                    \end{array}\right.
                \end{array}$ & / &  /    \\  \\
            $\mathcal{B}_{q}^{1}$  &  $\begin{array}{@{}r@{}c@{}c@{}c@{}c@{}c@{}l@{}}
                    & 0 & 1 & 2 & 3 & 4   \\
                   \left.\begin{array}{c}

                    \end{array}\right[
                    & \begin{array}{c}  0  \end{array}
                    & \begin{array}{c}  0  \end{array}
                    & \begin{array}{c}  0  \end{array}
                    & \begin{array}{c}  0  \end{array}
                    & \begin{array}{c}  0  \end{array}
                    & \left]\begin{array}{c}
                    \end{array}\right.
                \end{array}$ & / &  /    \\  \\
            $\mathcal{L}_{q}^{1+4}$  & $\left[\begin{array}{ccccc}
                 3  &  -1  &  -1  &  -1  &  0  \\
                -1  &   2  &  -1  &   0  &  0  \\
                -1  &  -1  &   4  &  -1  & -1  \\
                -1  &   0  &  -1  &   2  &  0  \\
                 0  &   0  &  -1  &   0  &  1
            \end{array}\right]$   & / &  /    \\  \\
            $\beta_{q}^{1+4}$                        & 1               & / & /    \\  \\
            $\text{dim}(\mathcal{L}_{q}^{1+4})$      & 5               & / & /    \\  \\
            $\text{rank}(\mathcal{L}_{q}^{1+4})$     & 4               & / & /    \\  \\
            $\text{nullity}(\mathcal{L}_{q}^{1+4})$  & 1               & / & /  \\  \\
            $\text{Spectra}(\mathcal{L}_{q}^{1+4})$    & $\{0,1,2,4,5\}$   & / & /  \\ \hline
            \end{tabular}
            \\
            \label{table:ex 7}
        \end{table}


        \begin{table}[H]
            \centering
            \setlength\tabcolsep{10pt}
            \captionsetup{margin=0.9cm}
            \caption{$K_1 \to K_6$}
            \begin{tabular}{c|ccc}
            \hline
            $q$  & $q=0$ & $q=1$ & $q=2$ \\ \hline \\
            $\mathcal{B}_{q+1}^{1+5}$  & $\begin{array}{@{}r@{}c@{}c@{}c@{}c@{}c@{}c@{}c@{}l@{}}
                    & 01 & 12 & 23 & 03 & 24 & 02 & 13 \\
                   \left.\begin{array}{c}
                    0 \\
                    1 \\
                    2 \\
                    3 \\
                    4
                    \end{array}\right[
                    & \begin{array}{c} -1 \\  1  \\  0  \\  0 \\ 0 \end{array}
                    & \begin{array}{c}  0 \\ -1  \\  1  \\  0 \\ 0 \end{array}
                    & \begin{array}{c}  0 \\  0  \\ -1  \\  1 \\ 0 \end{array}
                    & \begin{array}{c} -1 \\  0  \\  0  \\  1 \\ 0 \end{array}
                    & \begin{array}{c}  0 \\  0  \\ -1  \\  0 \\ 1 \end{array}
                    & \begin{array}{c} -1 \\  0  \\  1  \\  0 \\ 0 \end{array}
                    & \begin{array}{c}  0 \\ -1  \\  0  \\  1 \\ 0 \end{array}
                    & \left]\begin{array}{c} \\ \\  \\ \\  \\
                    \end{array}\right.
                \end{array}$ & / &  /    \\  \\
            $\mathcal{B}_{q}^{1}$  & $\begin{array}{@{}r@{}c@{}c@{}c@{}c@{}c@{}l@{}}
                    & 0 & 1 & 2 & 3 & 4   \\
                   \left.\begin{array}{c}

                    \end{array}\right[
                    & \begin{array}{c}  0  \end{array}
                    & \begin{array}{c}  0  \end{array}
                    & \begin{array}{c}  0  \end{array}
                    & \begin{array}{c}  0  \end{array}
                    & \begin{array}{c}  0  \end{array}
                    & \left]\begin{array}{c}
                    \end{array}\right.
                \end{array}$ & / &  /    \\  \\
            $\mathcal{L}_{q}^{1+5}$  & $\left[\begin{array}{ccccc}
                 3  &  -1  &  -1  &  -1  &  0  \\
                -1  &   3  &  -1  &  -1  &  0  \\
                -1  &  -1  &   4  &  -1  &  -1  \\
                -1  &  -1  &  -1  &   3  &  0  \\
                 0  &   0  &  -1  &   0  &  1
            \end{array}\right]$   & / &  /    \\  \\
            $\beta_{q}^{1+5}$                        & 1               & / & /    \\  \\
            $\text{dim}(\mathcal{L}_{q}^{1+5})$      & 5               & / & /    \\  \\
            $\text{rank}(\mathcal{L}_{q}^{1+5})$     & 4               & / & /    \\  \\
            $\text{nullity}(\mathcal{L}_{q}^{1+5})$  & 1               & / & /  \\  \\
            $\text{Spectra}(\mathcal{L}_{q}^{1+5})$    & $\{0,1,4,4,5\}$   & / & /  \\ \hline
            \end{tabular}
            \\
            \label{table:ex 8}
        \end{table}


        \begin{table}[H]
            \centering
            \setlength\tabcolsep{8pt}
            \captionsetup{margin=0.9cm}
            \caption{$K_2 \to K_3$}
            \begin{tabular}{c|ccc}
            \hline
            $q$  & $q=0$ & $q=1$ & $q=2$ \\ \hline \\
            $\mathcal{B}_{q+1}^{2+1}$  & $\begin{array}{@{}r@{}c@{}c@{}c@{}c@{}l@{}}
                    & 01 & 12 & 23 & 03   \\
                   \left.\begin{array}{c}
                    0 \\
                    1 \\
                    2 \\
                    3 \\
                    4
                    \end{array}\right[
                    & \begin{array}{c} -1 \\  1  \\  0  \\  0 \\ 0 \end{array}
                    & \begin{array}{c}  0 \\ -1  \\  1  \\  0 \\ 0 \end{array}
                    & \begin{array}{c}  0 \\  0  \\ -1  \\  1 \\ 0 \end{array}
                    & \begin{array}{c} -1 \\  0  \\  0  \\  1 \\ 0 \end{array}
                    & \left]\begin{array}{c} \\ \\  \\ \\ \\
                    \end{array}\right.
                \end{array}$ & / &  /    \\  \\
            $\mathcal{B}_{q}^{2}$  & $\begin{array}{@{}r@{}c@{}c@{}c@{}c@{}c@{}l@{}}
                    & 0 & 1 & 2 & 3 & 4   \\
                   \left.\begin{array}{c}

                    \end{array}\right[
                    & \begin{array}{c}  0  \end{array}
                    & \begin{array}{c}  0  \end{array}
                    & \begin{array}{c}  0  \end{array}
                    & \begin{array}{c}  0  \end{array}
                    & \begin{array}{c}  0  \end{array}
                    & \left]\begin{array}{c}
                    \end{array}\right.
                \end{array}$ & $\begin{array}{@{}r@{}c@{}c@{}c@{}c@{}l@{}}
                    & 01   \\
                   \left.\begin{array}{c}
                    0 \\
                    1 \\
                    2 \\
                    3 \\
                    4
                    \end{array}\right[
                    & \begin{array}{c} -1 \\  1  \\  0  \\  0 \\ 0 \end{array}
                    & \left]\begin{array}{c} \\ \\  \\ \\ \\
                    \end{array}\right.
                \end{array}$ &  /    \\  \\
            $\mathcal{L}_{q}^{2+1}$  & $\left[\begin{array}{ccccc}
                 2  &  -1  &   0  &  -1  &  0  \\
                -1  &   2  &  -1  &   0  &  0  \\
                 0  &  -1  &   2  &  -1  &  0  \\
                -1  &   0  &  -1  &   2  &  0  \\
                 0  &   0  &   0  &   0  &  0
            \end{array}\right]$  & [2] &  /    \\  \\
            $\beta_{q}^{2+1}$                        & 2               & 0 & /    \\  \\
            $\text{dim}(\mathcal{L}_{q}^{2+1})$      & 5               & 1 & /    \\  \\
            $\text{rank}(\mathcal{L}_{q}^{2+1})$     & 3               & 1 & /    \\  \\
            $\text{nullity}(\mathcal{L}_{q}^{2+1})$  & 2               & 0 & /  \\  \\
            $\text{Spectra}(\mathcal{L}_{q}^{2+1})$    & $\{0,0,2,2,4\}$   & ${2}$ & /  \\ \hline
            \end{tabular}
            \\
            \label{table:ex 9}
        \end{table}


        \begin{table}[H]
            \centering
            \setlength\tabcolsep{8pt}
            \captionsetup{margin=0.9cm}
            \caption{$K_2 \to K_4$}
            \begin{tabular}{c|ccc}
            \hline
            $q$  & $q=0$ & $q=1$ & $q=2$ \\ \hline \\
            $\mathcal{B}_{q+1}^{2+2}$  & $\begin{array}{@{}r@{}c@{}c@{}c@{}c@{}c@{}l@{}}
                    & 01 & 12 & 23 & 03 & 24   \\
                   \left.\begin{array}{c}
                    0 \\
                    1 \\
                    2 \\
                    3 \\
                    4
                    \end{array}\right[
                    & \begin{array}{c} -1 \\  1  \\  0  \\  0 \\ 0 \end{array}
                    & \begin{array}{c}  0 \\ -1  \\  1  \\  0 \\ 0 \end{array}
                    & \begin{array}{c}  0 \\  0  \\ -1  \\  1 \\ 0 \end{array}
                    & \begin{array}{c} -1 \\  0  \\  0  \\  1 \\ 0 \end{array}
                    & \begin{array}{c}  0 \\  0  \\ -1  \\  0 \\ 1 \end{array}
                    & \left]\begin{array}{c} \\ \\  \\ \\ \\
                    \end{array}\right.
                \end{array}$ & / &  /    \\  \\
            $\mathcal{B}_{q}^{2}$  & $\begin{array}{@{}r@{}c@{}c@{}c@{}c@{}c@{}l@{}}
                    & 0 & 1 & 2 & 3 & 4   \\
                   \left.\begin{array}{c}

                    \end{array}\right[
                    & \begin{array}{c}  0  \end{array}
                    & \begin{array}{c}  0  \end{array}
                    & \begin{array}{c}  0  \end{array}
                    & \begin{array}{c}  0  \end{array}
                    & \begin{array}{c}  0  \end{array}
                    & \left]\begin{array}{c}
                    \end{array}\right.
                \end{array}$  & $\begin{array}{@{}r@{}c@{}c@{}c@{}c@{}l@{}}
                    & 01   \\
                   \left.\begin{array}{c}
                    0 \\
                    1 \\
                    2 \\
                    3 \\
                    4
                    \end{array}\right[
                    & \begin{array}{c} -1 \\  1  \\  0  \\  0 \\ 0 \end{array}
                    & \left]\begin{array}{c} \\ \\  \\ \\ \\
                    \end{array}\right.
                \end{array}$ &  /    \\  \\
            $\mathcal{L}_{q}^{2+2}$  & $\left[\begin{array}{ccccc}
                 2  &  -1  &   0  &  -1  &  0  \\
                -1  &   2  &  -1  &   0  &  0  \\
                 0  &  -1  &   3  &  -1  & -1  \\
                -1  &   0  &  -1  &   2  &  0  \\
                 0  &   0  &  -1  &   0  &  1
            \end{array}\right]$  & [2] &  /    \\  \\
            $\beta_{q}^{2+2}$                        & 1               & 0 & /    \\  \\
            $\text{dim}(\mathcal{L}_{q}^{2+2})$      & 5               & 1 & /    \\  \\
            $\text{rank}(\mathcal{L}_{q}^{2+2})$     & 4               & 1 & /    \\  \\
            $\text{nullity}(\mathcal{L}_{q}^{2+2})$  & 1               & 0 & /  \\  \\
            $\text{Spectra}(\mathcal{L}_{q}^{2+2})$    & $\{0,0.8299,2,2.6889,4.4812\}$   & ${2}$ & /  \\ \hline
            \end{tabular}
            \\
            \label{table:ex 10}
        \end{table}


        \begin{table}[H]
            \centering
            \setlength\tabcolsep{8pt}
            \captionsetup{margin=0.9cm}
            \caption{$K_2 \to K_5$}
            \begin{tabular}{c|ccc}
            \hline
            $q$  & $q=0$ & $q=1$ & $q=2$ \\ \hline \\
            $\mathcal{B}_{q+1}^{2+3}$  & $\begin{array}{@{}r@{}c@{}c@{}c@{}c@{}c@{}c@{}l@{}}
                    & 01 & 12 & 23 & 03 & 24 & 02  \\
                   \left.\begin{array}{c}
                    0 \\
                    1 \\
                    2 \\
                    3 \\
                    4
                    \end{array}\right[
                    & \begin{array}{c} -1 \\  1  \\  0  \\  0 \\ 0 \end{array}
                    & \begin{array}{c}  0 \\ -1  \\  1  \\  0 \\ 0 \end{array}
                    & \begin{array}{c}  0 \\  0  \\ -1  \\  1 \\ 0 \end{array}
                    & \begin{array}{c} -1 \\  0  \\  0  \\  1 \\ 0 \end{array}
                    & \begin{array}{c}  0 \\  0  \\ -1  \\  0 \\ 1 \end{array}
                    & \begin{array}{c} -1 \\  0  \\  1  \\  0 \\ 0 \end{array}
                    & \left]\begin{array}{c} \\ \\  \\ \\ \\
                    \end{array}\right.
                \end{array}$ & $\begin{array}{@{}r@{}c@{}c@{}c@{}c@{}c@{}l@{}}
                    & 012 & & 023  \\
                   \left.\begin{array}{c}
                    01
                    \end{array}\right[
                    & \begin{array}{c}  1  \end{array}
                    & \begin{array}{c}  \   \end{array}
                    & \begin{array}{c}  0  \end{array}
                    & \left]\begin{array}{c} \\
                    \end{array}\right.
                \end{array}$ &  /    \\  \\
            $\mathcal{B}_{q}^{2}$  & $\begin{array}{@{}r@{}c@{}c@{}c@{}c@{}c@{}l@{}}
                    & 0 & 1 & 2 & 3 & 4   \\
                   \left.\begin{array}{c}

                    \end{array}\right[
                    & \begin{array}{c}  0  \end{array}
                    & \begin{array}{c}  0  \end{array}
                    & \begin{array}{c}  0  \end{array}
                    & \begin{array}{c}  0  \end{array}
                    & \begin{array}{c}  0  \end{array}
                    & \left]\begin{array}{c}
                    \end{array}\right.
                \end{array}$  & $\begin{array}{@{}r@{}c@{}c@{}c@{}c@{}l@{}}
                    & 01   \\
                   \left.\begin{array}{c}
                    0 \\
                    1 \\
                    2 \\
                    3 \\
                    4
                    \end{array}\right[
                    & \begin{array}{c} -1 \\  1  \\  0  \\  0 \\ 0 \end{array}
                    & \left]\begin{array}{c} \\ \\  \\ \\ \\
                    \end{array}\right.
                \end{array}$ &  /    \\  \\
            $\mathcal{L}_{q}^{2+3}$  & $\left[\begin{array}{ccccc}
                 3  &  -1  &  -1  &  -1  &  0  \\
                -1  &   3  &  -1  &  -1  &  0  \\
                -1  &  -1  &   4  &  -1  & -1  \\
                -1  &  -1  &  -1  &   3  &  0  \\
                 0  &   0  &  -1  &   0  &  1
            \end{array}\right]$  & [3] &  /    \\  \\
            $\beta_{q}^{2+3}$                        & 1               & 0 & /    \\  \\
            $\text{dim}(\mathcal{L}_{q}^{2+3})$      & 5               & 1 & /    \\  \\
            $\text{rank}(\mathcal{L}_{q}^{2+3})$     & 4               & 1 & /    \\  \\
            $\text{nullity}(\mathcal{L}_{q}^{2+3})$  & 1               & 0 & /  \\  \\
            $\text{Spectra}(\mathcal{L}_{q}^{2+3})$    & $\{0,1,2,4,5\}$   & ${3}$ & /  \\ \hline
            \end{tabular}
            \\
            \label{table:ex 11}
        \end{table}


        \begin{table}[H]
            \centering
            \setlength\tabcolsep{8pt}
            \captionsetup{margin=0.9cm}
            \caption{$K_2 \to K_6$}
            \begin{tabular}{c|ccc}
            \hline
            $q$  & $q=0$ & $q=1$ & $q=2$ \\ \hline \\
            $\mathcal{B}_{q+1}^{2+4}$  & $\begin{array}{@{}r@{}c@{}c@{}c@{}c@{}c@{}c@{}c@{}l@{}}
                    & 01 & 12 & 23 & 03 & 24 & 02 & 13 \\
                   \left.\begin{array}{c}
                    0 \\
                    1 \\
                    2 \\
                    3 \\
                    4
                    \end{array}\right[
                    & \begin{array}{c} -1 \\  1  \\  0  \\  0 \\ 0 \end{array}
                    & \begin{array}{c}  0 \\ -1  \\  1  \\  0 \\ 0 \end{array}
                    & \begin{array}{c}  0 \\  0  \\ -1  \\  1 \\ 0 \end{array}
                    & \begin{array}{c} -1 \\  0  \\  0  \\  1 \\ 0 \end{array}
                    & \begin{array}{c}  0 \\  0  \\ -1  \\  0 \\ 1 \end{array}
                    & \begin{array}{c} -1 \\  0  \\  1  \\  0 \\ 0 \end{array}
                    & \begin{array}{c}  0 \\ -1  \\  0  \\  1 \\ 0 \end{array}
                    & \left]\begin{array}{c} \\ \\  \\ \\  \\
                    \end{array}\right.
                \end{array}$ & $\begin{array}{@{}r@{}c@{}c@{}c@{}c@{}c@{}c@{}c@{}l@{}}
                    & 012 & & 023 & & 013 & & 123  \\
                   \left.\begin{array}{c}
                    01
                    \end{array}\right[
                    & \begin{array}{c}  1  \end{array}
                    & \begin{array}{c}  \   \end{array}
                    & \begin{array}{c}  0  \end{array}
                    & \begin{array}{c}  \   \end{array}
                    & \begin{array}{c}  1  \end{array}
                    & \begin{array}{c}  \   \end{array}
                    & \begin{array}{c}  0  \end{array}
                    & \left]\begin{array}{c} \\
                    \end{array}\right.
                \end{array}$ &  /    \\  \\
            $\mathcal{B}_{q}^{2}$  & $\begin{array}{@{}r@{}c@{}c@{}c@{}c@{}c@{}l@{}}
                    & 0 & 1 & 2 & 3 & 4   \\
                   \left.\begin{array}{c}

                    \end{array}\right[
                    & \begin{array}{c}  0  \end{array}
                    & \begin{array}{c}  0  \end{array}
                    & \begin{array}{c}  0  \end{array}
                    & \begin{array}{c}  0  \end{array}
                    & \begin{array}{c}  0  \end{array}
                    & \left]\begin{array}{c}
                    \end{array}\right.
                \end{array}$ & $\begin{array}{@{}r@{}c@{}c@{}c@{}c@{}l@{}}
                    & 01   \\
                   \left.\begin{array}{c}
                    0 \\
                    1 \\
                    2 \\
                    3 \\
                    4
                    \end{array}\right[
                    & \begin{array}{c} -1 \\  1  \\  0  \\  0 \\ 0 \end{array}
                    & \left]\begin{array}{c} \\ \\  \\ \\ \\
                    \end{array}\right.
                \end{array}$ &  /    \\  \\
            $\mathcal{L}_{q}^{2+4}$  & $\left[\begin{array}{ccccc}
                 3  &  -1  &  -1  &  -1  &  0  \\
                -1  &   2  &  -1  &   0  &  0  \\
                -1  &  -1  &   4  &  -1  &  -1  \\
                -1  &   0  &  -1  &   2  &  0  \\
                 0  &   0  &  -1  &   0  &  1
            \end{array}\right]$  & [4] &  /    \\  \\
            $\beta_{q}^{2+4}$                        & $1$               & $0$ & /    \\  \\
            $\text{dim}(\mathcal{L}_{q}^{2+4})$      & $5$               & $1$ & /    \\  \\
            $\text{rank}(\mathcal{L}_{q}^{2+4})$     & $4$               & $1$ & /    \\  \\
            $\text{nullity}(\mathcal{L}_{q}^{2+4})$  & $1$               & $0$ & /  \\  \\
            $\text{Spectra}(\mathcal{L}_{q}^{2+4})$    & $\{0,1,4,4,5\}$     & ${4}$ & /  \\ \hline
            \end{tabular}
            \\
            \label{table:ex 12}
        \end{table}


        \begin{table}[H]
            \centering
            \setlength\tabcolsep{8pt}
            \captionsetup{margin=0.9cm}
            \caption{$K_3 \to K_5$}
            \begin{tabular}{c|ccc}
            \hline
            $q$  & $q=0$ & $q=1$ & $q=2$ \\ \hline \\
            $\mathcal{B}_{q+1}^{3+2}$  & $\begin{array}{@{}r@{}c@{}c@{}c@{}c@{}c@{}c@{}l@{}}
                    & 01 & 12 & 23 & 03 & 24 & 02  \\
                   \left.\begin{array}{c}
                    0 \\
                    1 \\
                    2 \\
                    3 \\
                    4
                    \end{array}\right[
                    & \begin{array}{c} -1 \\  1  \\  0  \\  0 \\ 0 \end{array}
                    & \begin{array}{c}  0 \\ -1  \\  1  \\  0 \\ 0 \end{array}
                    & \begin{array}{c}  0 \\  0  \\ -1  \\  1 \\ 0 \end{array}
                    & \begin{array}{c} -1 \\  0  \\  0  \\  1 \\ 0 \end{array}
                    & \begin{array}{c}  0 \\  0  \\ -1  \\  0 \\ 1 \end{array}
                    & \begin{array}{c} -1 \\  0  \\  1  \\  0 \\ 0 \end{array}
                    & \left]\begin{array}{c} \\ \\  \\ \\ \\
                    \end{array}\right.
                \end{array}$ & $\begin{array}{@{}r@{}c@{}c@{}l@{}}
                    & 012 & 023  \\
                   \left.\begin{array}{c}
                    01 \\
                    12 \\
                    23 \\
                    03
                    \end{array}\right[
                    & \begin{array}{c}  1 \\  1  \\  0  \\  0   \end{array}
                    & \begin{array}{c}  0 \\  0  \\  1  \\ -1   \end{array}
                    & \left]\begin{array}{c} \\ \\  \\ \\
                    \end{array}\right.
                \end{array}$ &  /    \\  \\
            $\mathcal{B}_{q}^{3}$  & $\begin{array}{@{}r@{}c@{}c@{}c@{}c@{}c@{}l@{}}
                    & 0 & 1 & 2 & 3 & 4   \\
                   \left.\begin{array}{c}

                    \end{array}\right[
                    & \begin{array}{c}  0  \end{array}
                    & \begin{array}{c}  0  \end{array}
                    & \begin{array}{c}  0  \end{array}
                    & \begin{array}{c}  0  \end{array}
                    & \begin{array}{c}  0  \end{array}
                    & \left]\begin{array}{c}
                    \end{array}\right.
                \end{array}$ & $\begin{array}{@{}r@{}c@{}c@{}c@{}c@{}l@{}}
                    & 01 & 12 & 23 & 03  \\
                   \left.\begin{array}{c}
                    0 \\
                    1 \\
                    2 \\
                    3 \\
                    4
                    \end{array}\right[
                    & \begin{array}{c} -1 \\  1  \\  0  \\  0 \\ 0  \end{array}
                    & \begin{array}{c}  0 \\ -1  \\  1  \\  0 \\ 0  \end{array}
                    & \begin{array}{c}  0 \\  0  \\ -1  \\  1 \\ 0  \end{array}
                    & \begin{array}{c} -1 \\  0  \\  0  \\  1 \\ 0  \end{array}
                    & \left]\begin{array}{c} \\ \\  \\ \\ \\
                    \end{array}\right.
                \end{array}$ &  /    \\  \\
            $\mathcal{L}_{q}^{3+2}$  & $\left[\begin{array}{ccccc}
                 3  &  -1  &  -1  &  -1  &  0  \\
                -1  &   2  &  -1  &   0  &  0  \\
                -1  &  -1  &   4  &  -1  & -1  \\
                -1  &   0  &  -1  &   2  &  0  \\
                 0  &   0  &  -1  &   0  &  1
            \end{array}\right]$  & $\left[\begin{array}{cccc}
                 3  &   0  &   0  &   1   \\
                 0  &   3  &  -1  &   0   \\
                 0  &  -1  &   3  &   0   \\
                 1  &   0  &   0  &   3
            \end{array}\right]$ &  /    \\  \\
            $\beta_{q}^{3+2}$                        & 1               & 0 & /    \\  \\
            $\text{dim}(\mathcal{L}_{q}^{3+2})$      & 5               & 4 & /    \\  \\
            $\text{rank}(\mathcal{L}_{q}^{3+2})$     & 4               & 4 & /    \\  \\
            $\text{nullity}(\mathcal{L}_{q}^{3+2})$  & 1               & 0 & /  \\  \\
            $\text{Spectra}(\mathcal{L}_{q}^{3+2})$    & $\{0,1,2,4,5\}$   & $\{2,2,4,4\}$ & /  \\ \hline
            \end{tabular}
            \\
            \label{table:ex 13}
        \end{table}


        \begin{table}[H]
            \centering
            \setlength\tabcolsep{8pt}
            \captionsetup{margin=0.9cm}
            \caption{$K_3 \to K_6$}
            \begin{tabular}{c|ccc}
            \hline
            $q$  & $q=0$ & $q=1$ & $q=2$ \\ \hline \\
            $\mathcal{B}_{q+1}^{3+3}$  & $\begin{array}{@{}r@{}c@{}c@{}c@{}c@{}c@{}c@{}c@{}l@{}}
                    & 01 & 12 & 23 & 03 & 24 & 02 & 13 \\
                   \left.\begin{array}{c}
                    0 \\
                    1 \\
                    2 \\
                    3 \\
                    4
                    \end{array}\right[
                    & \begin{array}{c} -1 \\  1  \\  0  \\  0 \\ 0 \end{array}
                    & \begin{array}{c}  0 \\ -1  \\  1  \\  0 \\ 0 \end{array}
                    & \begin{array}{c}  0 \\  0  \\ -1  \\  1 \\ 0 \end{array}
                    & \begin{array}{c} -1 \\  0  \\  0  \\  1 \\ 0 \end{array}
                    & \begin{array}{c}  0 \\  0  \\ -1  \\  0 \\ 1 \end{array}
                    & \begin{array}{c} -1 \\  0  \\  1  \\  0 \\ 0 \end{array}
                    & \begin{array}{c}  0 \\ -1  \\  0  \\  1 \\ 0 \end{array}
                    & \left]\begin{array}{c} \\ \\  \\ \\  \\
                    \end{array}\right.
                \end{array}$ & $\begin{array}{@{}r@{}c@{}c@{}c@{}c@{}l@{}}
                    & 012 & 023 & 013 & 123  \\
                   \left.\begin{array}{c}
                    01 \\
                    12 \\
                    23 \\
                    03
                    \end{array}\right[
                    & \begin{array}{c}  1 \\  1  \\  0  \\  0   \end{array}
                    & \begin{array}{c}  0 \\  0  \\  1  \\ -1   \end{array}
                    & \begin{array}{c}  1 \\  0  \\  0  \\ -1   \end{array}
                    & \begin{array}{c}  0 \\  1  \\  1  \\  0   \end{array}
                    & \left]\begin{array}{c} \\ \\  \\ \\
                    \end{array}\right.
                \end{array}$ &  /    \\  \\
            $\mathcal{B}_{q}^{3}$  &  $\begin{array}{@{}r@{}c@{}c@{}c@{}c@{}c@{}l@{}}
                    & 0 & 1 & 2 & 3 & 4   \\
                   \left.\begin{array}{c}

                    \end{array}\right[
                    & \begin{array}{c}  0  \end{array}
                    & \begin{array}{c}  0  \end{array}
                    & \begin{array}{c}  0  \end{array}
                    & \begin{array}{c}  0  \end{array}
                    & \begin{array}{c}  0  \end{array}
                    & \left]\begin{array}{c}
                    \end{array}\right.
                \end{array}$ & $\begin{array}{@{}r@{}c@{}c@{}c@{}c@{}l@{}}
                    & 01 & 12 & 23 & 03  \\
                   \left.\begin{array}{c}
                    0 \\
                    1 \\
                    2 \\
                    3 \\
                    4
                    \end{array}\right[
                    & \begin{array}{c} -1 \\  1  \\  0  \\  0 \\ 0  \end{array}
                    & \begin{array}{c}  0 \\ -1  \\  1  \\  0 \\ 0  \end{array}
                    & \begin{array}{c}  0 \\  0  \\ -1  \\  1 \\ 0  \end{array}
                    & \begin{array}{c} -1 \\  0  \\  0  \\  1 \\ 0  \end{array}
                    & \left]\begin{array}{c} \\ \\  \\ \\ \\
                    \end{array}\right.
                \end{array}$ &  /    \\  \\
            $\mathcal{L}_{q}^{3+3}$  & $\left[\begin{array}{ccccc}
                 3  &  -1  &  -1  &  -1  &  0  \\
                -1  &   3  &  -1  &  -1  &  0  \\
                -1  &  -1  &   4  &  -1  &  -1  \\
                -1  &  -1  &  -1  &   3  &  0  \\
                 0  &   0  &  -1  &   0  &  1
            \end{array}\right]$  & $\left[\begin{array}{cccc}
                 4  &   0  &   0  &   0   \\
                 0  &   4  &   0  &   0   \\
                 0  &   0  &   4  &   0   \\
                 0  &   0  &   0  &   4
            \end{array}\right]$ &  /    \\  \\
            $\beta_{q}^{3+3}$                        & 1               & 0 & /    \\  \\
            $\text{dim}(\mathcal{L}_{q}^{3+3})$      & 5               & 4 & /    \\  \\
            $\text{rank}(\mathcal{L}_{q}^{3+3})$     & 4               & 4 & /    \\  \\
            $\text{nullity}(\mathcal{L}_{q}^{3+3})$  & 1               & 0 & /  \\  \\
            $\text{Spectra}(\mathcal{L}_{q}^{3+3})$    & $\{0,1,4,4,5\}$   & $\{4,4,4,4\}$ & /  \\ \hline
            \end{tabular}
            \\
            \label{table:ex 14}
        \end{table}


        \begin{table}[H]
            \centering
            \setlength\tabcolsep{8pt}
            \captionsetup{margin=0.9cm}
            \caption{$K_4 \to K_6$}
            \begin{tabular}{c|ccc}
            \hline
            $q$  & $q=0$ & $q=1$ & $q=2$ \\ \hline \\
            $\mathcal{B}_{q+1}^{4+2}$  & $\begin{array}{@{}r@{}c@{}c@{}c@{}c@{}c@{}c@{}c@{}l@{}}
                    & 01 & 12 & 23 & 03 & 24 & 02 & 13 \\
                   \left.\begin{array}{c}
                    0 \\
                    1 \\
                    2 \\
                    3 \\
                    4
                    \end{array}\right[
                    & \begin{array}{c} -1 \\  1  \\  0  \\  0 \\ 0 \end{array}
                    & \begin{array}{c}  0 \\ -1  \\  1  \\  0 \\ 0 \end{array}
                    & \begin{array}{c}  0 \\  0  \\ -1  \\  1 \\ 0 \end{array}
                    & \begin{array}{c} -1 \\  0  \\  0  \\  1 \\ 0 \end{array}
                    & \begin{array}{c}  0 \\  0  \\ -1  \\  0 \\ 1 \end{array}
                    & \begin{array}{c} -1 \\  0  \\  1  \\  0 \\ 0 \end{array}
                    & \begin{array}{c}  0 \\ -1  \\  0  \\  1 \\ 0 \end{array}
                    & \left]\begin{array}{c} \\ \\  \\ \\  \\
                    \end{array}\right.
                \end{array}$ & $\begin{array}{@{}r@{}c@{}c@{}c@{}c@{}l@{}}
                    & 012 & 023 & 013 & 123  \\
                   \left.\begin{array}{c}
                    01 \\
                    12 \\
                    23 \\
                    03 \\
                    24
                    \end{array}\right[
                    & \begin{array}{c}  1 \\  1  \\  0  \\  0  \\ 0  \end{array}
                    & \begin{array}{c}  0 \\  0  \\  1  \\ -1  \\ 0   \end{array}
                    & \begin{array}{c}  1 \\  0  \\  0  \\ -1  \\ 0   \end{array}
                    & \begin{array}{c}  0 \\  1  \\  1  \\  0  \\ 0   \end{array}
                    & \left]\begin{array}{c} \\ \\  \\ \\ \\
                    \end{array}\right.
                \end{array}$ &  /    \\  \\
            $\mathcal{B}_{q}^{4}$  &  $\begin{array}{@{}r@{}c@{}c@{}c@{}c@{}c@{}l@{}}
                    & 0 & 1 & 2 & 3 & 4   \\
                   \left.\begin{array}{c}

                    \end{array}\right[
                    & \begin{array}{c}  0  \end{array}
                    & \begin{array}{c}  0  \end{array}
                    & \begin{array}{c}  0  \end{array}
                    & \begin{array}{c}  0  \end{array}
                    & \begin{array}{c}  0  \end{array}
                    & \left]\begin{array}{c}
                    \end{array}\right.
                \end{array}$ & $\begin{array}{@{}r@{}c@{}c@{}c@{}c@{}c@{}l@{}}
                    & 01 & 12 & 23 & 03 & 24 \\
                   \left.\begin{array}{c}
                    0 \\
                    1 \\
                    2 \\
                    3 \\
                    4
                    \end{array}\right[
                    & \begin{array}{c} -1 \\  1  \\  0  \\  0 \\ 0  \end{array}
                    & \begin{array}{c}  0 \\ -1  \\  1  \\  0 \\ 0  \end{array}
                    & \begin{array}{c}  0 \\  0  \\ -1  \\  1 \\ 0  \end{array}
                    & \begin{array}{c} -1 \\  0  \\  0  \\  1 \\ 0  \end{array}
                    & \begin{array}{c}  0 \\  0  \\ -1  \\  0 \\ 1  \end{array}
                    & \left]\begin{array}{c} \\ \\  \\ \\ \\
                    \end{array}\right.
                \end{array}$ &  /    \\  \\
            $\mathcal{L}_{q}^{4+2}$  & $\left[\begin{array}{ccccc}
                 3  &  -1  &  -1  &  -1  &  0  \\
                -1  &   3  &  -1  &  -1  &  0  \\
                -1  &  -1  &   4  &  -1  &  -1  \\
                -1  &  -1  &  -1  &   3  &  0  \\
                 0  &   0  &  -1  &   0  &  1
            \end{array}\right]$  & $\left[\begin{array}{ccccc}
                 4  &   0  &   0  &   0  &   0  \\
                 0  &   4  &   0  &   0  &  -1  \\
                 0  &   0  &   4  &   0  &   1  \\
                 0  &   0  &   0  &   4  &   0  \\
                 0  &  -1  &   1  &   0  &   2
            \end{array}\right]$ &  /    \\  \\
            $\beta_{q}^{4+2}$                        & 1               & 0 & /    \\  \\
            $\text{dim}(\mathcal{L}_{q}^{4+2})$      & 5               & 5 & /    \\  \\
            $\text{rank}(\mathcal{L}_{q}^{4+2})$     & 4               & 5 & /    \\  \\
            $\text{nullity}(\mathcal{L}_{q}^{4+2})$  & 1               & 0 & /  \\  \\
            $\text{Spectra}(\mathcal{L}_{q}^{4+2})$    & $\{0,1,4,4,5\}$   & $\{1.2679,4,4,4,4.7321\}$ & /  \\ \hline
            \end{tabular}
            \\
            \label{table:ex 15}
        \end{table}




\section{Parameters in the protein B-factor prediction}\label{app:application}
    \begin{table}[H]
            \centering
            \setlength\tabcolsep{5pt}
            \captionsetup{margin=0.9cm}
            \caption{Fitting parameters from $w_0$ to $w_5$.  }
            \begin{tabular}{ccccccc}
            \hline
             $r$          & $0$         & $1$       & $2$        & $3$      & $4$         & $5$          \\ \hline
             $w_r$        & $10.6102$    & $0.2026$   & $-0.0031$   & $0.2169$   & $0.3127$   & $0.2815$     \\ \hline
            \end{tabular}
            \label{table:weights1}
        \end{table}

    \begin{table}[H]
            \centering
            \setlength\tabcolsep{5pt}
            \captionsetup{margin=0.9cm}
            \caption{Fitting parameters from $w_6$ to $w_{11}$. }
            \begin{tabular}{ccccccc}
            \hline
             $r$          & $6$          & $7$         & $8$         & $9$          & $10$         & $11$          \\ \hline
             $w_r$        & $-0.4623$    & $1.0203$    & $0.6110$    & $-0.6872$    & $-1.0695$    & $4.4257$     \\ \hline
            \end{tabular}
            \label{table:weights2}
        \end{table}

\end{appendices}

\end{document}